\documentclass[leqno]{article}
\usepackage{latexsym,amsmath,amsfonts,mathrsfs,amssymb,amsthm}
\usepackage{yhmath}
\usepackage[usenames]{color}
\usepackage{graphicx}
\def\R{\mathbb R}
\def\N{\mathbb N}

\def\Z{\mathbb Z}
\def\a{\alpha}
\def\e{\epsilon}
\def\d{\delta}
\def\pa{\partial}
\def\lg{\langle}
\def\rg{\rangle}
\def\Y{\mathbb Y}
\def\T{\mathbb T}
\def\P{\mathbb P}
\def\H{{\cal H} }

\def\mA{\mathfrak A}
\def\mc{\mathfrak c}

\def\mL{\mathfrak L}
\def\mS{\mathfrak S}
\def\ms{\mathfrak s}
\def\mZ{\mathfrak Z}
\def\cU{{\cal U}}
\def\cT{{\cal T}}
\def\cS{{\cal S}}
\def\mt{{\mathfrak t}}
\def\mT{{\mathfrak T}}
\def\mX{{\mathfrak X}}
\def\mI{\mathfrak I}
\def\mK{\mathfrak K}
\def\cG{{\cal G}}
\def\ml{\mathfrak l}

\def\g{\gamma}
\def\G{\Gamma}
\def\s{\sigma}
\def\S{\Sigma}
\def\F{{\cal F} }
\def\oF{\overline{\cal F} }
\def\tb{\textcolor{blue}}

\def\th{\theta}
\def\t{\tilde}
\def\x{\xi}

\def\z{\zeta}

\def\wg{\wedge}
\def\be{\begin{equation}}
\def\ee{\end{equation}}
\def\bs{\backslash}
\def\qed{\hfill$\Box$\bigskip}
\def\nd{\noindent\textbf{Proof.} }
\numberwithin{equation}{section}
\newtheorem{thm}{Theorem}[section]
\newtheorem{lem}[thm]{Lemma}
\newtheorem{pro}[thm]{Proposition}
\newtheorem{defn}[thm]{Definition}
\newtheorem{cor}[thm]{Corollary}
\newtheorem{rem}[thm]{Remark}
\begin{document}
\bigskip

\centerline{\Large \textbf{Sliding stability and uniqueness for the set ${Y}\times {Y}$}}

\bigskip

\centerline{\large Xiangyu Liang}


\begin{center}School of mathematical sciences, Beihang University \\
maizeliang@gmail.com \end{center}

\centerline {\large\textbf{Abstract.}}This article is dedicated to discuss the sliding stability and the uniqueness property for the 2-dimensional minimal cone $Y\times Y$ in $\R^4$. This problem is motivated by the classification of singularities for Almgren minimal sets, a model for Plateau's problem in the setting of sets. Minimal cones are blow up limits of Almgren minimal sets, thus the list of all minimal cones gives all possible types of singularities that can occur for minimal sets.  

As proved in \cite{2T}, when several 2-dimensional Almgren (resp. topological) minimal cones are Almgren (resp. topological) sliding stable, and Almgren (resp. topological) unique, the almost orthogonal union of them stays minimal. Hence if several minimal cones admit sliding stability and uniqueness properties, then we can use their almost orthogonal unions to generate new families of minimal cones. One then naturally ask which minimal cones admit these two properties.

This list of known 2-dimensional minimal cones in arbitrary ambient dimension is not long, and the stability and uniqueness properties for all the known 2-dimensional minimal cones, except for $Y\times Y$, have already been established in the previous works \cite{stablePYT, uniquePYT}. Among all the known 2-dimensional minimal cones,  $Y\times Y$ is the only one whose stability and uniqueness properties were left unsolved. This is due to two main reasons : 1) $Y\times Y$ is the only known minimal cone which is essentially of codimension larger than 1---that is, we cannot decompose it into transversal unions of minimal cones of codimension 1. 2) $Y\times Y$ lives in $\R^4$, where we know very little about which types of singularities can occur in a minimal set. This makes it difficult to control and estimate the measures of all possible competitors. Due to the above two issues, new ideas are required here for solving the problem.

We give affirmative answers to this problem for the stability and uniqueness properties for $Y\times Y$ in this paper: we prove that the set $Y\times Y$ is both Almgren sliding stable, and Almgren unique; for the topological case, we prove its topological sliding stability and topological uniqueness for the coefficient group $\Z_2$. This result, along with the results in \cite{2T, stablePYT, uniquePYT}, allows us to use all the known 2-dimensional minimal cones to generate new 2-dimensional minimal cones by taking almost orthogonal unions.

\bigskip


\bigskip


\section{Introduction}

In this article we prove the Almgren (resp. topological) sliding stability and uniqueness property for the 2-dimensional minimal cone in $Y\times Y$. The very original motivation of these results comes from the classification of singularities for minimal sets.

The notion of minimal sets (in the sense of Almgren \cite{Al76}, Reifenberg \cite{Rei60}. See David \cite{DJT}, Liang \cite{topo}, etc., for other variances) is a way to try to solve Plateau's problem in the setting of sets. Plateau's problem, as one of the main interests in geometric measure theory, aims at understanding the existence, regularity and local structure of physical objects that minimize the area while spanning a given boundary, such as soap films.

Roughly speaking, given integers $d<n$, and a region $U\subset \R^n$, a relatively closed set $E$ in $U$ of dimension $d$ is said to be minimal in $U$, if we have 
\be \H^d(\varphi(E))\ge \H^d(E)\ee
for any Lipschitz deformation $\varphi$ in $U$.

See Definition \tb{2.3} for the precise definition.

It is known (cf. Almgren \cite{Al76}, David \& Semmes \cite{DS00}) that a $d$-dimensional minimal set $E$ admits a unique tangent plane at almost every point $x$. In this case the local structure around a such point is very clear:  the set $E$ is locally a minimal surface (and hence real analytic) around a such point, due to the famous result of Allard \cite{All72}. 

So we are mostly interested in what happens around points that admit no tangent plane, namely, the singular points. 

In \cite{DJT}, David proved that the blow-up limits (''tangent objects'') of $d$-dimensional minimal sets at a point are $d$-dimensional minimal cones (minimal sets that are cones in the means time). Blow-up limits of a set at a point reflect the asymptotic behavior of the set at infinitesimal scales around this point. As a consequence, a first step to study local structures of minimal sets, is to classify all possible types of singularities--that is to say, minimal cones. 

The plan for the list of $d$-dimensional minimal cones in $\R^n$ is very far from clear. Even for $d=2$, we know very little, except for the case in $\R^3$, where Jean Taylor \cite{Ta} gave a complete classification in 1976, and the list is in fact already known a century ago in other circumstances (see \cite{La} and \cite{He}). They are, modulo isomorphism: a plane, a $\Y$ set (the union of three half planes that meet along a straight line where they make angles of $120^\circ$), and a $\T$ set (the cone over the 1-skeleton of a regular tetrahedron centred at the origin). See the pictures below.

\centerline{\includegraphics[width=0.16\textwidth]{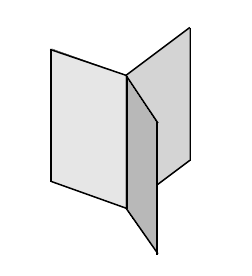} \hskip 2cm\includegraphics[width=0.2\textwidth]{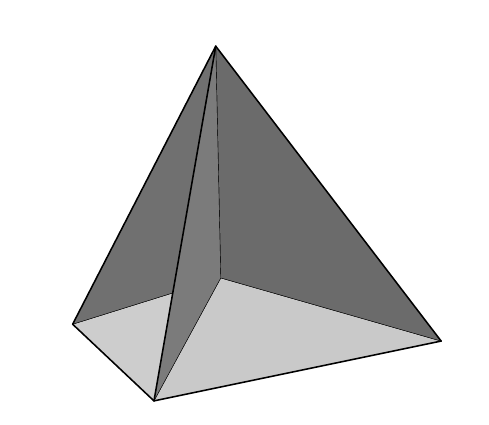}}
\nopagebreak[4]
\hskip 4.5cm a $\Y$ set\hskip 3.9cm  a $\T$ set

Based on the above, a natural way to find new types of singularities, is by taking unions and products of known minimal cones.

For unions: The minimality of the union of two orthogonal minimal sets of dimension $d$ can be obtained easily from a well known geometric lemma (cf. for example Lemma 5.2 of \cite{Mo84}). Thus one suspects that if the angle between two minimal sets is not far from orthogonal, the union of them might also be minimal.

In case of planes, the author proved in \cite{2p} and \cite{2ptopo}, that the almost orthogonal unions of several $d$-dimensional planes are Almgren and topological minimal. When the number of planes is two, this is part of Morgan's conjecture in \cite{Mo93} on the angle condition under which a union of two planes is minimal. 

As for minimal cones other than unions of planes, since they are all with non isolated singularities (after the structure Theorem \tb{2.22}), the situation is much more complicated, as briefly stated in the introduction of \cite{2T}. Up to now we are able to treat a big part of 2 dimensional cases: in \cite{2T} we proved that the almost orthogonal union of several 2-dimensional minimal cones (in any dimension) is minimal, provided that they all admit sliding stability and uniqueness property. See Definitions \tb{3.6, 3.7 and 7.1} for the precise definitions of these two properties. Moreover, this union also satisfies the same sliding stability and uniqueness property. This enables us to continue obtaining infinitely many new families of minimal cones by taking a finite number of iterations of almost orthogonal unions.

As a consequence, the following question arises naturally: which minimal cones admit these two properties?

This list of known 2-dimensional minimal cones in arbitrary ambient dimension is not long, and the stability and uniqueness properties for all the known 2-dimensional minimal cones, except for $Y\times Y$, have already been established in the previous works. See the account in the introduction of \cite{stablePYT}. 

The cone $Y\times Y$ is the product of two 1-dimensional $\Y$ sets. A 1-dimensional $\Y$ set is the union of three half lines issued from the origin and making angles of $120^\circ$. It is a subset of $\R^2$, hence $Y\times Y$ is a 2-dimensional minimal cone in $\R^4$.

Among all the known 2-dimensional minimal cones, $Y\times Y$ is the only one whose stability and uniqueness properties were left unsolved. This is due to the following main reasons: 

1) The cone $Y\times Y$ is the only known 2-dimensional minimal cone which is essentially of codimension larger than 1---that is, we cannot decompose it into transversal unions of 2-dimensional minimal cones of codimension 1. As consequence, 
the old idea for proving the sliding stability for all other 2-dimensional minimal cones does not work for $Y\times Y$. Briefly, in $\R^3$, we are estimating the sum of the $\H^2$ measure of disjoint regions in the unit sphere $S^2$, hence the measures that we calculate are of full dimension; while in $\R^4$, we have to estimate the sum of the $\H^2$ measure of disjoint sets of dimension 2 in the 3-dimensional space $S^3$. 
See the beginning of Section \tb{5} for a detailed description of this difficulty.

2) Another issue is due to the lack of knowledge of local structure for 2-dimensional minimal cones in $\R^4$. Recall that in $\R^3$, we have a complete description for possible local structures for minimal sets. But in $\R^4$, we do not even know which kind of singularities can occur in a minimal set. This makes it difficult to control and estimate the measure of competitors for $Y\times Y$, hence in the proof of uniqueness property for $Y\times Y$, we have to find new ideas to exclude all unknown types of singularities for possible competitors.

In this paper, we manage to give solutions for the above two main difficulties, and give affirmative answers to the problem of the topological and Almgren sliding stability and uniqueness properties (Theorem \tb{6.7, 7.4}): 

--we prove that the set $Y\times Y$ is both Almgren sliding stable, and Almgren unique. This result, along with the results in \cite{2T, stablePYT, uniquePYT}, allows us to use all the known 2-dimensional minimal cones to generate new 2-dimensional minimal cones by taking almost orthogonal unions.;

--for the topological case, we prove its topological sliding stability and topological uniqueness for the coefficient group $\Z_2$ (while all the other known 2-dimensional minimal cone are topological sliding stable and topological unique for all abelian groups). Hence when we take a almost orthogonal union of several 2-dimensional miniabelianmal cones, if one of them is $Y\times Y$, then we only know that this union is topological minimal for the coefficient group $\Z_2$.

\begin{rem}$1^\circ$ The readers, especially geometric analysts, will probably be puzzled about the claim that one of the main results is of the type: some minimizer is stable, which is usually immediate. In fact, here in our circumstance, the minimality is with respect to a fixed boundary, while the stability means that the measure stays minimal even when we allow the boundary of the sets to move. See Section \tb{2.3} of \cite{stablePYT} for some descriptions and examples. 

$2^\circ$ The notion of sliding stability is somehow a quatitative version of that of the ''minimality with sliding boundary'' introduced by David \cite{GD13, GD14}. The model of minimal sets with sliding boundary gives a general frame under which people can study the local structure for minimal sets at the boundary; while the sliding stabilities introduced in this article, as described above, is used to study a different and probably more technical problem, which is part of the study of local structure for minimal sets around interior points. 
\end{rem}

The plan for the rest of the article is the following: 

In Section \tb{2} we introduce basic definitions and preliminaries for minimal sets, and regularity properties for 2-dimensional minimal cones. 

We treat the sliding stabilities for $Y\times Y$ in Sections \tb{3-6}. 

In Section \tb{3} we introduce the associated convex domain $\cU$ and give the definition of various sliding stabilities for general 2-dimensional minimal cones, and then give some corresponding specifications for $Y\times Y$.

We  do necessary simplifications for topological sliding competitors for $Y\times Y$ in $\cU$ in Section \tb{4}. We prove that, to prove the topological sliding stability for $Y\times Y$, it is enough to consider the class $\F(\eta,\d,\nu,L)$ of competitors with some special regularity property.

In Section \tb{5} we decompose any competitor $F\in \F(\eta,\d,\nu,L)$ for $Y\times Y$ into 9 parts $F_{ij}, 1\le i,j\le 3$, and these 9 parts can only meet each other in some prescribed way. Then we apply the product of paired calibrations introduced in \cite{YXY} to these nine parts, and give a lower bound of the measure of $F$ that is a  quantity that depends only on $F\cap \pa \cU$.

In Section \tb{6} we use the regularity property for sets $F\in\F(\eta,\d,\nu,L)$ to calculate the above quantity that depends on $F\cap \pa \cU$. We will prove that this quantity will be uniformly bounded below by the measure of $(Y\times Y)\cap \cU$, and thus complete the proof of the Almgren and $\Z_2$-topological sliding stabilities for $Y\times Y$.

We discuss the Almgren and $\Z_2$-topological uniqueness for $Y\times Y$ in Section \tb{7}.

\textbf{Acknowledgement:} This work is supported by the National Natural Science Foundation of China (Grant Nos. 11871090, 12271018).

\section{Definitions and preliminaries}

\noindent\textbf{Some useful notation}

$[a,b]$ is the line segment with end points $a$ and $b$;

$\overrightarrow{ab}$ is the vector $b-a$;

$R_{a,b}$ denotes the half line issued from the point $a$ and passing through $b$;

$B(x,r)$ is the open ball with radius $r$ and centered on $x$;

$\overline B(x,r)$ is the closed ball with radius $r$ and center $x$; 

$\H^d$ is the Hausdorff measure of dimension $d$ ;

$d_H(E,F)=\max\{\sup\{d(y,F):y\in E\},\sup\{d(y,E):y\in F\}\}$ is the Hausdorff distance between two sets $E$ and $F$. 

For any subset $K\subset \R^n$, the local Hausdorff distance in $K$ $d_K$ between two sets $E,F$ is defined as $d_K(E,F)=\max\{\sup\{d(y,F):y\in E\cap K\},\sup\{d(y,E):y\in F\cap K\}\}$;

For any open subset $U\subset \R^n$, let $\{E_n\}_n$, $F$ be closed sets in $U$, we say that $F$ is the Hausdorff limit of $\{E_n\}_n$, if for any compact subset $K\subset U$, $\lim_n d_K(E_n,F)=0$;

$d_{x,r}$ : the relative distance with respect to the ball $B(x,r)$, is defined by
$$ d_{x,r}(E,F)=\frac 1r\max\{\sup\{d(y,F):y\in E\cap B(x,r)\},\sup\{d(y,E):y\in F\cap B(x,r)\}\}.$$

A polyhedral comples in $\R^n$ is a complex composed of convex polygons. For any polyhedral complex $\cS$ of dimension $d$ in $\R^n$, denote by $|\cS|$ the support $\cup\{\sigma:\sigma\in \cS\}$ of $\cS$. And for any $0\le j\le d$, $\cS^d$ denotes the complex composed of all $d$-faces of $\cS$.

For any (affine) subspace $Q$ of $\R^n$, and $x\in Q$, $r>0$, $B_Q(x,r)$ stands for $B(x,r)\cap Q$, the open ball in $Q$; and $\pi_Q$ stands for the orthogonal projection from $\R^n$ to $Q$. If $\{x_i\}_{1\le i\le m}$ is a family of vectors in $\R^n$, denote by  $Q_{x_1\wg\cdots \wg x_m}$ the linear subspace of $\R^n$ generated by this family, then $B_{x_1\wg\cdots \wg x_m}(x,r)$ stands for $B_{Q_{x_1\wg\cdots \wg x_m}}(x,r)$. Also denote by $\pi_{x_1\wg\cdots \wg x_m}$ the orthogonal projection to $Q_{x_1\wg\cdots \wg x_m}$.

For any subset $E$ of $\R^n$ and any $r>0$, we call $B(E,r):=\{x\in \R^n: dist (x,E)<r\}$ the $r$ neighborhood of $E$.

If $E$ is a $d$-rectifiable set, denote by $T_xE$ the tangent plane (if it exists and is unique) of $E$ at $x$.

If $E$ is a triangulable space (that is, it is homeomorphic to a simplicial complex), for any abelian group $G$, let $H_d(E, G)$ denote the $d$-dimensional simplicial homology group of $E$ with coefficient in $G$ (this group does not depend on the triangulation, since it is always natually isomorphic to the singular homology of $E$ with coefficient in $G$). For a simplicial $d$-chain $\sigma$ in $E$, let $[\s]$ denote its homology class in $H_d(E,G)$, and let $|\s|$ denote its $d$-dimensional support. We write $\sigma\sim \sigma'$ if $[\sigma-\sigma']=0$ in $H_d(E,G)$.

For any piecewise smooth simple curves or surfaces, they also stands for the $\Z_2$-chains naturally associated to them with corresponding dimensions.

Let $P$ be a 2-dimensional plane, and let $\s$ be any closed simplicial $\Z_2$-1-chain in $P$. Since there is no non trivial 3-chain in $P$, we know that there is a unique $\Z_2$-2-chain in $P$ whose boundary is $\s$. For each $\s$, let $\S_\s$ denote the $\Z_2$-2-chain in $P$ so that $\pa \S_\s=\s$. Set $|\s|_2=\H^2(|\S_\s|)$ the 2 dimensional measure of the support $|\S_\s|$ of $\S_\s$. In some sense, $|\s|_2$ is the area that $\s$ encloses.

\subsection{Basic definitions and notations about minimal sets}

In the next definitions, fix integers $0<d<n$. We first give a general definition for minimal sets. Briefly, a minimal set is a closed set which minimizes the Hausdorff measure among a certain class of competitors. Different choices of classes of competitors give different kinds of minimal sets.

\begin{defn}[Minimal sets]Let $0<d<n$ be integers. Let $U\subset \R^n$ be an open set. A relatively closed set $E\subset U$ is said to be minimal of dimension $d$ in $U$ with respect to the competitor class $\mathscr F$ (which contains $E$) if 
\be \H^d(E\cap B)<\infty\mbox{ for every compact ball }B\subset U,\ee
and
\be \H^d(E\bs F)\le \H^d(F\bs E)\ee
for any competitor $F\in\mathscr F$.
\end{defn}

\begin{defn}[Almgren competitor] Let $E$ be relatively closed in an open subset $U$ of $\R^n$. An Almgren competitor for $E$ is a relatively closed set $F\subset U$ that can be written as $F=\varphi_1(E)$, where $\varphi_t:U\to U,t\in [0,1]$ is a family of continuous mappings such that 
\be \varphi_0(x)=x\mbox{ for }x\in U;\ee
\be\mbox{ the mapping }(t,x)\to\varphi_t(x)\mbox{ of }[0,1]\times U\mbox{ to }U\mbox{ is continuous;}\ee
\be\varphi_1\mbox{ is Lipschitz,}\ee
  and if we set $W_t=\{x\in U\ ;\ \varphi_t(x)\ne x\}$ and $\widehat W=\bigcup_{t\in[0.1]}[W_t\cup\varphi_t(W_t)]$,
then
\be \widehat W\mbox{ is relatively compact in }U.\ee
 
Such a $\varphi_1$ is called a deformation in $U$, and $F$ is also called a deformation of $E$ in $U$.
\end{defn}

\begin{defn}[Almgren minimal sets and minimal cones]
Let $0<d<n$ be integers, $U$ be an open set of $\R^n$. An Almgren minimal set $E$ in $U$ is a minimal set defined in Definition \tb{2.1} while taking the competitor class $\mathscr F$ to be the class of all Almgren competitors for $E$.

An Almgren minimal set which is a cone is called a minimal cone.\end{defn}

For future convenience, we also have the following more general definition:

\begin{defn}Let $U\subset \R^n$ be an open set, and let $E\subset \R^n$ be a closed set (not necessarily contained in $U$). We say that $E$ is Almgren minimal in $U$, if $E\cap U$ is Almgren minimal in $U$. A closed set $F\subset \R^n$ is called a  deformation of $E$ in $U$, if $F=(E\bs U)\cup \varphi_1(E\cap U)$, where $\varphi_1$ is a deformation in $U$.
\end{defn}

\begin{rem}Since Almgren minimal sets are more often used, we usually omit the word ''Almgren'' and call them minimal sets.
\end{rem}

\begin{defn}[Topological competitors] Let $G$ be an abelian group. Let $E$ be a closed set in an open domain $U$ of $\R^n$. We say that a closed set $F$ is a $G$-topological competitor of dimension $d$ ($d<n$) of $E$ in $U$, if there exists an open convex subset $B$ with $\bar B\subset U$ such that

1) $F\bs B=E\bs B$;

2) For all Euclidean $n-d-1$-topological sphere $S\subset U\bs(B\cup E)$, if $S$ represents a non-zero element in the singular homology group $H_{n-d-1}(U\bs E;G)$, then it is also non-zero in $H_{n-d-1}(U\bs F;G)$.
We also say that $F$ is a $G$ topological competitor of $E$ in $B$.
\end{defn}

And Definition \tb{2.1} gives the definition of $G$-topological minimizers in a domain $U$ when we take the competitor class to be the classe of $G$-topological competitors of $E$.

The simplest example of a topological minimal set is a $d-$dimensional plane in $\R^n$.  

\begin{pro}[cf.\cite{topo} Proposition 3.7 and Corollary 3.17]  

$1^\circ$ Let $E\subset \R^n$ be closed. Then for any $d<n$, any abelian group $G$, and any convex set $B$, $B'$ such that $\bar B'\subset B^\circ$, every Almgren competitor of $E$ in $B'$ is a $G$-topological competitor of $E$ in $B$ of dimension $d$.

$2^\circ$ For any abelian group $G$, all $G$-topological minimal sets are Almgren minimal in $\R^n$.
\end{pro}

\begin{rem}The notion of (Almgren or topological) minimal sets does not depend much on the ambient dimension. One can easily check that $E\subset U$ is $d-$dimensional Almgren minimal in $U\subset \R^n$ if and only if $E$ is Almgren minimal in $U\times\R^m\subset\R^{m+n}$, for any integer $m$. The case of topological minimality is proved in \cite{topo} Proposition 3.18.\end{rem}

At last, we give the minimality for $Y\times Y$, which is the prerequisite of this article.

\begin{thm}[\cite{YXY} Theorems 6.1 and 6.6] The set $Y\times Y$ is Almgren minimal, and $\Z_2$-topological minimal.
\end{thm}

\subsection{Regularity results for minimal sets}

We now begin to list some regularity results for minimal sets that will be needed later. They are in fact regularity results for Almgren minimal sets, but they also hold for all $G$-topological minimizers, after Proposition \tb{2.7}. 

\begin{defn}[reduced set] Let $U\subset \R^n$ be an open set. For every closed subset $E$ of $U$, denote by
\be E^*=\{x\in E\ ;\ \H^d(E\cap B(x,r))>0\mbox{ for all }r>0\}\ee
 the closed support (in $U$) of the restriction of $\H^d$ to $E$. We say that $E$ is reduced if $E=E^*$.
\end{defn}

It is easy to see that
\be \H^d(E\bs E^*)=0.\ee
In fact we can cover $E\bs E^*$ by countably many balls $B_j$ such that $\H^d(E\cap B_j)=0$.

\begin{rem}
 It is not hard to see that if $E$ is Almgren minimal (resp. $G$-topological minimal), then $E^*$ is also Almgren minimal (resp. $G$-topological minimal). As a result it is enough to study reduced minimal sets. An advantage of reduced minimal sets is, they are locally Ahlfors regular (cf. Proposition 4.1 in \cite{DS00}). Hence any approximate tangent plane of them is a true tangent plane. Since minimal sets are rectifiable (cf. \cite{DS00} Theorem 2.11 for example), reduced minimal sets admit true tangent $d$-planes almost everywhere.
\end{rem}

If we regard two sets to be equivalent if they are equal modulo $\H^d$-null sets, then a reduced set is always considered to be a good (in the sense of regularity) representative of its equivalent class. 

\textbf{In the rest of the article, we only consider reduced sets.}

\begin{defn}[blow-up limit] Let $U\subset\R^n$ be an open set, let $E$ be a relatively closed set in $U$, and let $x\in E$. Denote by $E(r,x)=r^{-1}(E-x)$. A set $C$ is said to be a blow-up limit of $E$ at $x$ if there exists a sequence of numbers $r_n$, with $\lim_{n\to \infty}r_n=0$, such that the sequence of sets $E(r_n,x)$ converges to $C$ for the local Hausdorff distance in any compact set of $\R^n$.
\end{defn}

\begin{rem}
A set $E$ might have more than one blow-up limit at a point $x$. However it is not known yet whether this can happen to minimal sets. 
 
 When a set $E$ admits a unique blow-up limit at a point $x\in E$, denote this blow-up limit by $C_xE$.
\end{rem}

\begin{pro}[c.f. \cite{DJT} Proposition 7.31]Let $E$ be a reduced Almgren minimal set in an open set $U$ of $\R^n$, and let $x\in E$. Then every blow-up limit of $E$ at $x$ is a reduced Almgren minimal cone $F$ centred at the origin, and $\H^d(F\cap B(0,1))=\theta(x):=\lim_{r\to 0} r^{-d}\H^d(E\cap B(x,r)).$\end{pro}

An Almgren minimal cone is just a cone which is also Almgren minimal. We will call them minimal cones throughout this paper, since we will not talk about any other type of minimal cones. 

\begin{rem} Obviously, a cone in $\R^n$ is minimal if and only if it is minimal in the unit ball, if and only if it is minimal in any open subset containing the origin.
\end{rem}

We now focus on regularity results for Almgren minimal sets of dimension 2.

\begin{defn}[bi-H\"older ball for closed sets] Let $E$ be a closed set of Hausdorff dimension 2 in $\R^n$. We say that $B(0,1)$ is a bi-H\"older ball for $E$, with constant $\tau\in(0,1)$, if we can find a 2-dimensional minimal cone $Z$ in $\R^n$ centered at 0, and $f:B(0,2)\to\R^n$ with the following properties:

$1^\circ$ $f(0)=0$ and $|f(x)-x|\le\tau$ for $x\in B(0,2);$

$2^\circ$ $(1-\tau)|x-y|^{1+\tau}\le|f(x)-f(y)|\le(1+\tau)|x-y|^{1-\tau}$ for $x,y\in B(0,2)$;

$3^\circ$ $B(0,2-\tau)\subset f(B(0,2))$;

$4^\circ$ $E\cap B(0,2-\tau)\subset f(Z\cap B(0,2))\subset E.$  

In this case we also say that B(0,1) is of type $Z$ for $E$.

We say that $B(x,r)$ is a bi-H\"older ball for $E$ of type $Z$ (with the same parameters) when $B(0,1)$ is a bi-H\"older ball of type $Z$ for $r^{-1}(E-x)$.
\end{defn}

\begin{thm}[Bi-H\"older regularity for 2-dimensional Almgren minimal sets, c.f.\cite{DJT} Thm 16.1]\label{holder} Let $U$ be an open set in $\R^n$ and $E$ a reduced Almgren minimal set in $U$. Then for each $x_0\in E$ and every choice of $\tau\in(0,1)$, there is an $r_0>0$ and a minimal cone $Z$ such that $B(x_0,r_0)$ is a bi-H\"older ball of type $Z$ for $E$, with constant $\tau$. Moreover, $Z$ is a blow-up limit of $E$ at $x$.
\end{thm}

\begin{defn}[point of type $Z$] 

$1^\circ$ In the above theorem, we say that $x_0$ is a point of type $Z$ (or $Z$ point for short) of the minimal set $E$. The set of all points of type $Z$ in $E$ is denoted by $E_Z$. 

$2^\circ$ In particular, we denote by $E_P$ the set of regular points of $E$ and $E_Y$ the set of $\Y$ points of $E$.  Set $E_S:=E\bs E_P$ the set of all singular points in $E$.
\end{defn}

\begin{rem} Again, since we might have more than one blow-up limit for a minimal set $E$ at a point $x_0\in E$, the point $x_0$ might be of more than one type (but all the blow-up limits at a point are bi-H\"older equivalent). However, if one of the blow-up limits of $E$ at $x_0$ admits the``full-length'' property (see Remark \tb{2.21}), then in fact $E$ admits a unique blow-up limit at the point $x_0$. Moreover, we have the following $C^{1,\a}$ regularity around the point $x_0$. In particular, the blow-up limit of $E$ at $x_0$ is in fact a tangent cone of $E$ at $x_0$.
\end{rem}

\begin{thm}[$C^{1,\a}-$regularity for 2-dimensional minimal sets, c.f. \cite{DEpi} Thm 1.15]\label{c1} Let $E$ be a 2-dimensional reduced minimal set in the open set $U\subset\R^n$. Let $x\in E$ be given. Suppose in addition that some blow-up limit of $E$ at $x$ is a full length minimal cone (see Remark \tb{2.21}). Then there is a unique blow-up limit $X$ of $E$ at $x$, and $x+X$ is tangent to $E$ at $x$. In addition, there is a radius $r_0>0$ such that, for $0<r<r_0$, there is a $C^{1,\a}$ diffeomorphism (for some $\a>0$) $\Phi:B(0,2r)\to\Phi(B(0,2r))$, such that $\Phi(0)=x$ and $|\Phi(y)-x-y|\le 10^{-2}r$ for $y\in B(0,2r)$, and $E\cap B(x,r)=\Phi(X)\cap B(x,r).$ 

We can also ask that $D\Phi(0)=Id$. We call $B(x,r)$ a $C^1$ ball for $E$ of type $X$.
\end{thm}

\begin{rem} $1^\circ$ We are not going to give the precise definition of the full length property. Instead, for purpose of this paper, it is enough to know that the 2-dimensional planes and the 2-dimensional $\Y$ sets are minimal cones of full length (\cite{DEpi}), and hence points of type $\P$ and $\Y$ admit the above $C^1$ regularity.

$2^\circ$ As a result, after Theorem \tb{2.17}, a blow-up limit of a reduced minimal set $E$ at a point $x\in E$ is a plane if and only if the plane is the unique approximate tangent plane of $E$ at $x$. Same for $\Y$ points.
\end{rem}

After Remark \tb{2.21}, for any reduced minimal set $E$ of dimension $d$, and for any $x\in E$ at which an approximate tangent $d$-plane exists (which is true for a.e. $x\in E$), $T_xE$ also denotes the tangent plane of $E$ at $x$, and the blow-up limit of $E$ at $x$. 

\begin{thm}[Structure of 2-dimensional minimal cones in $\R^n$, cf. \cite{DJT} Proposition 14.1] Let $K$ be a reduced 2-dimensional minimal cone in $\R^n$, and let $X=K\cap \partial B(0,1)$. Then $X$ is a finite union of great circles and arcs of great circles $C_j,j\in J$. The arcs $C_j$ can only meet at their endpoints, and each endpoint is a common endpoint of exactly three $C_j$, which meet with $120^\circ$ angles (such an endpoint is called a $\Y$ point in $K\cap \pa B$). In addition, the length of each $C_j$ is at least $\eta_0$, where $\eta_0>0$ depends only on the ambient dimension $n$.
\end{thm}

An immediate corollary of the above theorem is the following:

\begin{cor}
$1^\circ$ If $C$ is a minimal cone of dimension 2, then for the set $C_P$ of regular points of $C$, each of its connected components is a sector. 

$2^\circ$ Let $E$ be a 2-dimensional minimal set in $U\subset \R^n$. Then $\bar E_Y=E_S$.

$3^\circ$ The set $E_S\bs E_Y$ is composed of isolated points. \end{cor}

As a consequence of the $C^1$ regularity for regular points and $\Y$ points, and Corollary \tb{2.23}, we have
\begin{cor}Let $E$ be an 2-dimensional Almgren minimal set in a domain $U\subset \R^n$. Then 

$1^\circ$ The set $E_P$ is open in $E$;

$2^\circ$ The set $E_Y$ is a countable union of $C^1$ curves. The endpoints of these curves are either in $E_S\bs E_Y$, or lie in $\partial U$.  
\end{cor}

%

\section{The associated convex domain and stabilities for 2-dimensional minimal cones}

In this section, we introduce the definition for various stabilities for 2-dimensional minimal cones. Some specifications will also be given for the convex domain associated to the set $Y\times Y$.

The idea of sliding stability can be traced back to the concept of ''minimal set with sliding boundary'' first proposed by Guy David in \cite{GD13, GD14}. Minimal set with sliding boundary is a new model for Plateau's problem, considering that the soap film can slide along a given boundary. This model enables people to study boundary regularity for solutions of Plateau's problem. 

The stabilities in our paper were defined in \cite{2T} for a different and more technical purpose: they are quantitative properties for 2-dimensional minimal cones used to control local regularities for 2-dimensional minimal sets.

\subsection{The associated convex domain $\cU(K,\eta)$}

We first give the general definition of the convex domain $\cU(K,\eta)$ associated to an arbitrary 2-dimensional minimal cone $K\subset \R^n$, based on the Theorem \tb{2.22} for structures of 2-dimensional minimal cones. 

So take any 2-dimensional minimal cone $K\subset \R^n$. Denote by $B$ the unit ball of $\R^n$. Then by Theorem \tb{2.22}, $K\cap\pa B$ is a union of circles $\{s_j,1\le j\le \mu\}$, and arcs of great circles with only $\Y$ type junctions. Let $\eta_0(K)$ denote the minimum of length of these arcs. It is positive, again by Theorem \tb{2.22}. 

Denote by $\{a_j, 1\le j\le m\}$ the set of $\Y$ points in $K\cap \pa B$. 

\begin{defn}For any $\eta>0$, the $\eta$-convex domain $\cU(K,\eta)$ for $K$ is defined as
\be \cU(K,\eta)=\{x\in B: \lg x, y\rg< 1-\eta, \forall y\in K\cap \pa B\mbox{ and }\lg x,a_j\rg<1-\sqrt\eta, \forall 1\le j\le m\}\subset \R^n. \ee 
\end{defn}

From the definition, we see directly that $\cU=\cU(K,\eta)$ is obtained by "cutting off" some small part of the unit ball $B$. More precisely, we first take the unit ball $B$, then just like peelling an apple, we use a knife to peel a thin band (with width about $2\sqrt \eta$) near the set $K\cap \pa B$. Then after this operation, the ball $B$ stays almost the same, except that near the set $K$, the boundary surface will be a thin cylindrical surface. This is the condition ''$\lg x, y\rg  < 1-\eta, \forall y\in K\cap \pa B$''. Next we turn to the singular points $a_j$: they are isolated, so we make one cut at each point, perpendicular to the radial direction, to get a small planar surface near each $a_j$, of diameter about $4\eta^\frac 14$. This follows from the condition ''$\lg x,a_j\rg  <1-\sqrt\eta, \forall 1\le j\le m $''. 
\medskip

Now let us introduce some notations. 

Fix $K$ and $\eta$. Let $\cU=\cU(K,\eta)$. 

For $1\le j,l\le m$, let $\gamma_{jl}$ denote the arc of great circle of $\pa B$ that connects $a_j$ and $a_l$, if it exists; otherwise set $\gamma_{jl}=\emptyset$. Set $J=\{(j,l): 1\le j,l\le m\mbox{ and }\gamma_{jl}\ne\emptyset\}$, which is exactly the set of pairs $(j,l)$ such that the $\Y$ points $a_j$ and $a_l$ are connected directly by an arc of great circle on $K$. 
 
Denote by $A_j$ the $(n-1)$-dimensional planar part centered at $(1-\sqrt\eta)a_j$ of $\partial \cU$. That is,  
\be A_j=\{x\in \bar \cU: \lg x,a_j\rg  =1-\sqrt\eta\}.\ee
Let $A=\cup_{1\le j\le m}A_j$.

Set 
\be\Gamma_{jl}=\{x\in \bar \cU, \lg x,y\rg  =1-\eta\mbox{ for some }y\in \gamma_{jl}\}\bs \mA,\ee
with $\mA$ being the intersection with $\bar\cU$ of the cone over $A$ centered at 0,
and
\be S_j=\{x\in \bar \cU,\lg x,y\rg  =1-\eta\mbox{ for some }y\in s_j\}.\ee
Then $\Gamma_{jl}$ is the band like part of $\pa \cU$ near each $(1-\eta)\gamma_{jl}$, and similar for $S_j$. The union
$\Gamma=\cup_{1\le j,l\le m}\Gamma_{jl}$ together with $S=\cup_{1\le j\le \mu}S_j$ are the whole cylindrical parts of $\pa \cU$.

Let $\mathfrak c_{jl}$ and $\ms_j$ denote the parts of the cone (centered at 0) included in $\bar\cU$ over $\gamma_{jl}$ and $s_j$ respectively.

\medskip

Let $\eta_1(K)\in (0,10^{-1}]$ be the supremum of the number $\eta$, such that $K\cap \pa\cU(K,\eta)$ is a deformation retract of $A\cup \G\cup S$, and such that $\eta_1(K)\le 10^{-1}$. The fact that $\eta_1(K)>0$ is because of the structure Theorem \tb{2.22}, which implies that $K\cap \pa\cU(K,\eta)$ is a finite union of piecewise smooth curves.

Note that on $\pa\cU(K,\eta)$, any 2 of the $\Gamma_{jl}, (j,l)\in J$ and $S_j,1\le j\le \mu$ never have a common point, and the $A_j,1\le j\le m$ are disjoint.

\medskip

For any fixed $\eta<\eta_1(K)$, let $R_1=R_1(\eta)=\sqrt{1-(1-\eta)^2}$, $R_2=R_2(\eta)=\sqrt{1-(1-\sqrt\eta)^2}$, $R_3=R_3(\eta)=\sqrt{(1-\eta)^2-(1-\sqrt\eta)^2}$. Then $R_2>R_3>2R_1$ since $\eta<10^{-1}$.

For a finer geometric description for a general $\cU(K,\eta)$, see \cite{stablePYT}. In this article we will only consider minimal cones of type $Y\times Y$, and will describe its associated convex domain in detail, in the next subsection.

\begin{defn}[$\d$-sliding deformation] Let $U$ be an open subset of $\R^n$, let $E\subset \R^n$ be closed. For $\d>0$, a $\d$-sliding Lipschitz deformation of $E$ in $\bar U$ is a set $F\subset \bar U$ that can be written as $F=\varphi_1(E\cap \bar U)$, where $\varphi_t: \bar U\to \bar U, 0\le t\le 1$ is a family of continuous mappings such that
\be \varphi_0(x)=x\mbox{ for }x\in \bar U;\ee
\be \mbox{ the mapping }(t,x)\to \varphi_t(x)\mbox{ of }[0,1]\times \bar U\mbox{ to }\bar U\mbox{ is continuous};\ee
\be \varphi_1\mbox{ is Lipschitz },\ee
\be \varphi_t(\partial U)\subset\partial U\mbox{ for every } 0\le t\le 1,\ee
and
\be|\varphi_t(x)-x|<\d\mbox{ for all }x\in E\cap\partial U\mbox{ and }0\le t\le 1.\ee

Such a $\varphi_1$ (that satisfies \tb{(3.5)-(3.8)}) is called a sliding deformation in $\bar U$, and $F$ is called a $\d$-sliding deformation of $E$ in $\bar U$.
\end{defn}

\begin{rem}Note that the ''$\d$-sliding condition'' \tb{(3.9)} only imposes condition on $\varphi_t$ on the boundary $E\cap \pa U$, but not the whole set $E$. Hence if $E\subset\R^n$ is closed, $\varphi_1$ is a sliding deformation in $\bar U$, and $\varphi_1(E\cap\bar U)$ is a $\d$-sliding deformation of $E$ in $\bar U$, then for any closed subset $F\subset \R^n$ so that $F\bs U=E\bs U$, $\varphi_1(F\cap\bar U)$ is also a $\d$-sliding deformation of $F$ in $\bar U$.
\end{rem}

\begin{defn}[$\d$-Almgren sliding competitors] Let $\F_\d(E,\bar U)$ denote the class of all $\d$-sliding deformations of $E$ in $\bar U$, and let $\oF_\d(E,\bar U)$ be the family of sets that are Hausdorff limits of sequences in $\F_\d(E,\bar U)$. That is: we set
\be \begin{split}\oF_\d(E,U)&=\{F\subset\bar U: \tb{(2.1)}\mbox{ holds for }F\mbox{, and }\exists \{E_n\}_n\subset\F_\d(E,U)\mbox{ such that }d_H(E_n,F)\to 0
\}.\end{split}\ee

Elements in $\oF_\d(E,\bar U)$ are called $\d$-Almgren sliding competitors of $E$ in $\bar U$.\end{defn}

Due to the specific structure of 2-dimensional minimal cones (Theorem \tb{2.22}), we also give the definition of $(\eta,\d,\nu, L)$-sliding competitor for 2-dimensional minimal cones. Roughly speaking, an $(\eta,\d,\nu, L)$-sliding competitor for a 2-dimensional minimal cone $K$, is just a $\d$-sliding competitor $F$ of $K$ in $\bar\cU(K,\eta)$, whose boundary $F\cap\pa\cU(K,\eta)$ is a $L$-Lipschitz graph of $K\cap \pa\cU(K,\eta)$ near regular parts of $K\cap \pa\cU(K,\eta)$.

So fix $\eta<\eta_1(K)$. Recall that for each $(j,l)\in J$, $\mc_{jl}$ is the part inside $\bar\cU(K,\eta)$ of the cone over $\g_{jl}$. Let $Q_{jl}$ denote the $n-2$-dimensional linear subspace orthogonal to $\mc_{jl}$ in $\R^n$. Then $I_{jl}:=\G_{jl}\cap A_j$ is a $n-2$ disk contained in $Q_{jl}$, with center $x_{jl}$. Set $l_{jl}:=\mc_{jl}\cap \pa\cU(K,\eta)$. It is a piecewise linear curve that connects the centers $(1-\sqrt\eta)a_j$ and $(1-\sqrt\eta)a_l$ of $A_j$ and $A_l$, which is the union of the arc $[(1-\eta)\g_{jl}]\cap\G_{jl}$ and the two segments $[(1-\sqrt\eta)a_j,x_{jl}]$ and $[(1-\sqrt\eta)a_l, x_{lj}]$. 

For each $\nu\le R_3$, let $l_{jl}^\nu=l_{jl}\bs B((1-\sqrt\eta)a_j, \nu)\cup B((1-\sqrt\eta)a_l,\nu)$. Then $l_{jl}^\nu$ is the union of the arc $[(1-\eta)\g_{jl}]\cap\G_{jl}$ and two segments contained in $[(1-\sqrt\eta)a_j,x_{jl}]$ and $[(1-\sqrt\eta)a_l, x_{lj}]$. In particular, $l_{jl}^0=l_{jl}$, and $l_{jl}^{R_3}=[(1-\eta)\g_{jl}]\cap\G_{jl}$. Let $\ml_{jl}^\nu$ denote the cone over $l_{jl}^\nu$.

Moreover, for $0\le \nu\le R_3$, 
\be\begin{split} \G_{jl}&=l_{jl}^{R_3}\times B_{Q_{jl}}(0,R_1)\subset l_{jl}^\nu\times B_{Q_{jl}}(0,R_1)\\
&\subset l_{jl}^0\times B_{Q_{jl}}(0,R_1)=\G_{jl}\cup A_j\cup A_l\subset \pa \cU\cap B(K\cap \pa\cU, R_1).
\end{split}\ee

Also, for $1\le j\le \mu$, let $Q_j$ denote the subspace orthogonal to $\ms_j$.

Note that all the notations above depend on the given $\eta$.

\begin{defn}Let $K$ be a 2-dimensional Almgren minimal cone in $\R^n$. Let $\eta<\eta_1(K)$. Take all the notations above (which depend on $\eta$ obviously). Let $\d \in(0,R_1)$, $\nu\in [0,R_3]$, $L>0$. A closed set $G\in \pa\cU(K,\eta)$ is called an $(\eta,\d,\nu, L)$-sliding boundary of $K$, if it is a $\d$-sliding deformation of $\pa\cU(K,\eta)\cap K$ in $\bar\cU(K,\eta)$, and it satisfies in addition,

$1^\circ$ For each $(j,l)\in J$, $G\cap (l_{jl}^\nu\times B_{Q_{jl}}(0,R_1))$ is the graph of an $L$-Lipschitz map from $l_{jl}^\nu$ to $B_{Q_{jl}}(0,R_1)$;

$2^\circ$, For each $1\le j\le \mu$, $G\cap (s_j\times B_{Q_{j}}(0,R_1))$ is the graph of an $L$-Lipschitz map from $\ms_j$ to $B_{Q_{j}}(0,R_1)$.

\end{defn}

\begin{defn}[sliding competitors]Let $K$ be a 2-dimensional Almgren minimal cone in $\R^n$. Let $\eta<\eta_1(K)$. Take all the notations above (which depend on $\eta$ obviously). Let $\d \in(0,R_1)$, $\nu\in [0,R_3]$, $L>0$.  

A closed set $F\subset \bar\cU(K,\eta)$ is called an $(\eta,\d)$-Almgren sliding competitor for $K$, if it is a $\d$-sliding competitor for $K$ in $\bar\cU(K,\eta)$; it is called an $(\eta,\d,\nu, L)$-Almgren sliding competitor for $K$, if it is a $\d$-sliding competitor for $K$ in $\bar\cU(K,\eta)$, and $F\cap \pa \cU(K,\eta)$ is an $(\eta,\d,\nu, L)$-sliding boundary of $K$.
\end{defn}

\begin{defn}[Stable minimal cones] Let $K$ be a 2-dimensional Almgren minimal cone in $\R^n$.

$1^\circ$ We say that $K$ is $(\eta,\d,\nu, L)$-Almgren sliding stable, if for some $\eta\in (0,\eta_1(K))$, $L>0$, $\d\in (0,R_1)$, and $\nu\in (0,R_3)$, \tb{(2.1)} holds, and for all $(\eta,\d,\nu, L)$-sliding competitors $F$ for $K$ we have
\be \H^2(K\cap \bar\cU(K,\eta))\le \H^2(F);\ee

$2^\circ$ We say that $K$ is $(\eta, \d)$-Almgren sliding stable, if for some $\eta\in (0,\eta_1(K))$, and $\d\in (0,R_1)$, \tb{(2.1)} holds, and \tb{(3.12)} holds for all $(\eta,\d)$-Almgren sliding competitors for $K$;

$3^\circ$ We say that $K$ is Almgren sliding stable if it is $(\eta,\d,\nu, L)$-Almgren sliding stable for some $\eta\in (0,\eta_1(K))$, $L>0$, $\d\in (0,R_1)$, and $\nu\in (0,R_3)$.
\end{defn}

Similarly for 2-dimensional topological minimal cones, we have

\begin{defn}Let $K$ be a 2-dimensional $G$-topological minimal cone in $\R^n$. Let $0<\eta<\eta_1(K)$, and $0<\d<R_1(\eta)$.

$1^\circ$ We say that a closed set $F$ is an $(\eta,\d)$-$G$-topological sliding competitor for $K$, if there exists a 2-dimensional $G$-topological competitor $E$ of $K$ in $\cU(K,\eta)$, such that $F$ is a $\d$-sliding deformation of $E$ in $\bar\cU(K,\eta)$;

$2^\circ$ For any $L>0$ and $\nu\in (0,R_3)$, we say that a closed set $F$ is an $(\eta,\d,\nu, L)$-$G$-topological sliding competitor for $K$, if it is an $(\eta,\d)$-$G$-topological sliding competitor for $K$, and $F\cap \pa \cU(K,\eta)$ is an $(\eta,\d,\nu, L)$-sliding boundary of $K$;

$3^\circ$ We say that $K$ is $(\eta,\d,\nu, L)$-$G$-(resp. $(\eta,\d)$-$G$-) topological sliding stable, if for all $(\eta,\d,\nu, L)$-$G$-(resp. $(\eta,\d)$-$G$-)topological sliding competitor $F$ of $K$, \tb{(3.12)} holds;

$4^\circ$ We say that $K$ is $G$-topological sliding stable, if it is $(\eta,\d,\nu, L)$-$G$-topological sliding stable for some $\eta\in (0,\eta_1(K))$, $\d\in (0, R_1)$, $L>0$ and $\nu\in (0,R_3)$.
\end{defn}

\begin{rem}We can see directly from the above definitions, that

$1^\circ$ $(\eta,\d, L,\nu)$-Almgren sliding competitors are $(\eta,\d)$-Almgren sliding competitors, hence $(\eta, \d)$-Almgren sliding stability implies $(\eta,\d, L,\nu)$-Almgren sliding stability;

$2^\circ$ Similarly, $(\eta,\d, L,\nu)$-$G$-topological sliding competitors are $(\eta,\d)$-$G$-topological sliding competitors, hence $(\eta, \d)$-$G$-topological sliding stability implies $(\eta,\d, L,\nu)$-$G$-topological sliding stability.
\end{rem}

\subsection{The convex domain $\cU(\eta)$ for $Y\times Y$}

In the previous section we gave necessary definitions for general 2-dimensional minimal cones. Since the cone $Y\times Y$ is the main object of this article, in this section we will give some specifications concerning the corresponding notions defined in the last subsection for $Y\times Y$.

Write $\R^4=P_1\times P_2$, where $P_1$ and $P_2$ are orthogonal 2-dimensional subspaces of $\R^4$. For any point $x$ in $\R^4$, write $x=(x_1,x_2)$, where $x_i\in P_i$.

Let $B=B(0,1)$. For $i=1,2$, $Y_i\subset P_i$ is a 1-dimensional $\Y$ set centered at the origin of $P_i$. Denote by $Z$ the set $Y_1\times Y_2\subset \R^4$. 

Denote by $o_1$ and $o_2$ the origin of $P_1$ and $P_2$. Let $B_i=B_{P_i}(0,1)$, $i=1,2$. Denote by $a_i\in P_1,i=1,2,3$ the three points of intersection of $Y_1$ with $B_1$, and let $b_j\in P_2,j=1,2,3$ denote the three points of intersections of $Y_2$ with $\pa B_2$. 

Then the 6 points of type $\Y$ in $Z\cap \pa B$ are $c_{1i}:=(a_j,o_2), i=1,2,3$ and $c_{2j}:=(o_1,b_j),j=1,2,3$. For $1\le i,j\le 3$, let $\gamma_{ij}$ denote the minor arc of great circle (or equivalently the geodesic) that connects $c_{1i}$ and $c_{2j}$. Then 
\be Z\cap \pa B=\cup_{1\le i\le 3}\cup _{1\le j\le 3}\gamma_{ij}.\ee

Set $\eta_0=\eta_0(Z)$, $\eta_1=\eta_1(Z)$. For any $\eta\in (0,\eta_1)$ small, recall that the $\eta$-convex domain for $Z$ is 
\be \cU(\eta):=\{x\in \bar B:\lg x,y\rg <1-\eta, \forall y\in Z\cap\pa B\mbox{ and }\lg x,c_\a\rg <1-\sqrt\eta,\a\in C\},\ee
where $C=\{11,12,13,21,22,23\}$.

For each $\a\in C$, denote by $A_\a=A_\a(\eta)$ the $3$-dimensional planar part centered at $(1-\sqrt\eta)c_\a$ of $\pa \cU(\eta)$. That is,
\be A_\a=\{x\in \bar B: \lg x,c_\a\rg =1-\sqrt\eta\mbox{ and }\lg x,y\rg \le 1-\eta, \forall y\in Z\cap \pa B\}.\ee
Let $A=A(\eta)=\cup_{\a\in C}A_\a(\eta)$.

Set 
\be\G_{ij}=\G_{ij}(\eta):=\{x\in \bar B, \lg x,y\rg =1-\eta\mbox{ for some }y\in \gamma_{ij}\}\bs \mA,\ee
with $\mA=\mA(\eta)$ being the cone over $A$ centered at 0. Then $\G_{ij}$ is the band like part of $\pa \cU(\eta)$ near each $(1-\eta)\g_{ij}$. The union $\G=\G(\eta):=\cup_{1\le i,j\le 3}\G_{ij}(\eta)$ is the whole cylinderical part of $\pa \cU(\eta)$.


\medskip

Fix an $\eta<\eta_1$. Let us now describe geometrically (after some simple calculation) the shape of each part $A_\a, \a\in C$, and $\G_{ij}, 1\le i,j\le 3$ (which depends on $\eta$ as before).

Fix an $1\le i\le 3$. Then the shape of $A_{1i}$ is the following: then condition $\lg x,c_{1i}\rg =1-\sqrt\eta$ gives the 3-dimensional disk $D_{1i}$ centered at $o_{D_{1i}}=(1-\sqrt\eta)c_{1i}$, perpendicular to $c_{1i}$, and with radius $R_2$. Then the intersection of $Z$ with $D_{1i}$ coincides with the part of $\{(1-\sqrt\eta)c_{1i}\}\times Y_2$ inside $D_{1i}$. The condition $\lg x,y\rg \le 1-\eta, \forall y\in Z$ is equivalent to the condition $\lg x,y\rg \le 1-\eta, \forall y\in \cup_{1\le j\le 3}\g_{ij}$. For any $1\le j\le 3$, the set $I_{ij}=\{x\in D_{1i}:\lg x,y\rg =1-\eta$ for some $y\in \g_{ij}\}$ is a 2-dimensional disk perpendicular to $c_{1i}$ and $c_{2j}$, centered at the point $o_{I_{ij}}=(1-\sqrt\eta)c_{1i}+R_3 c_{2j}=((1-\sqrt\eta)a_i,R_3b_j)$, with radius $R_1$. By definition, $I_{ij}$ separates $D_{1i}$ into two parts, and we keep the one that contains the center $o_{D_{1i}}$ of $D_{1i}$ and throw away the other part. We do this for $1\le j\le 3$, and get $A_{1i}$. Note that the boundary of $A_{1i}$ is the union of $I_{ij}, 1\le j\le 3$, and the rest of the sphere $\pa D_{1i}$. See the picture below.

\centerline{\includegraphics[width=0.7\textwidth]{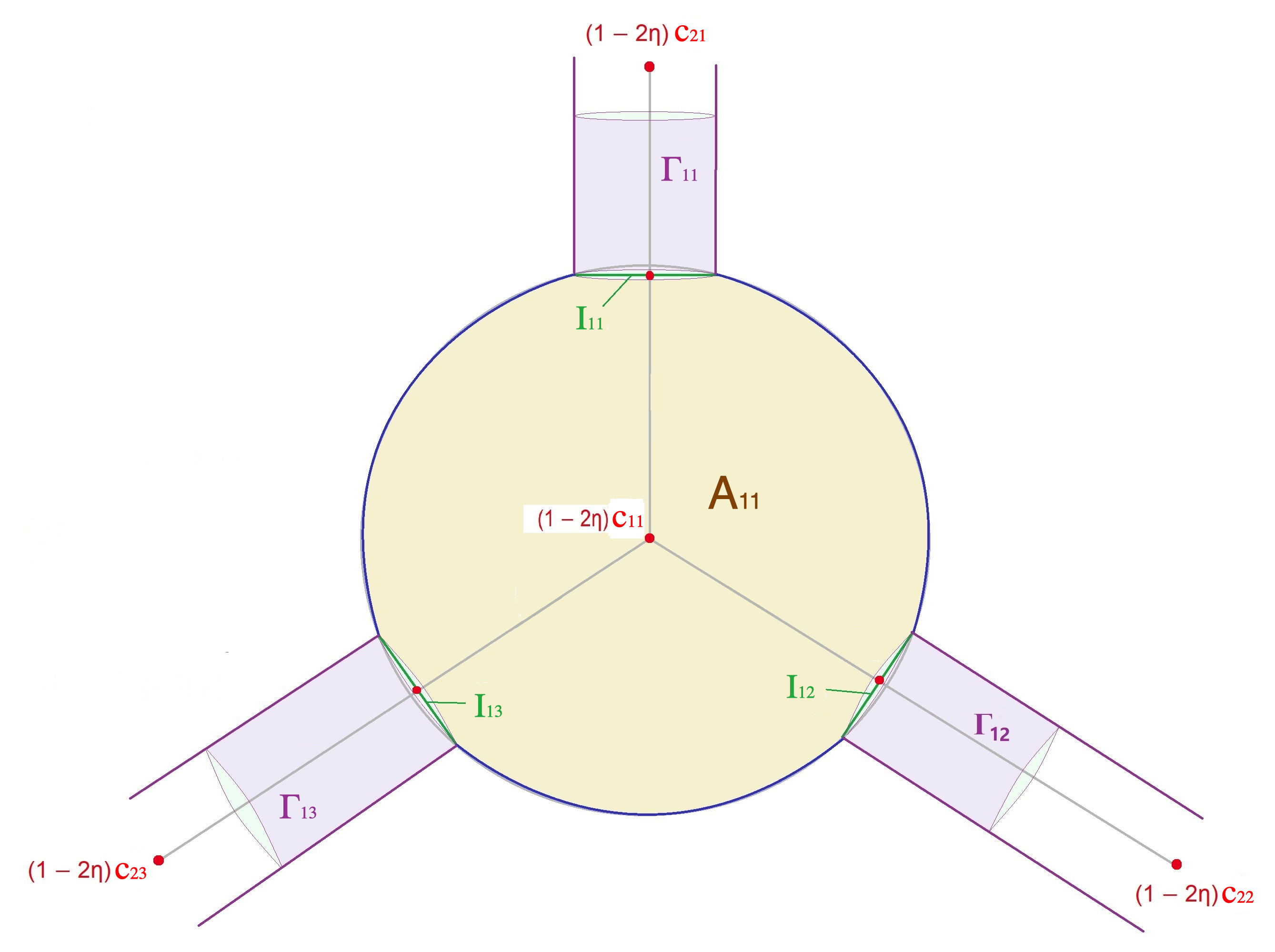}}
\nopagebreak[4]
\centerline{The 3-dimensional planar region $A_{11}$}

Similarly, fix any $1\le j\le 3$, to get $A_{2j}$, let $D_{2j}$ denote the 3-dimensional disk centered at $o_{D_{2j}}=(1-\sqrt\eta)c_{2j}$, perpendicular to $c_{2j}$, with radius $R_2$. The intersection of $Z$ with $D_{2j}$ coincides with the part of $Y_1\times\{(1-\sqrt\eta)c_{2j}\}$ inside $D_{2j}$. For each $1\le i\le 3$, set $J_{ij}=\{x\in D_{2j}:\lg x,y\rg =1-\eta$ for some $y\in \g_{ij}\}$. Then it is a 2-dimensional disk perpendicular to $c_{1i}$ and $c_{2j}$ (hence is parallel to $I_{ij}$), centered at the point $o_{J_{ij}}=R_3c_{1i}+(1-\sqrt\eta)c_{2j}$ with radius $R_1$. Each $I_{ij}$ separates $D_{2j}$ into two parts, and we throw away the part that does not contain the center $o_{D_{2j}}$ of $D_{2j}$. We do this for $1\le i\le 3$. The rest part of $D_{2j}$ is the set $A_{2j}$, whose boundary is the union of $J_{ij},1\le j\le 3$, and the rest of the sphere $\pa D_{2j}$.

\medskip

The structure of $\G_{ij}$ is easier to describe: let $\theta_0=\theta_0(\eta)=\arcsin \frac{R_3}{1-\eta}$. Then it is easy to see that $o_{I_{ij}}=(1-\eta)(\cos\theta_0 c_{1i}+\sin\theta_0c_{2j})=((1-\eta)\cos\theta_0a_i,(1-\eta)\sin\theta_0b_j)$ and $o_{J_{ij}}=(1-\eta)(\sin\theta_0 c_{1i}+\cos\theta_0c_{2j})=((1-\eta)\sin\theta_0a_i,(1-\eta)\cos\theta_0b_j)$. Let 
\be \g^0_{ij}:=\{(\cos\theta a_i,\sin\theta b_j),\theta\in [\theta_0,\frac\pi2-\theta_0]\}\subset \gamma_{ij}.\ee For each $1\le i\le 3$, let $x_i\in P_1$ be a unit normal vector to $a_i$, and let $y_i\in P_2$ be a unit normal vector to $b_i$, so that the $a_i\wedge x_i$ are the same, $1\le i\le 3$, and $b_i\wedge y_i$ are the same, $1\le i\le 3$. 

Withoutloss of generality, suppose that 
\be a_2=-\frac 12 a_1+\frac{\sqrt 3}{2} x_1, a_3=-\frac 12 a_1-\frac{\sqrt 3}{2} x_1\mbox{ and 
}b_2=-\frac 12 b_1+\frac{\sqrt 3}{2} y_1,b_3=-\frac 12 b_1-\frac{\sqrt 3}{2} y_1.\ee

Then
\be \G_{ij}=(1-\eta)\g^0_{ij}\times \bar B_{x_i\wedge y_j}(0,R_1).\ee

Let $P_{ij}$ be the plane generated by $a_i\wg b_j$, and $Q_{ij}$ be the plane generated by $x_i\wg y_j$. Then $P_{ij}\perp Q_{ij}$, $\g_{ij}\subset P_{ij}$.

Also note that $I_{ij}$ is the intersection of $A_{1i}$ with $\G_{ij}$, and $J_{ij}$ is the intersection of $\G_{ij}$ with $A_{2j}$.

Note that by definition, 

\be Z\cap A_{1i}=\cup_{j=1}^3[o_{D_{1i}}, o_{I_{ij}}], 1\le i\le 3,\  
Z\cap A_{2j}=\cup_{j=1}^3[o_{D_{2j}}, o_{J_{ij}}], 1\le j\le 3,\ee
and 
\be Z\cap \G_{ij}=(1-\eta)\G^o_{ij}, 1\le i,j\le 3.\ee
Moreover, we have the essentially disjoint union
\be Z\cap \pa \cU=[\cup_{\a\in C}(Z\cap A_\a)]\cup [\cup_{1\le i,j\le 3}(Z\cap \G_{ij})].\ee

Let $l_{ij}$ be the intersection of the cone over $\g_{ij}$ with $\pa \cU$. Then it is the essentially disjoint union
\be l_{ij}=[o_{D_{1i}}, o_{I_{ij}}]\cup(1-\eta)\g^o_{ij}\cup [o_{D_{2j}}, o_{J_{ij}}],\ee
and we have another essentially disjoint union
\be Z\cap\pa \cU=\cup_{1\le i,j\le 3}l_{ij}.\ee

For each $\nu\le R_3$, let $l_{ij}^\nu=l_{ij}\bs B((o_{D_{1i}}, \nu)\cup B(o_{D_{2j}},\nu)$. Then $l_{ij}^\nu$ is the union of the arc $[(1-\eta)\g_{ij}]\cap\G_{ij}$ and two segments contained in $[o_{D_{1i}},o_{I_{ij}}]$ and $[o_{D_{2j}}, o_{J_{ij}}]$. In particular, $l_{ij}^0=l_{ij}$, and $l_{ij}^{R_3}=[(1-\eta)\g_{ij}]\cap\G_{ij}$. Let $\ml_{ij}^\nu$ denote the cone in $\bar\cU$ over $l_{ij}^\nu$. Set $L_{ij}^\nu=l_{ij}^\nu\times B_{Q_{ij}}(0,R_1)$, let $\t\mL_{ij}^\nu$ denote the cone over $L_{ij}^\nu$, and let $\mL_{ij}^\nu=\t\mL_{ij}^\nu\cap \bar\cU$. Let $L^\nu=\cup_{1\le i,j\le 3}L_{ij}^\nu $, $\mL^\nu=\cup_{1\le i,j\le 3}\mL_{ij}^\nu $, and $\t\mL^\nu=\cup_{1\le i,j\le 3}\t\mL_{ij}^\nu$.

Moreover, for $0\le \nu\le R_3$, 
\be\G_{ij}=L_{ij}^{R_3}\subset L_{ij}^\nu\subset L_{ij}^0=\G_{ij}\cup A_j\cup A_l\subset \pa \cU\cap B(K\cap \pa\cU, R_1).\ee

Note that all the notations above depend on the given $\eta$.

Finally we have, for the cone $Z$, Definition \tb{3.5} becomes:

\begin{defn}Let $\eta<\eta_1$. Take all the notations above (which depend on $\eta$ obviously). Let $\d \in(0,R_1)$, $\nu\in [0,R_3]$, $L>0$. A closed set $G\in \pa\cU(Z,\eta)$ is called an $(\eta,\d,\nu, L)$-sliding boundary of $Z$, if it is a $\d$-sliding deformation of $\pa\cU(Z,\eta)\cap Z$ in $\bar\cU(Z,\eta)$, and it satisfies in addition, 
\be \mbox{for each }(i,j)\in J, G\cap L_{ij}^\nu\mbox{ is the graph of an }L\mbox{-Lipschitz map from }l_{ij}^\nu\mbox{ to }B_{Q_{ij}}(0,R_1).\ee

\end{defn}

\section{Simplification for competitors}

We now begin to prove the $\Z_2$-topological sliding stability for $Z$. From now on, we fix an $\eta<\eta_1$. Let $\cU=\cU(\eta)$. Take all the notations as in Subsection \tb{3.2}.

In this section we will do some simplifications on the $\Z_2$-topological sliding competitors for $Z$ in the following two propositions and corollary. 

\begin{defn}For $0\le k\le 4$, a closed subset $F$ of $\bar \cU$ is said to be $k$-regular in $\bar \cU$ if there exists a finite smooth triangulation of $\bar \cU$ such that $F$ is the support of a $k$-simplicial sub-complex of this triangulation. \end{defn}

\begin{pro}For any $\eta<\eta_1$, let $\cU$ denote $\cU(\eta)$. Then for any $\d\in(0,R_1)$, $0<\nu<\nu'<R_3$, and $L>0$, 
\be \begin{split}&\inf\{\H^2(F):F\mbox{ is a }(\eta,\d,\nu, L)-\Z_2-\mbox{topological sliding competitor for }Z\}\\
\ge &\inf\{\H^2(F):F\mbox{ is a 2-regular }(\eta,\d,\nu', L)-\Z_2-\mbox{topological sliding competitor for }Z\}.\end{split}\ee
\end{pro}

\begin{pro}If $F$ is a 2-regular $(\eta,\d,\nu, L)$-$\Z_2$-topological sliding competitor  for $Z$, then there exists a 2-regular $\Z_2$-topological competitor $E$ for $Z$ in $\cU$, and an $\d$-sliding deformation $\varphi_t,1\le t\le 1$ in $\bar\cU$, so that $\varphi_1(E\cap\bar\cU)\subset F$. Moreover, $\varphi_1(E\cap\bar\cU)$ is also a 2-regular $(\eta,\d,\nu, L)$-$\Z_2$-topological sliding competitor for $Z$.
\end{pro}

Now set 
\be \begin{split}\F(\eta,\d,\nu, L):=\{&F=\varphi_1(E\cap \bar\cU)\mbox{ is a }(\eta,\d,\nu, L)-\Z_2-\mbox{topological sliding competitor for }\\
&Z\mbox{, where }E\mbox{ is a 2-regular-}\Z_2-\mbox{topological competitor for }Z\mbox{ in }\cU\mbox{, and }\\
&\{\varphi_t,0\le t\le 1\}\mbox{ is a }\d\mbox{-sliding deformation in }\bar\cU\}.\end{split}\ee

Then we have

\begin{cor}For any $\eta<\eta_1$, let $\cU$ denote $\cU(\eta)$. Then for any $\d\in(0,R_1)$, $\nu\in (0,R_3)$, and $L>0$, 
\be \begin{split}\inf\{\H^2(F):&\ F\mbox{ is a }(\eta,\d,\nu, L)-\Z_2-\mbox{topological sliding competitor for }Z\}\\
\ge \inf\{\H^2(F):&\ F\in \F(\eta,\d,\nu', L)\}, \forall \nu'\in (\nu, R_3).\end{split}\ee
\end{cor}

\nd This is a direct corollary of Propositions \tb{4.2 and 4.3}.\qed

After Corollary \tb{4.4}, to prove the $(\eta,\d,\nu, L)-\Z_2-$topological stability for $Z$, it is enough to look at the classes $\F(\eta,\d,\nu, L)$ of competitors that admits good regularities.

\smallskip

The rest of this section will be devoted to the proofs of Propositions \tb{4.2 and 4.3}. The proofs are somehow technical, while the conclusions are not too surprising. So for readers who would first like to know the main clue of the proof of the stabilities of $Y\times Y$, they can admit Corollary \tb{4.4} and jump directly to the next section.

\medskip

Let us first give a technical lemma:

\begin{lem}Let $K$ be a 2-dimensional minimal cone in $\R^n$, let $\eta\in (0,\eta_1(K))$, $\d\in (0,R_1)$. Denote by $\cU=\cU(K,\eta)$. Suppose $\psi_1:\bar \cU\to \bar \cU$ is a Lipschitz map satisfying that $\psi_1(\pa\cU)\subset \pa\cU$ and $|\psi_1(x)-x|<\d$ for $x\in K\cap \pa\cU$. Then we can extend $\psi_1$ to a sliding deformation in $\cU$, so that for any closed set $E\subset\R^n$ with $E\bs \cU=K\bs\cU$, $\psi_1(E\cap\bar\cU)$ is a $\d$-sliding deformation of $E$ in $\bar\cU$. That is, we can extend $\psi_1$ to a family of Lipschitz maps $\psi_t,0\le t\le 1$ from $\bar \cU\to \bar \cU$ which satisfies the conditions in Definition \tb{3.2} for $U=\cU$ and replacing $\varphi_t$ by $\psi_t$. Moreover, we can ask that 
\be\psi_t(x)=\frac{(1-t)x+t\psi_1(x)}{r_{(1-t)x+t\psi_1(x)}}, \forall x\in \pa\cU,\ee
where for any $x\in \R^n$, $r_x$ denote the real number such that $x\in r_x\pa\cU$.
\end{lem}

\nd This is a direct corollary of Remark \tb{3.3}, \cite{stablePYT} Lemma \tb{4.3} and its proof.\qed

\noindent\textbf{Proof of Proposition \tb{4.2}.} It is enough to prove that, for any $(\eta,\d,\nu, L)-\Z_2-$topological sliding competitor $F$ for $Z$, and any $\e>0$, there exists a 2-regular $(\eta,\d,\nu', L)-\Z_2-$topological sliding competitor $F'$ for $Z$, such that 
\be \H^2(F')<\H^2(F)+\e.\ee

So take any $(\eta,\d,\nu, L)-\Z_2-$topological sliding competitor $F$ for $Z$, and fix any $\e>0$. Let $\e_0<\e$ and $\a<10^{-2}$ be very small, to be decided later.

Since $F$ is a $(\eta,\d)-\Z_2-$topological sliding competitor, there exists a $\Z_2$-topological competitor $E$ of $Z$ in $\cU$, and a $\d$-sliding deformation $\varphi_t$ in $\bar\cU$, such that $F=\varphi_1(E\cap\bar\cU)$. By Lemma \tb{4.5}, we can suppose that for $x\in \pa\cU$, we have $\varphi_t(x)=\frac{(1-t)x+t\varphi_1(x)}{r_{(1-t)x+t\varphi_1(x)}}$, the ''projection'' of the segment $[x,\varphi_1(x)]$ out to $\pa\cU$.

Then it is clear that $|x-\varphi_t(x)|<|x-\varphi_1(x)|<\d$ and $\varphi_t(x)\in \pa\cU$ for $x\in E\cap \pa\cU$ and all $ t<1$. Let $\d_0=\sup_{x\in E\cap \pa\cU}|\varphi_1(x)-x|$. Since the function $|\varphi_1(x)-x|<\d$ for all $x\in E\cap \pa\cU$, and $E\cap \pa \cU$ is compact, we know that $\d_0<\d$. 

Now since $F$ is a $(\eta,\d,\nu, L)-\Z_2-$topological sliding competitor, by definition, for $1\le i,j\le 3$, $F\cap L_{ij}^\nu$ is the graph of a $L$-Lipschitz map $r_{ij}$ from $l_{ij}^\nu$ to $B_{Q_{ij}}(0,\d)$.

For each pair of $(i,j)$, let $a_{ij}^\nu$ denote the endpoint of $l_{ij}^\nu$ that is closer to $a_i$, and $b_{ij}^\nu$ denote the endpoint of $l_{ij}^\nu$ closer to $b_j$. Let $u_{ij}$ be a smooth $L$-Lipschitz map from $l_{ij}^\nu$ to $B_{Q_{ij}}(0,\d)$, so that 
\be u_{ij}(a_{ij}^\nu)=r_{ij}(a_{ij}^\nu), u_{ij}(b_{ij}^\nu)=r_{ij}(b_{ij}^\nu),\ee
and
\be \sup_{z\in l_{ij}^\nu}||r_{ij}(z)-u_{ij}(z)||<\min\{\e_0,\frac 12(\d-\d_0)\}.\ee
Then the graph of $u_{ij}$ is contained in $l_{ij}^\nu\times B_{Q_{ij}}(0,\d)\subset L_{ij}^\nu$.

Define $f_{ij}: F\cap L_{ij}^\nu\to L_{ij}^\nu$: for each $x\in F\cap L_{ij}^\nu$, by definition, there exists $z\in l_{ij}^\nu$ so that $x=(z, r_{ij}(z))\in l_{ij}^\nu\times B_{Q_{ij}}(0,\d)\subset L_{ij}^\nu$. Set $f_{ij}(x)=(z, u_{ij}(z))\in l_{ij}^\nu\times B_{Q_{ij}}(0,\d)$. Then $f_{ij}$ is $\sqrt{1+L^2}$-Lipschitz, and for any $x\in F\cap L_{ij}^\nu$, the segment $[x,f_{ij}(x)]\in L_{ij}^\nu$. 

Let $f: \varphi_1(E\cap \pa \cU)\to \pa \cU$: 
\be f(x)=\left\{\begin{array}{rcl}f_{ij}(x),&\  &x\in L_{ij}^\nu;\\
x,&\  &x\in A\bs \t\mL^\nu.
\end{array}\right.\ee
Then $f$ is $\sqrt{1+L^2}$-Lipschitz, and for any $x\in E\cap \pa \cU$, the segment $[x,f(x)]\subset \pa\cU$, that is, for any $t\in [0,1]$, $(1-t)x+tf(x)\in \pa\cU$. Moreover, $|f(x)-x|\le \min\{\e_0,\frac 12(\d-\d_0)\}$ for all $x\in \varphi_1(E)\cap\pa\cU$.

We define $\psi_1:(1+9\a)\bar\cU\to (1+9\a)\bar\cU$, so that 
\be \psi_1(x)=\left\{\begin{array}{rcl}
\varphi_1(x)&,\ if\ &x\in \bar\cU;\\
(1+3t\a)[(1-t)\varphi_1(\frac{x}{1+3t\a})+tf\circ\varphi_1(\frac{x}{1+3t\a})]&,\ if\ &x\in (1+3t\a)\pa\cU, t\in [0,1];\\
r_x[f\circ\varphi_1(\frac{x}{r_x})]&,\ if\ &x\in (1+9\a)\bar\cU\bs (1+3\a)\bar\cU.
\end{array}\right.\ee

Then $\psi_1=\varphi_1$ in $\bar\cU$, $\psi_1(E\cap (1+t\a)\pa\cU)\subset 1+t\a)\pa\cU, \forall 0\le t\le 9$, and since $E\cap(1+9\a)\bar\cU\bs \bar\cU=K\cap(1+9\a)\bar\cU\bs \bar\cU$, which is the cone over $K\cap\pa\cU$, hence $\psi_1(E\cap(1+9\a)\bar\cU)\bs (1+3\a)\bar\cU$ coincides with the cone over $f\circ\varphi_1(K\cap \pa\cU)$. In particular, for any $x\in E\cap (1+9\a)\pa\cU=K\cap (1+9\a)\pa\cU$, we know that $\frac{x}{1+9\a}\in \pa\cU$, and hence 
\be \begin{split}|\psi_1(x)-x|& \le (1+9\a)|f\circ\varphi_1(\frac{x}{1+9\a})-\frac{x}{1+9\a}|\\
&\le (1+9\a)(|f\circ\varphi_1(\frac{x}{1+9\a})-\varphi_1(\frac{x}{1+9\a})|+|\varphi_1(\frac{x}{1+9\a})-\frac{x}{1+9\a}|)\\
&\le (1+9\a)[\d_0+\frac12(\d-\d_0)], \forall x\in E\cap (1+9\a)\pa\cU.
\end{split}\ee

Set $F_1=\psi_1(E\cap(1+9\a)\bar\cU)$, then a simple calculation yields that
\be \H^2(F_1)=\H^2(F)+\H^2(F_1\bs F)\le \H^2(F)+(100\a+2\e_0)\sqrt{1+L^2}\H^2(\varphi_1(E\cap \pa\cU)).\ee

By Lemma \tb{4.5}, $F_1$ is a $(1+9\a)\d$ sliding deformation of $E$ in $(1+9\a)\bar\cU$. It satisfies that
\be \begin{split} F_1\cap\t\mL^\nu\cap (1+9\a)\bar\cU\bs (1+3\a)\bar\cU\mbox{ coincides with the cone over }\\
f\circ\varphi_1(K\cap \pa\cU)\mbox{, and hence is 2-regular}.\end{split}\ee
In particular, if we set  $G_0:=F_1\cap \t\mL^\nu\cap (1+9\a)\bar\cU\bs (1+3\a)\cU$ and $G_1:=F_1\cap \t\mL^{\nu''}\cap (1+9\a)\bar\cU\bs (1+4\a)\cU$ with $\nu''=\frac{\nu+\nu'}{2}$, then for each pair $i,j$, 
\be G_0\mbox{ coincides with the cone over the graph of }u_{ij}\mbox{ in }\mL_{ij}^\nu\bs (1+3\a)\cU,\ee and 
\be \mbox{ and }G_1\mbox{ coincides with the cone over the graph of }u_{ij}\mbox{ in }\mL^{\nu''}_{ij}\bs (1+4\a)\cU.\ee Hence $G_1$ is a piecewise smooth Lipschitz surface with boundary $\pa G_1$. Note that $\pa G_1\cap \pa \cU=\cup_{1\le i,j\le 3}G(u_{ij})$, where $G(u_{ij})$ denotes the graph of $u_{ij}$, the endpoints of $G(u_{ij})$ being $a_{ij}^\nu, b_{ij}^\nu$. 

Since $G_0$ coincides with the cone over the union of the graphs of the $L$-Lipschitz maps $u_{ij}, 1\le i,j\le 3$, there exists a constant $M=M(L)>0$ that depends only on $L$, so that 
\be\H^2(G_0\cap B(x,r))\ge Mr^2, \forall x\in G_0, r\in (0, \a).\ee

 Let $W_1=\pa G_1\cap (1+9\a)\cU$, the part of $\pa G_1$ in $(1+9\a)\cU$. Then it is a union of disjoint piecewise smooth curves.

Take $\e_1>0$ so that $100\e_1<\min\{\nu'-\nu, R_1-\d,\a\}$, and $\H^2(G_0\cap \bar B(W_1, 10\e_1))<\e_0$. Let $G_3=(F_1\bs G_1)\cup \overline{W_1}$. Then $F_1=G_1\cup G_3$ and $G_1\cap G_3=\bar W_1$.

Set $V=\t\mL^{\nu'}\bs (1+5\a)\cU$. Let $G_2=G_1\cap V$.

We claim that
\be d(V, G_3)\ge 10\e_1.\ee

By definition $G_3\subset (1+9\a)\bar\cU$ and does not meet $\t\mL^\nu\cap (1+9\a)\bar\cU\bs (1+4\a)\bar\cU$. Let $G_3^1=G_3\cap (1+4\a)\bar \cU$, and $G_3^2=G_3\cap (1+9\a)\bar\cU\bs (1+4\a)\bar\cU$. Then $G_3=G_3^1\cup G_3^2$. By definition of $V$, we know that $V\cap (1+5\a)\cU=\emptyset$, hence $d(V, (1+4\a)\bar \cU)>\frac12 \a>10\e_1$, and since $G_3^1\subset (1+4\a)\bar \cU$, we know that $d(V, G_3^1)>10\e_1$; 

On the other hand, by definition, we know that both $G_3^2$ and $V$ are both part of cones: $G_3^2$ is contained in the cone over $A\bs L^\nu$, and $V$ is contained in the cone over $L^{\nu'}$. Note that both  $A\bs L^\nu$ and $L^{\nu'}$ are parts of $\pa\cU$, $d(A\bs L^\nu,L^{\nu'})\ge \frac 12(\nu-\nu')\}>50\e_1$. Hence outside $\bar\cU$, the distance between the cone over $A\bs L^\nu$ and over $L^{\nu'}$ is larger than $10\e_1$. Therefore $d(G_3^2, V)>10\e_1$.

Altogether we get Claim \tb{(4.16)}.

Set $U=B(G_3, 5\e_1)$. Then $U$ is open and $d(\bar U, G_2)\ge 5\e_1$. We apply \cite{Fv} Theorem 4.3.17  to $U$, $h=1$, $E=G_3$, and get a 4-dimensional polyhedral complex $\cS$, and a Lipschitz map $\phi: U\to U$ (where $|\cS|=\cup\{\sigma:\sigma\in \cS\}$ is the support of $\cS$) so that:

$1^\circ$ The maximal diameter of simplices in $\cS$ is less than $\e_1$;

$2^\circ$ $G_3\subset |\cS|^\circ$, and $\phi(F_1\cap |\cS|)$ is a union of some 2-faces of $\cS$;

$3^\circ$ For each $\sigma\in \cS$ we have $\phi(\sigma)\subset \sigma$, and $\phi=Id$ on the 2-skeleton of $\cS$.  Moreover, $\H^d(\phi(E\cap \sigma))\le M_0\H^d(E\cap \sigma), \forall\sigma\in \cS$, where $M_0$ is a constant (that only depend on the ambient dimension, which is 4 in our situation).

$4^\circ$ The ''surface'' of $\cS$ is a 3-dimensional dyadic complex: that is, there exists a dyadic complex $\cT$ so that
\be \forall \sigma\in \cS, \sigma\cap\pa|\cS|\ne\emptyset\Rightarrow\exists \sigma'\in \cT\mbox{ so that }\sigma\cap \pa|\cS|=\sigma'\cap\pa|\cT|.\ee

$5^\circ$ $\phi|_{|\cS_\pa|}$ is a Federer-Fleming projection of $F_1\cap|\cS_\pa|$ in $\cS_\pa$, where $\cS_\pa=\{\sigma\in \cS: \sigma\cap \pa|\cS|\ne\emptyset\}$.

$6^\circ$ $\H^2(\phi(F_1\cap U))\le \H^2(F_1\cap U)+\e_0$.

As a remark, here the dyadic complex can be chosen with arbitrary fixed coordinates of $\R^4$. Hence we can ask that 

$6^\circ$ Suppose that the dyadic cubes $\cT$ in $3^\circ$ are of sidelength $2^{-N}$. Let $\cT'$ be the set of all dyadic cubes $\sigma$ of sidelength $2^{-N}$ so that $\sigma\cap \pa|\cS|\ne\emptyset$ and $\sigma\cap |\cS|^\circ=\emptyset$. Then $G_1$ is transversal to all the $d$-faces of $\cT'$ for $d\le 3$.

Since $\phi$ is a Federer-Fleming projection of $F_1\cap |\cS_\pa|$ in $\cS_\pa$, the restriction of $\phi$ on $F_1\cap |\cS_\pa|$ is in fact a composition of two maps $\phi_1: F_1\cap |\cS_\pa|\to |\cS_\pa^3|$ and $\phi_2:\phi_1(F_1\cap |\cS_\pa^3|)\to |\cS^2|$, so that
\be \phi_1(\sigma)\subset\sigma, \forall \sigma\in \cS_\pa\mbox{ and }\phi_2(\sigma)\subset\sigma, \forall \sigma\in \cS_\pa^3,\ee
\be {\phi_1}|_{F_1\cap|\cS_\pa^3|}=Id.\ee

We define a new map $\phi': [F_1\cap |\cS|]\cup\pa|\cS|\to |\cS^2|\cup \pa|\cS|$, so that:
\be\phi'(x)=\left\{\begin{array}{rcl}\phi(x)&,\ &\mbox{if }\x\in |\cS|\bs |\cS_\pa|;\\
\phi(x)=\phi_2\circ\phi_1(x)&,\ &\mbox{if }\phi_1(x)\in |\cS_\pa^3|\cap |\cS|^\circ;\\
\phi_1(x)&,\ &\mbox{if }\phi_1(x)\in |\cS_\pa^3|\cap\pa|\cS|;\\
x&,\ &\mbox{if }x\in \pa|\cS|.
\end{array}
\right.\ee
Then $\phi'$ is Lipchitz and $\phi'(F_1\cap |\cS|)\bs \pa|\cS|$ is the union of some 2-faces of $\cS$.

we extend it to a Lipschitz map, still denoted by $\phi':\R^4\to \R^4$, so that
\be \phi'(\sigma)\subset\sigma, \forall \sigma\in \cS\mbox{, and }\phi'|_{\R^4\bs |\cS|^\circ}=Id.\ee 

By definition of $\phi'$, we know that $\phi'(F_1\cap |\cS|)\bs \pa|\cS|\subset \phi(F_1\cap |\cS|)\subset \phi(F_1\cap U)$, and hence by $5^\circ$, 
\be \H^2(\phi'(F_1\cap |\cS|)\bs \pa|\cS|)\le \H^2(\phi(F_1\cap U))\le \H^2(F_1\cap U)+\e_0.\ee

Suppose that the dyadic cubes $\cT$ in $3^\circ$ are of sidelength $2^{-N}$. Then $2^{-N}<\e_1$. Let $\cT'$ be the set of all dyadic cubes $\sigma$ of sidelength $2^{-N}$ so that $\sigma\cap U\ne\emptyset$ and $\sigma\cap |\cS|^\circ=\emptyset$. Then $|\cT'|\subset B(U, 2^{-N+1})\subset B(U, 2\e_1)$. Set $\cS'=\cS\cup \cT'$. 

By definition, $|\cS|\subset U$, and hence $|\cS'|\subset B(U, 2\e_1)$. By definition of $U$, we know that $|\cS'|\subset B(G_3, 5\e_1+2\e_1)=B(G_3, 7\e_1)$. By Claim \tb{(4.16)}, $|\cS'|\cap V=\emptyset$. Since $G_3\subset |\cS|^\circ\subset |\cS'|^\circ$, and $|\cS'|\cap V=\emptyset$, we know that $\pa|\cS'|\cap F_1=\pa|\cS'|\cap G_1\bs V$.

By definition of $\cS'$, we know that $\pa|\cS'|$ is the union of 3-dimensional dyadic faces of sidelength $2^{-N}$. By $6^\circ$, we know that $G_1\bs |\cS'|$ is a smooth 2-surface with piecewise smooth boundary. Hence $G_1\bs |\cS'|$ is 2-regular.

Now let $\phi'':\R^4\to \R^4$ be a Federer-Fleming projection of $\phi'(F_1)$ in $|\cS'|$. More precisely:
\be\phi''|_{|\cS'|^C}=Id\mbox{, and }\phi''|_{|\cS'^2|\cup \pa|\cS'|}=Id\ee
\be\phi''\circ\phi'(F_1)\cap |\cS'|\subset |\cS'^2|\cup \pa|\cS'|,\ee
\be\phi''(\sigma)\subset\sigma, \forall \sigma\in \cS',\ee
\be \phi''\circ\phi'(F_1\cap|\cS'|)\bs \pa|\cS'|\mbox{ is a union of 2-faces of }\cS'\ee
and there exists a constant $C_0$ that only depend on the dimension of $\R^4$, so that for any $\sigma\in \cT'$, 
\be \H^2(\phi''(\phi'(F_1)\cap \sigma))\le C_0\H^2(\phi'(F_1)\cap \sigma).\ee
The existence of $C_0$ is because $\cT'$ is a dyadic complex.

Also note that $\phi''|_{\phi'(F_1)\cap|\cS|^\circ}=Id$, because $\phi'(F_1)\cap|\cS|^\circ$ is a union of 2-faces of $\cS$.

We will prove that
\be \phi''\circ\phi'(F_1)\mbox{ is 2-regular},\ee
and 
\be \H^2(\phi''\circ\phi'(F_1))\le \H^2(F_1)+C_1\e_0,\ee
where $C_1$ is a constant that only depends on $L$.

For \tb{(4.28)}, we know that $\phi''\circ\phi'(F_1)\bs |\cS'|=F_1\bs |\cS'|=G_1\bs |\cS'|$ is 2-regular, and by 
\tb{(4.26)}, $\phi''\circ\phi'(F_1)\cap |\cS'|^\circ$ is 2-regular, hence it is enough to prove that $\phi''\circ\phi'(F_1)\cap \pa|\cS'|$ is 2-regular. By the process of a Federer-Fleming projeciton in a polyhedral complex with boundary, if we denote by $\cS''=\{\sigma\in \cS':\sigma\cap \pa|\cS'|\ne\emptyset\}\subset \cT'
$, then $\phi''\circ\phi'(F_1)\cap \pa|\cS'|$ comes from the image of $\phi''$ of $\phi'(F_1)\cap \cS''$. Moreover, for each $\sigma\in \cS''$, $\phi''\circ\phi'(F_1)\cap\sigma\cap \pa|\cS'|$ is the intersection of $\sigma\cap \pa|\cS'|$ with the image of $\phi'(F_1)\cap\sigma^\circ$ under a radial projection $\pi_{\sigma}: \sigma\to \pa\sigma$.

Take any $\sigma\in \cS''$. Then $\sigma\in \cT'$, and hence $\sigma^\circ\cap |\cS|=\emptyset$. By definition of $\phi'$, we know that $\phi'(|\cS|)\subset |\cS|$, and $\phi'|_{|\cS|^C}=Id$, hence $\phi(F_1)\cap \sigma^\circ=F_1\cap \sigma^\circ$. Since $G_3\subset |\cS|^\circ$, $\sigma^\circ\cap |\cS|=\emptyset$, and $F_1=G_1\cup G_3$, hence $F_1\cap \sigma^\circ=G_1\cap \sigma^\circ$. Note that $G_1$ is a 2-dimensional smooth surface with smooth boundary, and is transversal to $\cT'$, hence is transversal to $\sigma$. As a result, $G_1\cap \sigma^\circ$ is a smooth 2-surface with piecewise smooth boundary. Hence $\pi_\sigma(G_1\cap \sigma^\circ)$ is a piecewise smooth 2-surface with piecewise smooth boundary, therefore so is $\pi_\sigma(G_1\cap \sigma^\circ)\cap [\sigma\cap \pa|\cS'|]$. Altogether we have that
\be \phi''\circ\phi'(F_1)\cap\sigma\cap \pa|\cS'|=\pi_\sigma(G_1\cap \sigma^\circ)\cap [\sigma\cap \pa|\cS'|]\ee is a piecewise smooth 2-surfacw with piecewise smooth boundary, for any $\sigma\in \cS''$.

Hence $\phi''\circ\phi'(F_1)\cap \pa|\cS'|$ is 2-regular, and we get \tb{(4.28)}.

Next let us look at the measure of $\phi''\circ\phi'(F_1)$. Note that $\phi''\circ\phi'|_{|\cS'|^C}=Id$, hence to prove \tb{(4.29)}, it is enough to look at $\phi''\circ\phi'(F_1\cap |\cS'|)$. By definition, 
\be\phi'(F_1\cap |\cS'|)=\phi'(F_1\cap |\cS|)\cup \phi'(F_1\cap |\cT'|)=[\phi'(F_1\cap |\cS|)\bs \pa|\cS|]\cup[\phi'(F_1\cap |\cS|)\cap\pa|\cS|]\cup \phi'(F_1\cap |\cT'|),\ee
and hence
\be \H^2(\phi''\circ\phi'(F_1\cap |\cS'|))\le \H^2(\phi''(\phi'(F_1\cap |\cS|)\bs \pa|\cS|))+\H^2(\phi''(\phi'(F_1\cap |\cS|)\cap\pa|\cS|))+\H^2(\phi''(\phi'(F_1\cap |\cT'|))).\ee
For the first term of the right-hand-side of \tb{(4.32)}, we know that $\phi'(F_1\cap |\cS|)\bs \pa|\cS|$ is the union of some 2-faces of $\cS$, hence  by \tb{(4.23)}, $\phi''(\phi'(F_1\cap |\cS|)\bs \pa|\cS|)=\phi'(F_1\cap |\cS|)\bs \pa|\cS|$, and therefore by \tb{(4.22)}, 
\be \H^2(\phi''(\phi'(F_1\cap |\cS|)\bs \pa|\cS|))=\H^2(\phi'(F_1\cap |\cS|)\bs \pa|\cS|)\le \H^2(F_1\cap U)+\e_0.\ee

For the second term of the right-hand-side of \tb{(4.32)}, we know that
\be \phi'(F_1\cap |\cS|)\cap\pa|\cS|\subset \pa|\cS|\subset |\cT'|\cap|\cS|,\ee
and hence by \tb{(4.25)} we know that $\phi''\circ\phi'(F_1\cap |\cS|)\cap\pa|\cS|\subset |\cT'|\cap|\cS|\subset (\pa|\cS'|)^C$, hence by \tb{(4.25) and (4.26)}, $\phi''(\phi'(F_1\cap |\cS|)\cap\pa|\cS|)$ is a union of 2-faces of $\cT'$, that is,
\be \phi''(\phi'(F_1\cap |\cS|)\cap\pa|\cS|)\subset |\cT'^2|.\ee

By \tb{(4.21) and (4.25)}, we know that 
\be \phi''(\phi'(F_1\cap |\cS|)\cap\pa|\cS|)\subset |\cT''^2|,\ee where $\cT''=\{\sigma\in \cT':\sigma\cap F_1\ne\emptyset\}$. Since $G_3\cap |\cT'|=\emptyset$, hence $|\cT'|\cap F_1=|\cT'|\cap G_1$. On the other hand, we know that $|\cT'|\subset B(U,2\e_1)$, hence $|\cT''|\cap F_1\subset|\cT'|\cap F_1\subset G_1\cap B(U,2\e_1)$. By definition of $U$, we know that
\be |\cT''|\cap F_1\subset G_1\cap B(G_3, 7\e_1)\subset B(W_1, 7\e_1).\ee
Also note that $B(W_1, 7\e_1)\cap F_1\subset G_0$, hence for any $\sigma\in \cT''$, we know that $\sigma\cap F_1=\sigma\cap G_0$, and $\sigma\subset B(W_1, 7\e_1)$. Now by \tb{(4.15)} and the fact that $\H^2(G_0\cap B(W_1, 7\e_1))\le \e_0$, we know that 
\be \H^2(|\cT''^2|)\le C(M)\e_0,\ee where $C(M)$ is a constant that only depends on $M$, and hence on $L$. Combine with \tb{(4.36)}, we get
\be \H^2(\phi''(\phi'(F_1\cap |\cS|)\cap\pa|\cS|))\le C(M)\e_0.\ee

Finally for the third term of the right-hand-side of \tb{(4.32)}, by \tb{(4.27)}, we have
\be \H^2(\phi''(\phi'(F_1\cap |\cT'|)))\le \sum_{\sigma\in \cT'}\H^2(\phi''(\phi'(F_1\cap \sigma)))\le C_0\sum _{\sigma\in \cT'}\H^2(\phi'(F_1\cap \sigma)).\ee
By \tb{(4.21)}, we know that $\phi'(F_1\cap |\cT'|)=F_1\cap |\cT'|$, hence for any $\sigma\in \cT'$ we have $\phi'(F_1\cap \sigma)=F_1\cap \sigma$, and therefore
\be \H^2(\phi''(\phi'(F_1\cap |\cT'|)))\le C_0\sum_{\sigma\in \cT'}\H^2(F_1\cap \sigma)\le C_4C_0\H^2(F_1\cap|\cT'|),\ee
here $C_4=\sup_{x\in \R^4}\sum_{\sigma\in \Delta_N}1_{\sigma}(x)$ is a constant that depend only on the dimension of $\R^4$ ($\Delta_N$ being the family of dyadic cubes of sidelength $2^{-N}$).

Again since $|\cT'|\cap F_1\subset G_1\cap B(U,2\e_1)\subset G_1\cap B(G_3, 7\e_1)\subset G_1\cap B(W_1, 7\e_1)$, we have $\H^2(|\cT'|\cap F_1)\le \H^2(G_1\cap B(W_1, 7\e_1))\le \H^2(G_0\cap B(W_1,10\e_1))<\e_0$. Combine with \tb{(4.41)} we get
\be \H^2(\phi''(\phi'(F_1\cap |\cT'|)))\le C_4C_0\e_0.\ee

Summing up \tb{(4.32), (4.33), (4.39) and (4.42)} we get
\be \H^2(\phi''\circ\phi'(F_1\cap |\cS'|)\le \H^2(F_1\cap U)+(1+C(M)+C_0C_4))\e_0.\ee

Now since $\phi''\circ\phi'|_{|\cS'|^C}=Id$, we know that $\phi''\circ\phi'(F_1\bs |\cS'|)=F_1\bs |\cS'|$, and hence
\be\begin{split} \H^2(\phi''\phi'(F_1))&\le \H^2(\phi''\phi'(F_1\cap |\cS'|))+\H^2(\phi''\phi'(F_1\bs |\cS'|))\\
&=\H^2(\phi''\phi'(F_1\cap |\cS'|))+\H^2(F_1\bs |\cS'|)\\
&\le \H^2(F_1\cap U)+(1+C(M)+C_0C_4))\e_0+\H^2(F_1\bs |\cS'|)\\
&=\H^2(F_1\cap U)+(1+C(M)+C_0C_4))\e_0+\H^2(F_1\bs |\cS'|)+\H^2(F_1\bs U)-\H^2(F_1\bs U)\\
&=\H^2(F_1)+(1+C(M)+C_0C_4))\e_0+[\H^2(F_1\bs |\cS'|)-\H^2(F_1\bs U)].
\end{split}\ee

Since $G_3\subset F_1\cap |\cS'|$ and $U=B(G_3, 5\e_1)$, we have
\be \begin{split}&\H^2(F_1\bs |\cS'|)-\H^2(F_1\bs U)\le \H^2(F_1\bs G_3)-\H^2(F_1\bs B(G_3,5\e_1))\\
&=\H^2(G_1)-\H^2(G_1\bs B(G_3,5\e_1))=\H^2(G_1\cap B(W_1,5\e_1))\\
&\le \H^2(G_0\cap B(W_1,10\e_1))<\e_0.
\end{split}\ee

Combine with \tb{(4.44)}, and set $C_1=2+C(M)+C_0C_4$, we get
\be \H^2(\phi''\phi'(F_1))\le \H^2(F_1)+(2+C(M)+C_0C_4))\e_0=C_1\e_0,\ee
which yields \tb{(4.29)}.

Set $\pi:\R^4\to (1+7\a)\bar\cU$:
\be \pi(x)=\left\{\begin{array}{rcl}\frac{1}{1+7\a}(x)&,\ if\ &x\in \R^4\bs(1+7\a)\bar\cU;\\
x&,\ if\ &x\in (1+7\a)\bar\cU.
\end{array}\right.\ee

Then $\pi$ is $C(\eta)$-Lipschitz, where $C(\eta)$ is a constant that depends only on $\eta$. Let $F_2=\pi\circ\phi''\circ\phi'(F_1)$. Then $F_2$ is 2-regular. Since $\pi(x)=x$ in $(1+7\a)\bar\cU$, we know that
\be \begin{split}\H^2(F_2)&\le\H^2(\phi''\phi'(F_1)\cap (1+7\a)\bar\cU)+\H^2(\pi (\phi''\phi'(F_1)\bs (1+7\a)\bar\cU))\\
&\le \H^2(\phi''\phi'(F_1))+C(\eta)^2\H^2(\phi''\phi'(F_1)\bs (1+7\a)\bar\cU)\\
&\le  \H^2(F_1)+C_1\e_0+C(\eta)^2\H^2(\phi''\phi'(F_1)\bs (1+7\a)\bar\cU).
\end{split}\ee 
By \tb{(4.21), (4.23) and (4.25)}, we know that
\be \phi''\phi'(F_1)\bs (1+7\a)\bar\cU\subset \phi''\phi'((F_1)\bs (1+6\a)\cU),\ee
hence by \tb{(4.27)} and property \tb{$3^\circ$} of $\phi$,
\be \H^2(\phi''\phi'(F_1)\bs (1+7\a)\bar\cU)\le \H^2(\phi''\phi'((F_1)\bs (1+6\a)\cU))\le M_0C_0\H^2(F_1\bs (1+6\a)\cU).\ee

Since $F_1$ coincide with the cone over $f\circ\varphi_1(E\cap \pa\cU)$ in $(1+9\a)\bar\cU\bs (1+3\a)\cU$, we know that
\be \H^2(F_1\bs (1+6\a)\cU)\le 3\a(1+9\a)\H^2(f\circ\varphi_1(E\cap \pa\cU))\le 3\a(1+9\a)\sqrt{1+L^2}\H^2(\varphi_1(E\cap \pa\cU)).\ee

Combine with \tb{(4.48) and (4.50)} we have
\be \H^2(F_2)\le\H^2(F_1)+C_1\e_0+C(\eta)^2M_0C_03\a(1+9\a)\sqrt{1+L^2}\H^2(\varphi_1(E\cap \pa\cU)).\ee
Set $F'=\frac{1}{1+7\a}F_2$, then
\be \H^2(F')\le \frac{1}{1+7\a}[\H^2(F_1)+C_1\e_0+C(\eta)^2M_0C_03\a(1+9\a)\sqrt{1+L^2}\H^2(\varphi_1(E\cap \pa\cU))].\ee

By \tb{(4.11)}, we get
\be \begin{split}\H^2(F')\le &\frac{1}{1+7\a}[\H^2(F)+(100\a+2\e_0)\sqrt{1+L^2}\H^2(\varphi_1(E\cap \pa\cU))\\
&+C_1\e_0+C(\eta)^2M_0C_03\a(1+9\a)\sqrt{1+L^2}\H^2(\varphi_1(E\cap \pa\cU))]\\
\le &\H^2(F)+C_1\e_0+[(100\a+2\e_0)+C(\eta)^2M_0C_03\a(1+9\a)]\sqrt{1+L^2}\H^2(\varphi_1(E\cap \pa\cU)).
\end{split}\ee

We would like to prove that $F'$ is an $(\eta,\d,\nu',L)$-$\Z_2$-topological sliding competitor for $Z$. 

Let $E_1=\frac{1}{1+9\a}E$, then 
$E_1$ is a deformation of $E$ in $\cU$, hence by Proposition \tb{2.7}, $E_1$ is a $\Z_2$-topological competitor for $E$ in $\cU$. 

For every $r>0$, let $\d_r$ denote the map $\d_r(x)=rx, \forall x\in \R^4$. Then $g_1:=\d_\frac{1}{1+7\a}\circ\pi\circ\psi''\circ\psi'\circ\psi_1\circ\d_{1+9\a}$ is a Lipschitz map from $\bar\cU\to \bar\cU$. And $F'=g_1(E_1\cap\bar\cU)$.

Let us first prove that $g_1(E_1\cap\bar\cU)$ is a $\d$-sliding deformation of $E_1$ in $\bar\cU$. By Lemma \tb{4.5}, it is enough to prove that 
\be g_1(\pa\cU)\subset \pa\cU\ee and \be|g_1(x)-x|<\d\mbox{ for }x\in Z\cap \pa\cU.\ee

For \tb{(4.55)}, take any $x\in \pa\cU$. Then since $\psi_1:(1+9\a)\bar\cU\to (1+9\a)\bar\cU$ is a sliding deformation, we know that $\psi_1\circ\d_{1+9\a}(x)\in (1+9\a)\pa\cU$. Now by \tb{(4.21), (4.23) and (4.25)}, since the maximal diameters of simplices in $\cS, \cS'$ are less than $\e_1$, we know that 
\be |\psi''\circ\psi'\circ\psi_1\circ\d_{1+9\a}(x)-\psi_1\circ\d_{1+9\a}(x)|<2\e_1<\frac12 \a,\ee
and hence
\be\psi''\circ\psi'\circ\psi_1\circ\d_{1+9\a}(x)\in B((1+9\a)\pa\cU, \frac12 \a)\subset \R^4\bs (1+8\a)\bar\cU,\ee
which imples that $\pi\circ\psi''\circ\psi'\circ\psi_1\circ\d_{1+9\a}(x)\in \pa{(1+7\a)}\cU$, and hence 
\be g_1(x)=\d_\frac{1}{1+7\a}\circ\pi\circ\psi''\circ\psi'\circ\psi_1\circ\d_{1+9\a}(x)\in \pa\cU.\ee
This yields \tb{(4.55)};

For \tb{(4.56)}, take any $x\in Z\cap \pa\cU=E_1\cap \pa\cU$. Then $\psi_1\circ\d_{1+9\a}(x)\in F_1$, and by \tb{(4.10)}, $|\psi_1\circ\d_{1+9\a}(x)-\d_{1+9\a}(x)|<(1+9\a)[\d_0+\frac 12(\d-\d_0)]$. 

By \tb{(4.21), (4.23) and (4.25)}, we know that 
\be|\psi''\circ\psi'\circ\psi_1\circ\d_{1+9\a}(x)-\psi_1\circ\d_{1+9\a}(x)|<4\e_1,\ee
and \be\psi''\circ\psi'\circ\psi_1\circ\d_{1+9\a}(x)\in (1+10\a)\bar\cU\bs (1+8\a)\cU.\ee
By definition of $\pi$, \tb{(4.61)} tells that
\be |\pi\circ\psi''\circ\psi'\circ\psi_1\circ\d_{1+9\a}(x)-\psi''\circ\psi'\circ\psi_1\circ\d_{1+9\a}(x)|\le 3\a.\ee
Hence \tb{(4.10), (4.60) and (4.62)} yields
\be \begin{split}&|\pi\circ\psi''\circ\psi'\circ\psi_1\circ\d_{1+9\a}(x)-\d_{1+7\a}(x)|\\
\le &|\pi\circ\psi''\circ\psi'\circ\psi_1\circ\d_{1+9\a}(x)-\psi''\circ\psi'\circ\psi_1\circ\d_{1+9\a}(x)|\\
&+|\psi''\circ\psi'\circ\psi_1\circ\d_{1+9\a}(x)-\psi_1\circ\d_{1+9\a}(x)|\\
&+|\psi_1\circ\d_{1+9\a}(x)-\d_{1+9\a}(x)|+|\d_{1+9\a}(x)-\d_{1+7\a}(x)|\\
\le &(1+9\a)[\d_0+\frac 12(\d-\d_0)]+4\e_1+3\a+2\a\\
\le &(1+9\a)[\d_0+\frac 12(\d-\d_0)]+6\a.
\end{split}\ee

As consequence,
\be \begin{split}|g_1(x)-x|&=|\d_\frac{1}{1+7\a}\circ\pi\circ\psi''\circ\psi'\circ\psi_1\circ\d_{1+9\a}(x)-\d_\frac{1}{1+7\a}\circ\d_{1+7\a}(x)|\\
&=\frac{1}{1+7\a}|\pi\circ\psi''\circ\psi'\circ\psi_1\circ\d_{1+9\a}(x)-\d_{1+7\a}(x)|\\
&\le \frac{1}{1+7\a}\{(1+9\a)[\d_0+\frac 12(\d-\d_0)]+6\a\}.
\end{split}\ee

Now take $\a$ and $\e_0$ so that 
\be \frac{1}{1+7\a}\{(1+9\a)[\d_0+\frac 12(\d-\d_0)]+6\a\}<\d\ee
and
\be C_1\e_0+[(100\a+2\e_0)+C(\eta)^2M_0C_03\a(1+9\a)]\sqrt{1+L^2}\H^2(\varphi_1(E\cap \pa\cU))<\e.\ee

Then we get \tb{(4.56)}, and \tb{(4.66)} gives \tb{(4.5)}. \tb{(4.55) and (4.56)} yields that $F'=g_1(E_1\cap\bar\cU)$ is a $\d$-sliding deformation of $E_1$ in $\bar\cU$.

We still have to prove that $F'\cap \pa\cU=g_1(E_1\cap\bar\cU)\cap\pa\cU$ is a $(\eta,\d,\nu', L)$-sliding boundary of $Z$. That is, we have to show that 
\be \forall 1\le i,j\le 3, F'\cap L_{ij}^{\nu'}\mbox{ is the graph of an }L\mbox{-Lipschitz map from }l_{ij}^{\nu'}\mbox{ to }B_{Q_{ij}}(0, R_1).
\ee
Since $F'=\frac{1}{1+7\a}F_2$, hence \tb{(4.60)} is equivalent to say that $F_2\cap (1+7\a)L_{ij}^{\nu'}$ is the graph of an $L$-Lipschitz map from $(1+7\a)l_{ij}^{\nu'}$ to $(1+7\a)B_{Q_{ij}}(0, R_1)$. 

We claim that
\be F_2\cap (1+7\a)\mL_{ij}^{\nu'}\bs (1+5\a)\bar\cU \mbox{ coincide with the cone over the graph of }u_{ij}\mbox{ in }(1+7\a)\mL^{\nu'}_{ij}\bs (1+5\a)\cU.\ee

Since $(1+7\a)\mL_{ij}^{\nu'}\bs (1+5\a)\bar\cU\subset V$, it is enough to look at $F_2\cap V$. Recall that $F_2=\pi\circ\psi''\circ\psi'\circ\psi_1\circ\d_{1+9\a}(E_1\cap\bar\cU)=\pi\circ\psi''\circ\psi'(F_1)$.

We know that $|\cS'|\cap V=\emptyset$, $\phi''\circ\phi'|_{|\cS'|^C}=Id$, and $\phi''\circ\phi'(|\cS'|)\subset |\cS'|$, hence $\psi''\circ\psi'\circ(F_1\cap V)=F_1\cap V=G_1\cap V$. By definition of $V$, $\pi(V)\subset V$ and $\pi(V^C)\subset V^C$, thus 
\be F_2\cap V=\pi\circ\psi''\circ\psi'(F_1)\cap V=\pi (\psi''\circ\psi'(F_1)\cap V)=\pi (F_1\cap V)=\pi(G_1\cap V).\ee

By \tb{(4.14)}, $G_1\cap V$ coincide with a cone in $V$. Hence by the definition of $V$, we know that $\pi(G_1\cap V)=G_1\cap V\cap(1+7\a)\bar\cU$. In particular, 
\be \begin{split}&F_2\cap(1+7\a)\mL_{ij}^{\nu'}\bs (1+5\a)\bar\cU=F_2\cap V\cap (1+7\a)\mL_{ij}^{\nu'}\bs (1+5\a)\bar\cU\\
=&G_1\cap V\cap(1+7\a)\bar\cU\cap (1+7\a)\mL_{ij}^{\nu'}\bs (1+5\a)\bar\cU
=G_1\cap (1+7\a)\mL_{ij}^{\nu'}\bs (1+5\a)\bar\cU.
\end{split}\ee
Again by \tb{(4.14)}, we get Claim \tb{(4.68)}. As a result, \tb{(4.67)} holds, and hence $F'\cap \pa\cU$ is a $(\eta,\d,\nu', L)$-sliding boundary of $Z$. Since $F'$ is a $\d$-sliding deformation of $E_1$ in $\bar\cU$, and $E_1$ is a $\Z_2$-topological competitor for $E$ in $\cU$, hence $F'$ is an $(\eta,\d,\nu', L)$-$\Z_2$-topological sliding competitor for $Z$. We have alreay proved that \tb{(4.5)} holds. Thus we complete the proof of Proposition \tb{4.2}.\qed

\noindent\textbf{Proof of Proposition \tb{4.3}.}  By definition, since $F$ is a 2-regular $(\eta,\d,\nu, L)$-$\Z_2$-topological sliding competitors for $Z$, there exists a $\Z_2$-topological competitor $E_0$ for $Z$ in $\cU$, and a $\d$-sliding deformation $\psi_t,1\le t\le 1$ in $\bar\cU$, so that $\psi_1(E_0\cap \bar\cU)=F$. The problem is just that $E_0$ has no reason to be 2-regular.

Since $F$ is 2-regular, there exists $\e\in (0,\d)$, so that there exists a $\e$-neighborhood retract to $F$ in $\bar\cU$. That is, there exists a Lipschitz deformation retract $g_t, 0\le t\le 1$ from $B(F,\e)\cap \bar\cU\to B(F,\e)\cap \bar\cU$, so that $g_t(x)=x,\forall x\in F, \forall t\in [0,1]$, $g_0=id$, $g_1(B(F,\e)\cap \bar\cU)=F$, and $g_t(\pa\cU)\subset \pa\cU$. We extend $g_t$ to a $\d$-sliding deformation in $\bar\cU$.

Let $L$ be the Lipschitz constant for $\psi_1$. Fix an $n\in \N$ such that $2^{-n-4}L<\e$. Let $\Delta$ denote the set of all dyadic cubes of size $2^{-n}$ in $\R^4$, and let $\Delta^d,0\le d\le 4$ be the $d$-skeleton of $\Delta$. Let $f: E_0\cap\bar\cU\to \Delta^2$ be the Federer Fleming projection. Then it is easy to see that 
\be |f(x)-x|\le 2^{-n-1}, \forall x\in E_0\cap\bar\cU,\ee
hence $f(E_0\cap\bar\cU)\subset (1+2^{-n-2})\bar\cU$.

Let $\a=2^{-n}$, and set
\be h: E_0\to [(1+\a)\bar\cU]\cup Z: h(x)=\left\{\begin{array}{rcl} f(x),&\ & x\in E_0\cap \bar\cU;\\
\frac{t-1}{\a}x+\frac{\a-t}{\a}f(x),&\ &x\in E_0\cap t\pa\cU, 1\le t\le \a;\\
x,&\ & x\in E_0\bs (1+\a)\cU.\end{array}\right.\ee

We extend $h$ to a Lipschitz deformation $h_t: (1+2\a)\cU\to (1+2\a)\cU$, so that $h_0=id$, $h_1=h$ on $E_0$.
Then by Proposition \tb{2.7}, $h(E_0)$ is a $\Z_2$-topological competitor for $E_0$ in $(1+2\a)\cU$, and hence is a $\Z_2$-topological competitor for $Z$ in $(1+2\a)\cU$. Moreover, $h(E_0)$ is 2-regular, and by \tb{(4.71) and (4.72)}, we know that 
\be d_H(h(E_0), E_0)\le 2^{-n-1}.\ee

Set $E=\frac{1}{1+2\a}h(E_0)$. Then $E$ is a 2-regular $\Z_2$-topological competitor for $Z$ in $\cU$. and $d_H(E, h(E_0))\le 2\a$. Combine with \tb{(4.73)} we know that 
\be d_H(E, E_0)\le 2\a+2^{-n-2}=2^{-n-3}.\ee

As a result, we know that 
\be \psi_1(E\cap\bar\cU)\subset B(\psi_1(E_0\cap\bar\cU),2^{-n-3}L)=B(F, 2^{-n-3}L)\subset B(F,\frac \e2).\ee

Thus we have
\be g_1( \psi_1(E\cap\bar\cU))\subset F.\ee

All together, we get the 2-regular $\Z_2$-topological competitor $E$ for $Z$ in $\cU$, and the $(\eta,\d)$-sliding deformation $\varphi_t:=g_t\circ \psi_t$, so that $\varphi_1(E\cap \bar\cU)\subset F$.

The rest is to prove that $\varphi_1(E\cap \bar\cU)$ is a $(\eta,\d,\nu, L)$-$\Z_2$-topological sliding competitor for $Z$. That is, if we set $G:=\varphi_1(E\cap \bar\cU)\cap \pa \cU$, then $G$ satisfies \tb{(3.26)} in Definition \tb{3.10}. 

By definition, $\varphi_t$ sends $\pa \cU$ to $\pa \cU$, and is a $\d$-sliding deformation in $\bar\cU$, hence $\varphi_t(\pa Z)\subset B(\pa Z, \d)\cap \pa \cU$, and $\varphi_1(\pa Z)=\varphi_1(\pa E)\subset G\subset \pa F:=F\cap \pa\cU$. 

Now for each $1\le j\le l\le 3$. Let $x_{jl}$ denote the midpoint of $l_{jl}$. Let $\sigma=\{x_{jl}\}\times \pa B_{Q_{jl}}(0, R_1)$. Then $\sigma\subset \pa\cU$, and it represents a non zero element in $H_1(\pa\cU\bs\pa Z,\Z_2)$. Moreover, since $G\subset \pa F\subset B(\pa Z, \d)$, hence $G\cap \sigma=\emptyset$.

Now since $\varphi_t(\pa Z)\subset B(\pa Z, \d)\cap \pa \cU$, and $\d<R_1$, we know that $\{\varphi_t(\pa Z), 1\le t\le 1\}$ does not meet $\sigma$. 

We claim that
\be [\sigma]\ne 0\mbox{ in } H_1(\pa\cU\bs \varphi_1(\pa Z),\Z_2).\ee

Suppose not, then there is a simplicial $\Z_2$-2 chain $\Gamma\subset \pa\cU\bs \varphi_1(\pa Z)$ so that $\pa\Gamma=\sigma$. Let $\e>0$ be such that $dist(\Gamma, \varphi_1(\pa Z))>3\e$.

Let $\varphi:\pa\cU\to \pa\cU$ be such that $\varphi=\varphi_1$ on $\pa Z$, $\varphi=id$ on $\pa\cU\bs  B(\pa Z, \e)$.

As a result, we can find a smooth map $k:\pa\cU\to \pa \cU$, transverse to $\Gamma$, $k=id$ on $\pa\cU\bs  B(\pa Z, 2\e)$, and 
\be ||k-\varphi||_\infty\le \e.\ee

Then $k^{-1}(\Gamma)$ is a simplicial $\Z_2$-2 chain in $\bar\cU$, $\pa(k^{-1}(\Gamma))=\sigma$, and for any $x\in k^{-1}(\Gamma),y\in \pa Z$, we know that $k(x)\in \Gamma, k(y)\in k(\pa Z)$, hence
\be |k(x)-k(y)|\ge |k(x)-\varphi(y)|-|\varphi(y)-k(y)|\ge \mbox{dist}(\Gamma, \varphi(\pa Z))-\e\ge 3\e-\e=2\e>\e.\ee
Hence $x\ne y$. As a result, we know that $k^{-1}(\Gamma)\cap \pa Z=\emptyset$. 

Thus we have found a a simplicial $\Z_2$-2 chain in $\bar\cU\bs \pa \Z$, whose boundary is $\sigma$. This contradicts the fact that $[\sigma]\ne 0\mbox{ in } H_1(\pa\cU\bs \pa \Z, \Z_2).$ Thus we have proved Claim \tb{(4.77)}.

Now since $\varphi(\pa Z)\subset G$, hence $\bar\cU\bs G\subset \bar\cU\bs \varphi(\pa Z)$. As a result, we know that 
\be [\sigma]\ne 0\mbox{ in } H_1(\pa\cU\bs G,\Z_2).\ee

Finally we claim that, for $1\le i,j\le 3$.
\be G\cap L_{ij}^\nu=\pa F\cap L_{ij}^\nu.\ee
In fact, we know that $\pa F\cap L_{ij}^\nu$ is a graph of a lipschitz map $\xi$ from $l_{ij}^\nu$ to $B_{Q_{ij}}(0,R_1)$, and that $G\subset \pa F$. So if $G\cap L_{ij}^\nu\ne \pa F\cap L_{ij}^\nu$, there exists $p\in (\pa F\cap L_{ij}^\nu)\bs(G\cap L_{ij}^\nu)$. Since $\pa F\cap L_{ij}^\nu$ is the graph of $\xi$, there exists $z\in l_{ij}^\nu$ so that $p=(z,\xi(z))$, and that $\pa F\cap \{z\}\times B_{Q_{ij}}(0,R_1)=\{z,\xi(z)\}=\{p\}$. As a result, we know that $G\cap \{z\}\times B_{Q_{ij}}(0,R_1)=\emptyset$, and hence $ \{z\}\times B_{Q_{ij}}(0,R_1)\subset L_{ij}^\nu\bs G$. Note that $\pa [\{z\}\times B_{Q_{ij}}(0,R_1)]=\{z\}\times \pa B_{Q_{ij}}(0,R_1)\sim \sigma$ in $l_{ij}^\nu\times \pa B_{Q_{ij}}(0,R_1)$, hence $[\s]=0$ in $H_1(l_{ij}^\nu\times \pa B_{Q_{ij}}(0,R_1), \Z_2)$. But since $F$ is an $(\eta,\d)-\Z_2-$topological sliding competitor for $Z$, we know that $\pa F\cap L_{ij}^\nu \subset l_{ij}^\nu\times \bar B_{Q_{ij}}(0, \d)\subset l_{ij}^\nu\times B_{Q_{ij}}(0, R_1)$, hence $l_{ij}^\nu\times \pa B_{Q_{ij}}(0,R_1)\subset \pa\cU\bs \pa F\subset \pa\cU\bs G$, thus $[\s]=0$ in $H_1(\pa\cU\bs G, \Z_2)$. This contradicts \tb{(4.80)}.

Thus we get \tb{(4.81)}, and the proof is completed.\qed

\section{Decomposition of a competitor, and estimate its measure by projections}
 
 After the simplifications in the last section, we begin to carry on the proof of sliding stabilities for $Z$ in this and the next section.

Recall that in \cite{YXY}, we proved the minimality of $Z$ in $B=B(0,1)$ using product of paired calibrations. There are two main steps : 

\textbf{Step 1: Decomposition.} since we are treating objects in codimension 2, we lose the separation conditions needed for normal paired calibrations. This leads us to work in $\Z_2$ homology groups, so that we can properly decompose any regular competitor $F$ for $Z$ into 9 parts $F_{ij}, 1\le i,j\le 3$, and these 9 parts can only meet each other in some prescribed way. 

\textbf{Step 2: measures of projections.} Then we use 9 different projections $\pi_{ij}$ to project these 9 parts to the boundary $\pa \cU$, and use the sum of the measures of $\pi_{ij}(F_{ij}), 1\le i,j\le 3$ of these nine parts to give a lower bound for the measure of competitors for $Z$. 

Due to Step 1, we only managed to prove the topological minimality in the coefficient group $\Z_2$. This will be the same for the topological sliding stability and uniqueness property: the coefficient group has to be $\Z_2$. On the other hand, in the proof of the minimality of  $Z$, since for each Almgren or $\Z_2$ topological competitor $F$, we have $F\cap \pa \cU=Z\cap \pa \cU$, hence the measure of the projections $\pi_{ij}(F_{ij}), 1\le i,j\le 3$ are more or less fixed. But in the sliding case, we allow $F\cap \pa \cU$ to be different from $Z\cap \pa \cU$, hence we have one more thing to do : to give a uniform lower bound for the sum of the projections. 

We have done similar argument in the proof of the sliding stability for $\T$ and $\Y$ sets in dimension 3. In their case, the unit sphere is of dimension 2, and hence no matter how we do the sliding deformation, the above projections are just disjoint regions in the unit sphere, whose union is the unit sphere. Therefore the sum of the measure of the projections are always the same. This allows us to prove the $(\eta,\d)$-topological sliding stability
for them.

But for the case of $Z$, the boundary $\pa B$ is of dimension 3, and our projections of sets are of dimension 2, which is not of full dimension, and this makes the estimate of the sum of the measures after projections much more involved. And at last, we only arrive to prove the $(\eta,\d,\nu,L)$-$\Z_2$-topological sliding stability.

As a final remark, we do not know whether the set $Z$ is also $(\eta,\d)-Z_2$-topological sliding stable. But as we mentioned in the introduction, this $(\eta,\d,\nu,L)$-$\Z_2$-topological sliding stability is enough for our purpose of finding new singularities by taking almost orthogonal unions with other minimal cones.

\subsection{The decomposition}

Fix any $\eta<\eta_1$, $\d\in(0,R_1)$, $\nu\in (0,R_3)$, and $L>0$. Take a set $F\in \F=\F(\eta,\d,\nu, L)$. 


Set $\pa F=F\cap \pa\cU$. Then it is a deformation of $\pa Z:=Z\cap\pa\cU$.

For $1\le i,j\le 3$, set $L_i:=Y_1\bs R_{o_1a_i}$ ($R_{ab}$ denotes the ray issued from $a$ and passing through $b$, for $a,b\in \R^4$), and $S_j:=Y_2\bs R_{o_2b_j}$. Then $L_i$ and $S_j$ are cones in $P_1$ and $P_2$ respectively. For $1\le i,j\le 3$, set 
\be Z_{ij}=(L_i\times S_j)\cap \bar\cU,\ee
Then it is a $\Z_2$-chain, whose boundary is the essentially disjoint union
\be z_{ij}:=\pa Z_{ij}=\cup_{k\ne i, l\ne j}l_{kl}.\ee
Let $z_{ij}$ and $l_{ij}$ also denote the $\Z_2$-simplicial chain supported on $z_{ij}$ and $l_{ij}$ in the obvious sense.

By \tb{(5.2)}, we have
\be \sum_{k\ne i,l\ne j}l_{kl}=z_{ij}\sim 0\mbox{ in }H_1(Z\cap\bar\cU,\Z_2), \forall 1\le i,j\le 3.\ee

\begin{pro}If $E$ is a 2-regular $\Z_2$-topological competitor for $Z$ in $\cU$, then $[z_{ij}]=0$ in $H_1(E\cap\bar\cU, \Z_2)$, $\forall 1\le i,j\le 3$.
\end{pro}

\nd Let $E$ be a 2-regular $\Z_2$-topological competitor for $Z$ in $\cU$.
 
Recall that the set $Z=Y\times Y$ is the union of nine $\frac14-$planes $H_{ij}:=R_{o_1,a_i}\times R_{o_2,b_j}, 1\le i,j,\le 3$. Set $e_{ij}=(a_i,b_j)\in\R^4$. Set $s_{ij}=e_{ij}+\{x\in Q_{ij}: |x|=\frac {1}{10}\}$. Then $s_{ij}$ is a small circle outside $\cU$ that links $H_{ij}$. We denote also by $s_{ij}$ the corresponding element in homology groups with coefficients in $\Z_2$ for short. Then in $H_1(\R^4\bs E,\Z_2)$
 \be s_{i1}+s_{i_2}+s_{i3}=0=s_{1j}+s_{2j}+s_{3j}, \mbox{ for }1\le i,j\le 3,\ee
 and, in the special case of $Z=Y\times Y$, $H_1(\R^4\bs Z,\Z_2)$ is the Abelian group generated by $s_{ij},1\le i,j\le 3$ with the relations \tb{(5.4)}. Notice that the relations \tb{(5.4)} has in fact only 5 independent relations. Thus $H_1(\R^4\bs Z,\Z_2)$ is in fact a vector space (since $\Z_2$ is a field) with basis $\{s_{ij},1\le i,j\le 2\}$.
 
 Take $i=1,j=1$ for example. We want to show that $[z_{11}]=0$ in $H_1(E\cap \bar\cU,\Z_2)$.
 
 Set $F=[E\cap\cU]\cup [(L_3\times S_3)\bs \cU]$. Note that the topological plane $L_3\times S_3$ is the essentially disjoint union $\cup_{i,j\ne 3}H_{ij}$; hence the four circles $s_{11},s_{12},s_{21}$ and $s_{22}$ represent the same element in $H_1(\R^4\bs F,\Z_2)$. Denote by $s$ this element in $H_1(\R^4\bs F,\Z_2)$.
 
 We want to show first that $s\ne 0$ in $H_1(\R^4\bs F,\Z_2)$. Suppose not, that is, $s=0$ in $H_1(\R^4\bs F,\Z_2)$. Then there exists a $C^1$ simplicial 2-chain $\Gamma$ in $\R^4\bs F$ such that $\partial \Gamma=s_{11}$. Since $E$ is a topological competitor of $Z$ in $\cU$, $s_{11}\ne 0$ in $H_1(\R^4\bs E,\Z_2)$, hence $\Gamma\cap E\ne\emptyset$. But $\Gamma\subset\R^4\bs F$, hence $\Gamma$ can only meet $E$ at $E\bs F=[H_{13}\cup H_{23}\cup H_{31}\cup H_{32}\cup H_{33}]\bs \cU$. We can also ask that $\Gamma$ meet these five $\frac14$ planes transversally, and do not meet any of their intersections. This gives that in $H_1(\R^4\bs E,\Z_2)$,
 \be s_{11}+\d_{13}s_{13}+\d_{23}s_{23}+\d_{33}s_{33}+\d_{31}s_{31}+\d_{32}s_{32}=0\ee
 with $\d_{ij}\in \Z_2$, and at least one of the five $\d_{ij}$ is not zero.
 
Combine with \tb{(5.4)}, we have
 \be s_{11}+\d_{13}[s_{11}+s_{12}]+\d_{23}[s_{21}+s_{22}]+(\d_{31}+\d_{33})[s_{21}+s_{11}]+(\d_{32}+\d_{33})[s_{12}+s_{22}]=0.\ee
 We simplify and get
 \be [1+\d_{13}+\d_{31}+\d_{33}]s_{11}+[\d_{13}+\d_{32}+\d_{33}]s_{12}+[\d_{23}+\d_{31}+\d_{33}]s_{21}+[\d_{23}+\d_{32}+\d_{33}]s_{22}=0,\ee
in $H_1(\R^4\bs E,\Z_2)$.

But $s_{11},s_{12},s_{21},s_{22}$ are independent elements in $H_1(\R^4\bs Z,\Z_2)$, hence they are also independent in $H_1(\R^4\bs E,\Z_2)$, because $E$ is a $\Z_2$-topological competitor for $Z$. Thus \tb{(5.7)} gives
\be 1+\d_{13}+\d_{31}+\d_{33}=\d_{13}+\d_{32}+\d_{33}=\d_{23}+\d_{31}+\d_{33}=\d_{23}+\d_{32}+\d_{33}=0.\ee
However, the sum of the four numbers gives
\be \begin{split}0&=[1+\d_{13}+\d_{31}+\d_{33}]+[\d_{13}+\d_{32}+\d_{33}]+[\d_{23}+\d_{31}+\d_{33}]+[\d_{23}+\d_{32}+\d_{33}]\\
&=1+2(\d_{13}+\d_{31}+\d_{32}+\d_{23}+2\d_{33})=1,
 \end{split}\ee
 which is impossible. Hence we have a contradiction.
 
 Hence $s_{11}=s\ne 0$ in $H_1(\R^4\bs F,\Z_2)$.

 Now because $F$ is 2-regular in $\cU$ and $F\bs\cU=(L_1\times S_1)\bs \cU$, $H_1(\R^4\bs F,\Z_2)$ is a finite dimensional vector space (because $\Z_2$ is a field), and hence $H^1(\R^4\bs F,\Z_2)=[H_1(\R^4\bs F,\Z_2)]^*$. For the non zero element $s_{11}\in H_1(\R^4\bs F,\Z_2)$, denote by $\zeta$ the dual element of $s_{11}$ in the cohomology group $H^1(\R^4\bs F,\Z_2)$. Denote by $\varphi$ the isomorphism of the Poincar\'e-Lefschetz duality
 
 \be\a:H^1(\R^4\bs F,\Z_2)\cong H_3(\R^4, F ,\Z_2),\ee
 Then $\a(\zeta)$ can be represented by a smooth simplicial 3-chain $\Sigma$ with boundary in $F$, and $|\Sigma|\cap s_{11}$ is a single point. Denote by $\xi=\partial\Sigma$, then this is a 2-chain with support in $F$ such that $s_{11}$ is non-zero in $H_1(\R^4\bs |\xi|,\Z_2)$. 
 
 Notice that outside $\cU$, the set $F$ is topologically a plane, which is linked by $s_{11}$, hence if $s_{11}$ is non-zero in $H_1(\R^4\bs |\xi|,\Z_2)$, then $|\xi|\supset (F\bs \cU)$. Since $|\xi|\subset F$, we know that 
 \be |\xi|\bs \cU=F\bs \cU.\ee
  
 But $F\cap\partial \cU=z_{11}$, hence $z_{11}=\partial(F\bs \cU)=\pa\pa\Sigma+\partial(F\bs \cU)=\partial (\xi+(F\bs \cU))$ (here we regard $F\bs \cU$ as a 2-chain). By \tb{(5.11)},  the support of $\xi+(F\bs \cU)$ is contained in $F\cap \cU$, hence $z_{11}$ is a boundary in $F\cap \bar\cU$, which yields that $z_{11}$ represents a zero element in $H_1(F\cap \bar\cU,\Z_2)$, and thus in $H_1(E\cap \bar\cU,\Z_2)$.
 
 The same arguments holds for all $z_{ij},1\le i,j\le 3$. Thus the proof of Proposition \tb{5.1} is completed.\qed
 
After Proposition \tb{5.1}, we know that if $F\in \F(\eta,\d,\nu, L)$, then 
there exists a 2-regular $\Z_2$-topological competitor $E$ for $Z$ in $\cU$, and a $\d$-sliding deformation $\varphi_1$ of $E$ in $\cU$, such that $F=\varphi_1(E\cap \bar\cU)$. As a result, if we denote by $\varphi=\varphi_1$, then
\be \s_{ij}:=\varphi_*(z_{ij})\sim 0\mbox{ in }H_1(F,\Z_2),\ee  
and 
\be \sum_{k\ne i,l\ne j}\varphi_*(l_{kl})=\s_{ij}\sim 0\mbox{ in }H_1(F,\Z_2), \forall 1\le i,j\le 3,\ee
by \tb{(5.3)}.

Recall that $Q_{ij}$ denote the 2-subspace of $\R^4$ generated by $x_i$ and $y_j$, and let $q_{ij}$ be the orthogonal projection from $\R^4$ to $Q_{ij}$, $1\le i,j\le 3$. Let $P_{ij}$ denote the 2-subspace containing $\mc_{ij}$ (or equivalently, $\g_{ij}$). Then $P_{ij}$ is just the 2-subspace of $\R^4$ generated by $a_i$ and $b_j$, and is the 2-subspace of $\R^4$ orthogonal to $Q_{ij}$.

\begin{pro}Let $F\in \F(\eta,\d,\nu, L)$ for some $\eta<\eta_1$, $\d\in (0,R_1)$, $\nu\in (0,R_3)$ and $L>0$. Let $\s_{ij}$ be defined as above. Then there exists subsets $F_{ij}, 1\le i,j\le 3$ of $F$, so that the following holds:

$1^\circ$ $\H^2(q_{ij}(F_{ij}))\ge|q_{ij*}(\s_{ij})|_2$;

$2^\circ$ For each $1\le i\le 3$, $\H^2-$almost every point in the union $F_{i1}\cup F_{i2}\cup F_{i3}$ belongs to exactly two of $F_{i1}, F_{i2}$ and $F_{i3}$; similarly, for each $1\le j\le 3$, $\H^2-$almost every point in the union $F_{1j}\cup F_{2j}\cup F_{3j}$ belongs to exactly two of $F_{1j}, F_{2j}$ and $F_{3j}$.
\end{pro}

\nd
For $i,j=1,2$, by \tb{(5.12)}, let $E_{ij}$ be a $\Z_2$ chain in $F$ so that $\pa E_{ij}=\s_{ij}$. Set $E_{3j}=E_{1j}+E_{2j}$, for $j=1,2$, and then set $E_{i3}=E_{i1}+E_{i2}$ for $1\le i\le 3$.
Then  we know that 
\be \pa E_{3j}=\pa E_{1j}+\pa E_{2j}=\s_{1j}+\s_{2j}=\s_{3j},j=1,2,\ee
and hence
\be \pa E_{i3}=\pa E_{i1}+\pa E_{i2}=\s_{i3}.\ee

Moreover we have
\be E_{33}=E_{31}+E_{32}=(E_{11}+E_{21})+(E_{12}+E_{22})=(E_{11}+E_{12})+(E_{21}+E_{22})=E_{13}+E_{23}.\ee

Altogether we have
\be \sum_{1\le i\le 3}E_{ij}=0, \forall 1\le j\le 3\mbox{ and }\sum_{1\le j\le 3}E_{ij}=0, \forall 1\le 1\le 3,\ee
and
\be\pa E_{ij}=\s_{ij}, \forall 1\le i,j\le 3.\ee

It is important to point out that, by \tb{(5.17)}, for any $1\le i\le 3$, modulo a $\H^2$ negligible set, the supports of the three 2-chains satisfy that :
\be |E_{i3}|=(|E_{i2}|\cup|E_{i1}|)\bs(|E_{i2}|\cap|E_{i1}|).\ee
In other words, $\H^2-$almost every point in the union $|E_{i1}|\cup|E_{i2}|\cup|E_{i3}|$ of the three supports  belongs to exactly two of them.

Similarly, we know that for each $1\le j\le 3$, $\H^2-$almost every point in the union $|E_{1j}|\cup|E_{2j}|\cup|E_{3j}|$ of the three supports  belongs to exactly two of them.

Now let us turn to the projections of the nine $\Z_2$-chains $E_{ij}$.

For every $1\le i,j\le 3$, we know that $\s_{ij}$ is a closed $\Z_2$ chain, and hence $q_{ij*}(\s_{ij})$ is also a closed $\Z_2$ chain in the 2-plane $Q_{ij}$. Since $H_3(Q_{ij},\Z_2)=0$, hence there exists a unique $\Z_2$-2-chain $\S_{ij}$ in $Q_{ij}$, so that $\pa\S_{ij}=q_{ij*}(\s_{ij})$.

As a result, for any $\Z_2$-simplicial 2 chain $\G\subset \R^4$, if $\pa\G=\s_{ij}$, then $\pa (q_{ij*}(\G))=q_{ij*}(\pa\G)=q_{ij*}(\s_{ij})$, and $q_{ij*}(\G)\subset Q_{ij}$. By the uniqueness of $\S_{ij}$, we know that $q_{ij*}(\G)=\S_{ij}$. That is
\be \forall \G\subset \R^4, q_{ij*}(\pa \G)=q_{ij*}(\s_{ij})\Rightarrow q_{ij*}(\G)=\S_{ij}.\ee

In particular, we know that $q_{ij*}(E_{ij})=\S_{ij}$, and hence
\be \H^2(q_{ij}(|E_{ij}|))\ge \H^2(|q_{ij*}(E_{ij})|)=\H^2(|\S_{ij}|)=|q_{ij*}(\s_{ij})|_2.\ee

It is now enough to set $F_{ij}=|E_{ij}|$. \qed

\subsection{The calibrations and projections}

\begin{pro}Let $F\in \F(\eta,\d,\nu, L)$ for some $\eta<\eta_1$, $\d\in (0,R_1)$, $\nu\in (0,R_3)$ and $L>0$. Let $\s_{ij}$ be defined as before, then
\be \H^2(F)\ge \frac 13 \sum_{1\le i\le 3}\sum_{1\le j\le 3}|q_{ij*}(\s_{ij})|_2.\ee 
\end{pro}

\nd For $1\le i,j\le 3$, set $v_{ij}=a_i\wedge b_j$, a unit simple 2-vector in $\R^4$.

We define the function $f_{ij}$ on the set of simple 2-vectors in $\R^4$ : for any simple 2-vector $\xi\in \wedge_2(\R^4)$, $f_{ij}(\xi):=|\xi\wedge v_{ij}|=|\det_{e_1,e_2,e_3,e_4}\xi\wedge v_{ij}|$, with $\{e_j\}_{1\le j\le 4}$ the canonical orthonormal basis of $\R^4$. Now for any unit (with respect to the $L^2$ norm $|\cdot|$ for the orthonormal basis $\{e_i\wedge e_j\}_{1\le i<j\le 4}$ of $\wedge_2(\R^4)$) simple 2-vector $\xi$, we can associate to it a plane $P(\xi)\in G(4,2)$, where $G(4,2)$ is the set of all 2-dimensional subspaces of $\R^4$ :
\be P(\xi)=\{v\in \R^4, v\wedge\xi=0\}.\ee
In other words, $P(x\wedge y)$ is the subspace generated by $x$ and $y$.

Now denote also by $g_{ij}$ the function from $G(4,2)$ to $\R$: for any $P=P(x\wedge y)\in G(4,2)$ with $x\wedge y$ a unit simple 2-vector, $g_{ij}(P)=f_{ij}(x\wedge y)$. Since the definition of $f_{ij}$ on $\wedge_2(\R^4)$ is to take the absolute value of the determinant, the function $g_{ij}$ is well defined.

Let $F_{ij}, 1\le i,j\le 3$ be as in Proposition \tb{5.2}. Now since each $F_{ij}$ is the support of a smooth simplicial 2-chain in $\R^4$, it is 2-rectifiable, and the tangent plane $T_xF_{ij}$ of $F_{ij}$ at  $x$ exists for $\H^2-$almost all $x\in F_{ij}$. We want to estimate 
\be \int_{F_{ij}}g_{ij}(T_xF_{ij}) d\H^2(x).\ee

Notice that $F_{ij}$ is piecewise smooth, hence $g_{ij}(T_xF_{ij})$ is measurable. Note also that $|g_{ij}|\le 1$, hence the integral is well defined.

Denote by $E_{ij}$ the subset $\{x\in F_{ij} : J_2q_{ij}(x)\ne 0\}$ of $F_{ij}$, where $J_2q_{ij}(x)$ is the Jacobian of the restriction $q_{ij}|_{F_{ij}}:F_{ij}\to Q_{ij}$. Then $T_xE_{ij}=T_xF_{ij}$ for $\H^2$ almost all $x\in E_{ij}$. By the Sard theorem, we have
\be \H^2(q_{ij}(F_{ij})\bs q_{ij}(E_{ij}))=0,\ee
and hence by Proposition \tb{5.2}, 
\be \H^2(q_{ij}(E_{ij}))=\H^2(q_{ij}(F_{ij}))\ge |q_{ij*}(\s_{ij})|_2.\ee
 
 Now for $x\in E_{ij}$, define $h_{ij}(x)=g_{ij}(T_xF_{ij}) (J_2q_{ij}(x))^{-1}$. 
Recall that the projection $q_{ij}: F_{ij}\to Q_{ij}$ is a 1-Lipschitz function.  hence by the coarea formula for Lipschitz functions between two rectifiable set (cf.\cite{Fe} Theorem 3.2.22), we have
\be \int_{E_{ij}}h_{ij}(x) J_2q_{ij}(x) d\H^2(x)=\int_{q_{ij}(E_{ij})} d\H^2(y)[\sum_{q_{ij}(x)=y}h_{ij}(x)].\ee

By definition of $h_{ij}$, the left hand side of the above equality is just
\be \begin{split}\int_{E_{ij}}h_{ij}(x) J_2q_{ij}(x) d\H^2(x)&=\int_{E_{ij}}g_{ij}(T_xF_{ij}) (J_2q_{ij}(x))^{-1} J_2q_{ij}(x) d\H^2(x)\\
&=\int_{E_{ij}}g_{ij}(T_xF_{ij})d\H^2(x)\le \int_{F_{ij}}g_{ij}(T_xF_{ij})d\H^2(x),\end{split}\ee
the last inequality is because $g_{ij}$ is non negative.

For the right hand side of \tb{(5.28)}, note that for almost all $x\in E_{ij}$, the tangent plane $T_xF_{ij}$ exists. Suppose $T_xF_{ij}=P(u\wedge v)$, with $u,v$ an orthonormal basis of $T_xF_{ij}$. Hence
\be g_{ij}(T_xF_{ij})=f_{ij}(u\wedge v)=|v_{ij}\wedge u\wedge v|=|v_{ij}\wedge q_{ij}(u\wedge v)|.\ee

Notice that $q_{ij}(u\wedge v)\in\wedge_2(Q_{ij})$, hence if we take a unit simple two vector $\xi_{ij}$ of $Q_{ij}$, we have $q_{ij}(u\wedge v)=\pm|q_{ij}(u\wedge v)|\xi_{ij},$ and hence by \tb{(5.29)}
\be g_{ij}(T_xF_{ij})=|q_{ij}(u\wedge v)||v_{ij}\wedge\xi_{ij}|=|q_{ij}(u\wedge v)|=J_2q_{ij}(x),\ee
and thus 
\be h_{ij}(x)=g_{ij}(T_xF_{ij})J_2q_{ij}(x)^{-1}=1\ee
 for $\H^2-a.e.$ $x\in E_{ij}$.
 As a result, 
 \be \int_{q_{ij}(E_{ij})} d\H^2(y)[\sum_{q_{ij}(x)=y}h_{ij}(x)]=\int_{q_{ij}(E_{ij})} d\H^2(y)[\sum_{q_{ij}(x)=y}1]
 =\int_{q_{ij}(E_{ij})} d\H^2(y)\sharp \{q_{ij}^{-1}(x)\cap E_{ij}\}.\ee
But for all $y\in q_{ij}(E_{ij})$, $\sharp \{q_{ij}^{-1}(x)\cap E_{ij}\}\ge 1$, hence
\be \int_{q_{ij}(E_{ij})} d\H^2(y)[\sum_{q_{ij}(x)=y}h_{ij}(x)]\ge \int_{q_{ij}(E_{ij})} d\H^2(y)=\H^2(q_{ij}(E_{ij}))\ge |q_{ij*}(\s_{ij})|_2.\ee

Combine \tb{(5.27), (5.28) and (5.33)} we have
\be |q_{ij*}(\s_{ij})|_2\le \int_{F_{ij}}g_{ij}(T_xF_{ij})d\H^2(x), \mbox{ for }1\le i,j\le 3.\ee

We sum over $1\le i,j\le 3$, and have
\be \begin{split}\sum_{1\le i\le 3}\sum_{1\le j\le 3}|q_{ij*}(\s_{ij})|_2&\le \sum_{1\le i\le 3}\sum_{1\le j\le 3}\int_{F_{ij}}g_{ij}(T_xF_{ij})d\H^2(x)\\
&=\int_{\cup_{1\le i,j\le 3}F_{ij}}d\H^2(x) [\sum_{1\le j\le 3}\sum_{1\le l\le 3}g_{ij}(T_xF_{ij})1_{F_{ij}}(x)].\end{split}\ee

But $F$ is 2-rectifiable, and each $F_{ij}$ is its subset, hence we have for $\H^2-a.e. x\in F_{ij}$, $T_xF_{ij}=T_xF$. Hence we have
\be \sum_{1\le i\le 3}\sum_{1\le j\le 3}|q_{ij*}(\s_{ij})|_2\le \int_{\cup_{1\le i,j\le 3}F_{ij}}d\H^2(x) [\sum_{1\le i\le 3}\sum_{1\le j\le 3}g_{ij}(T_xF)1_{F_{ij}}(x)].\ee

Now we want to use Proposition \tb{5.2} to derive a essential upper bound for the function 
\be [\sum_{1\le i\le 3}\sum_{1\le j\le 3}g_{ij}(T_xF)1_{F_{ij}}(x)].\ee

Given a point $x\in \cup_{1\le j,l\le 3}F_{ij}$, by Proposition \tb{5.2} $2^\circ$, modulo a negligible set, there are two possibilities :

1) There exists $1\le i_1,i_2\le 3$ and $1\le j_1,j_2\le 3$ such that $x$ only belongs to the four pieces $F_{i_1j_1},F_{i_1j_2},F_{i_2j_1},F_{i_2j_2}.$

2) There exists a permutation $\sigma$ of $\{1,2,3\}$ such that $x$ belongs to all the nine $F_{ij}$ except for $F_{1\sigma(1)},F_{2\sigma(2)},F_{3\sigma(3)}$.

We will estimate the function $[\sum_{1\le i\le 3}\sum_{1\le j\le 3}g_{ij}(T_xF)1_{F_{ij}}(x)]$ in these two cases.

For 1), without loss of generality, we suppose that $i_1=j_1=1,i_2=j_2=2$.  Then 
\be \sum_{1\le i\le 3}\sum_{1\le j\le 3}g_{ij}(T_xF)1_{F_{ij}}(x)=\sum_{1\le i\le 2}\sum_{1\le j\le 2}g_{ij}(T_xF).\ee
 Suppose that $T_xF=P(\xi)$ with $\xi$ a unit simple 2-vector. Then for each $1\le i,j\le 2$, by definition of $g_{ij}$, $g_{ij}(T_xF)=|v_{ij}\wedge\xi|$.

Hence 
\be\begin{split}\sum_{1\le i\le 2}\sum_{1\le j\le 2}g_{ij}(T_xF)&=\sup_{\e}\sum_{1\le i,j\le 2}\e(i,j)\det(v_{ij}\wedge\xi)\\
&=\sup_{\e}\det[(\sum_{1\le i,j\le 2}\e(i,j)v_{ij})\wedge\xi]\le\sup_{\e}||\sum_{1\le i,j\le 2}\e(i,j)v_{ij}||.\end{split}\ee
where $\e$ run over all function from $\{1,2\}\times\{1,2\}\to\{1,-1\}$, and the norm $||\cdot ||$ on $\wedge_2(\R^4)$ is defined by 
\be ||\a||=\sup\{\det(\a\wedge \beta); \beta \mbox{ simple unit 2-vector}\}.\ee
Then the last inequality of \tb{(5.39)} is because $|\xi|=1$. Hence
\be \sum_{1\le i\le 3}\sum_{1\le j\le 3}g_{ij}(T_xF)1_{F_{ij}}(x)\le \sup\{||\sum_{1\le i,j\le 2}\e(i,j)v_{ij}||,\e:\{1,2\}\times\{1,2\}\to\{1,-1\}\}.\ee

Similarly, for 2), we have
\be \begin{split}\sum_{1\le i\le 3}&\sum_{1\le j\le 3}g_{ij}(T_xF)1_{F_{ij}}(x)=\sum_{1\le i,j\le 3, l\ne\sigma(i)}g_{ij}(T_xF)\\
&\le\sup\{||\sum_{1\le i,j\le 3, i\ne  j}\e(i,j)v_{ij}||,\e:[\{1,2,3\}\times\{1,2,3\}]\bs\{(1,1),(2,2),(3,3)\}\to\{1,-1\}\}.\end{split}\ee

The following lemma will lead to the conclusion of Proposition \tb{5.3}.

\begin{lem}[cf. \cite{YXY} Lemma 4.24]\be \sup\{||\sum_{1\le i,j\le 2}\e(i,j)v_{ij}||,\e:\{1,2\}\times\{1,2\}\to\{1,-1\}\}\le 3,\ee
and
\be \sup\{||\sum_{1\le i,j\le 3, i\ne j}\e(i,j)v_{ij}||,\e:[\{1,2,3\}\times\{1,2,3\}]\bs\{(1,1),(2,2),(3,3)\}\to\{1,-1\}\}\le 3.\ee
\end{lem}

By Lemma \tb{5.4}, \tb{(5.41) and (5.42)}, for $\H^2-$almost all $x\in \cup_{1\le i,j\le 3}F_{ij}$, 
\be\sum_{1\le i\le 3}\sum_{1\le j\le 3}g_{ij}(T_xF)1_{F_{ij}}(x)\le 3.\ee
Hence by \tb{(5.36)}, we have
\be \begin{split}\sum_{1\le i\le 3}\sum_{1\le j\le 3}|q_{ij*}(\s_{ij})|_2&\le \int_{\cup_{1\le i,j\le 3}F_{ij}}d\H^2(x) [\sum_{1\le i\le 3}\sum_{1\le j\le 3}g_{ij}(T_xF)1_{F_{ij}}(x)]\\
&\le 3\int_{\cup_{1\le i,j\le 3}F_{ij}}d\H^2(x)=3\H^2(\cup_{1\le i,j\le 3}F_{ij})\le 3\H^2(F),\end{split}\ee
which yields \tb{(5.22)}.\qed

%

\begin{cor}For the set $Z$, 
\be \H^2(Z\cap \bar\cU)=\frac 13 \sum_{1\le i\le 3}\sum_{1\le j\le 3}|q_{ij*}(z_{ij})|_2.\ee
\end{cor}

\nd For $F=Z$, we set $F_{ij}=Z_{ij}$, and then it is not hard to find that all inequalities in the proof of Proposition \tb{5.3} are equalities for $F=Z$.\qed

\section{Stability for $Y\times Y$}

In this section we are going to give a lower bound of the term $\sum_{1\le i\le 3}\sum_{1\le j\le 3}|q_{ij*}(z_{ij})|_2$ on the right-hand-side of \tb{(5.22)}. Note that the definition of this term depends only on $F\cap \pa \cU$. Hence in fact, it is defined for all $(\eta,\d,\nu,L)$-sliding boundary $G$ for $Z$. 

The aim of this section is to prove the following Theorem, and then deduce the $\Z_2$-topologically sliding stability for $Y\times Y$.

\begin{thm}Let $G$ be an $(\eta,\d,\nu, L)$-sliding boundary for $Z$ for some $\eta<\eta_1$, $\d<R_1$, $\nu\in (0,R_3)$, and $L>0$. Suppose that $\arctan\frac{R_3}{1-\sqrt\eta}<\frac{\pi}{16}$, $\d<\min\{\nu,\frac{R_3-\nu}{6}\}$, and $L<\frac{\d^2R_3}{(1-\sqrt\eta)(1-\eta)^2}$. Then we have
\be \sum_{1\le i\le 3}\sum_{1\le j\le 3}|q_{ij*}(\s_{ij})|_2\ge \sum_{1\le i\le 3}\sum_{1\le j\le 3}|q_{ij*}(z_{ij})|_2\ee
\end{thm}

The proof of Theorem 6.1 consists of a series of constructions and propositions. Let us first introduce some notation:
 
For a point $z\in P_{ij}$, let $\arg z\in \R/2\pi$ be the argument of $z$ under the basis $\{x_i,y_j\}$. Recall that 
\be\begin{split} l^\nu_{ij}=&\{(1-\sqrt\eta)a_i+tb_j:\nu\le t\le R_3\}\cup \{(1-\sqrt\eta)b_j+ta_i:\nu\le t\le R_3\}\\
&\cup \{(1-\eta)\cos\theta a_i+(1-\eta)\sin\theta b_j, \theta_0< \theta< \frac \pi 2-\theta_0\},
\end{split}\ee
where $\theta_0=\arccos\frac{1-\sqrt\eta}{1-\eta}$.

\centerline{\includegraphics[width=0.6\textwidth]{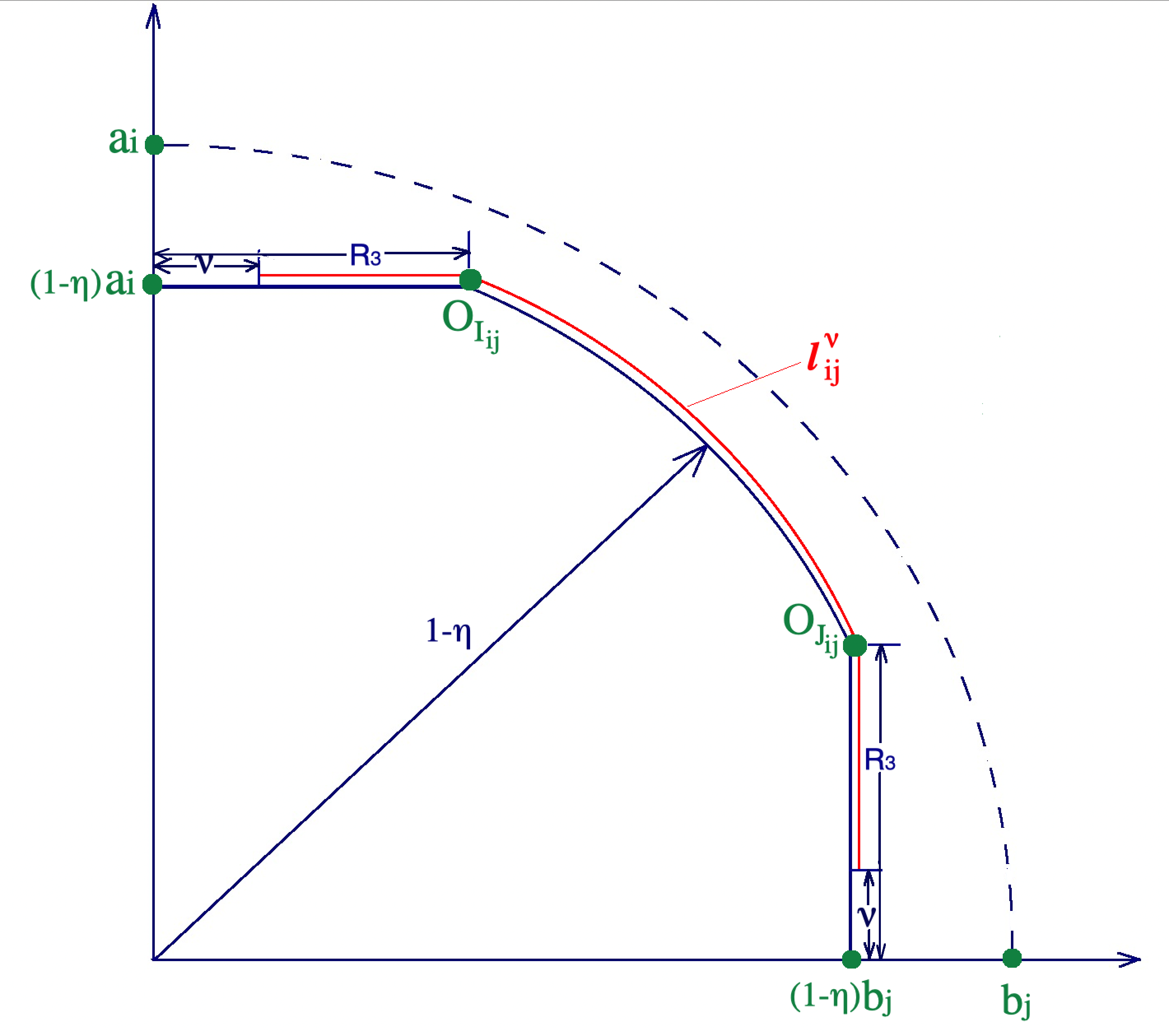}}
\nopagebreak[4]
\centerline{The 3-dimensional planar region $l^\nu_{ij}$}

\bigskip

Then it is easy to see that the map $\arg : l^\nu_{ij}\to \R/2\pi$ is injective. For $t\in [\nu, R_3]$, let $\theta(t)=\arctan\frac{t}{1-\sqrt\eta}$. Then $\{\arg z:z\in l^\nu_{ij}\}=[\theta(\nu), \frac\pi 2-\theta(\nu)]$, 
\be \arg [(1-\sqrt\eta)a_i+tb_j]=\theta(t), \arg[(1-\sqrt\eta)b_j+ta_i]=\frac\pi 2-\theta(t),\ee
and
\be \arg [(1-\eta)\cos\theta a_i+(1-\eta)\sin\theta b_j]=\theta.\ee

For any $\theta(\nu)\le\theta_1<\theta_2\le \frac\pi 2-\theta(\nu)$, let $l^{\theta_1,\theta_2}_{ij}$ denote the following subcurve of $l^\nu_{ij}$:
\be l^{\theta_1,\theta_2}_{ij}=\{z\in l^\nu_{ij}: \arg z\in [\theta_1,\theta_2]\}.\ee

For any map $f$ from $l^\nu_{ij}$ to $Q_{ij}$, let $G(f)$ denote its graph, and $G(f, \theta_1,\theta_2)$ denote the graph of the restriction of $f$ on $l^{\theta_1,\theta_2}_{ij}$: $G(f, \theta_1,\theta_2)=\{(x, f(x)):x\in l^{\theta_1,\theta_2}_{ij}\}$. Then $G(f)$ and $G(f, \theta_1,\theta_2)$ are subsets of $l^{\theta_1,\theta_2}_{ij}\times Q_{ij}$.

\medskip

So let $G$ be an $(\eta,\d,\nu,L)$-sliding boundary $G$ for $Z$. Then there exists a Lipschitz deformation $\varphi_G:\pa \cU\to \pa \cU$, so that $|\varphi_G(z)-z|<\d$ for all $z\in G$, and $G$ satisfies the condition \tb{(3.29)} in Definition \tb{3.10}. 

For $1\le i,j\le 3$, let $\s^G_{ij}$ denote $\varphi_{G*}(z_{ij})$ defined as before. Then 
\be \s^G_{ij}=\sum_{k\ne i,l\ne j}\varphi_{G*}(l_{kl}).\ee

Since $\s^G_{ij}=\varphi_{G*}(z_{ij})$, and $\varphi_G$ is a $\d$-sliding deformation, combine with \tb{(5.2)} we know that 
\be|\s^G_{ij}|\subset B(z_{ij}, \d)\subset \cup_{k\ne i,l\ne j}B(l_{kl}, \d).
\ee
Since \be B(l_{kl}, \d)\subset [A_{k2}\cup A_{1l}\cup \G_{kl}],\ee
we have
\be |\s^G_{ij}|\subset\cup_{k\ne i,l\ne j}[A_{k2}\cup A_{1l}\cup \G_{kl}].\ee

Note that for any $(k,l)\ne (i,j)$, since dist$(l_{kl}, L^\nu_{ij})\ge \nu>\d$, hence $|\varphi_G(l_{kl})|\cap L^\nu_{ij}\subset B(l_{kl}, \d)\cap L^\nu_{ij}=\emptyset$, and therefore

 \be G\cap L^\nu_{ij}=\varphi_G(l_{ij})\cap L^\nu_{ij}.\ee
 
 Let $\S^G_{ij}$ denote the unique $\Z_2$-2-chain in $Q_{ij}$ so that $\pa\S^G_{ij}=q_{ij*}(\s^G_{ij})$.
 
\begin{pro}Let $G$ be an $(\eta,\d,\nu, L)$-sliding boundary for $Z$ for some $\eta<\eta_1$, $\d<R_1$, $\nu\in (0,R_3)$, and $L>0$. Suppose that $\arctan\frac{R_3}{1-\sqrt\eta}<\frac{\pi}{16}$, $\d<\min\{\nu,\frac{R_3-\nu}{6}\}$, and $L<\frac{\d^2R_3}{(1-\eta)^2(1-\sqrt\eta)}$. Then there exists an $(\eta,\d,\nu,(1-\sqrt\eta)L/\d)$ sliding boundary $G'$ of $Z$, such that  $G'\cap L^{\frac{5\pi}{32},\frac\pi 2-\theta(R_3-3\d)}_{ij}$ is the image of an $(1-\sqrt\eta)L/\d$-Lipschitz graph from $l^{\frac{5\pi}{32},\theta(R_3-3\d)}_{ij}$ to $Q_{x_i}(0,\d)=[-\d x_i, \d x_i]$, and 
\be \sum_{1\le i\le 3}\sum_{1\le j\le 3}|q_{ij*}(\s^{G'}_{ij})|_2=\sum_{1\le i\le 3}\sum_{1\le j\le 3}|q_{ij*}(\s^G_{ij})|_2\ee
\end{pro}

\nd Since $G$ is an $(\eta,\d,\nu, L)$-sliding boundary for $Z$, by definition, for each $(i,j)$, $G\cap L^\nu_{ij}$ is the graph of an $L$-Lipschitz map $\xi^G_{ij}$ from $l^\nu_{ij}$ to $B_{Q_{ij}}(0,\d)$. That is,
\be G\cap L^\nu_{ij}=G(\xi^G_{ij}, \theta(\nu), \frac \pi2-\theta(\nu)).\ee

On the other hand, since $|\varphi_G(x)-x|<\d<\nu$, and the boundary $a_i$ and $b_j$ of $l_{ij}$ are of distance $\nu$ to $L^\nu_{ij}$, hence the boundary $\pa \varphi_{G*}(l_{ij})=\varphi_G(a_i)+\varphi_G(b_j)$ still does not touch $L^\nu_{ij}$. In other words, $\varphi_{G*}(l_{ij})$ admits no boundary in $L^\nu_{ij}$. Note that $|\varphi_{G*}(l_{ij})|\cap L^\nu_{ij}\subset \varphi_G(l_{ij})\cap L^\nu_{ij}$, and by \tb{(6.12)}, $\varphi_G(l_{ij})\cap L^\nu_{ij}$ is a Lipschitz graph of a curve, hence we know that 
\be|\varphi_{G*}(l_{ij})|\cap L^\nu_{ij}= \varphi_G(l_{ij})\cap L^\nu_{ij}.\ee

Combine with \tb{(6.12) and (6.13)}, we know that 
\be \varphi_{G*}(l_{ij})\cap L^\nu_{ij}\mbox{ is just the 1-chain represented by the graph }G(\xi_{ij}, \theta(\nu), \frac \pi2-\theta(\nu)).\ee

We want to construct another sliding boundary $G'$ based on $G$, so that the restriction of $G'$ on the part $L^{\frac{5\pi}{32},\theta(R_3-3\d)}_{ij}$ is a graph from $l^{\frac{5\pi}{32},\theta(R_3-3\d)}_{ij}$ to $Q_{x_i}$, instead of to both direction of $Q_{x_i\wg y_j}$. The idea is just to do some homotopy from $\xi^G_{ij}$ to $\pi_{x_i}\circ \xi^G_{ij}$, and guarantee that meanwhile \tb{(6.11)} holds.

So for $\theta\in [\frac\pi 2-\theta(R_3-3\d),\frac\pi 2-\theta(R_3-4\d)]$, let $t_\theta=\frac{\theta-(\frac\pi 2-\theta(R_3-3\d))}{\theta(R_3-3\d)-\theta(R_3-4\d)}$; and for $\theta\in [\frac\pi 8, \frac{5\pi}{32}]$, let $t_\theta=\frac{\theta-\frac\pi 8}{\frac {\pi}{32}}$. Then define $\xi^{G'}_{ij}: l^\nu_{ij}\to B_{Q_{ij}}(0,\d)$: 
\be \xi^{G'}_{ij}(z)=\left\{\begin{array}{rcl}
\xi^G_{ij}(z)&,\ &\arg z\in [\theta(\nu), \frac{\pi}{8}];\\
(1-t_{\arg z})\xi^G_{ij}(z)+t_{\arg z}\pi_{x_i}\circ \xi^G_{ij}(z)&,\ &\arg z\in [\frac{\pi}{8}, \frac{5\pi}{32}];\\
\pi_{x_i}\circ \xi^G_{ij}(z)&,\ &\arg z\in [\frac{5\pi}{32}, \frac\pi 2-\theta(R_3-3\d)];\\
(1-t_{\arg z})\pi_{x_i}\circ \xi^G_{ij}(z)+t_{\arg z}\xi^G_{ij}(z)&,\ &\arg z\in [\frac\pi 2-\theta(R_3-3\d), \frac\pi 2-\theta(R_3-4\d)];\\
\xi^G_{ij}(z)&,\ &\arg z\in [\frac\pi 2-\theta(R_3-4\d), \frac\pi 2-\theta(\nu)].
\end{array}\right.\ee
Then a simple calculate tells that $\xi^{G'}_{ij}$ is $\frac{L}{\min\{\frac{\pi}{32},\theta(R_3-3\d)-\theta(R_3-4\d) \}}$-Lipschitz. Note that 
\be\begin{split}\theta(R_3-3\d)-\theta(R_3-4\d)=\arctan(\frac{R_3-3\d}{1-\sqrt\eta})-\arctan(\frac{R_3-4\d}{1-\sqrt\eta})\ge (\frac{\d}{1-\sqrt\eta}),\end{split} \ee
Hence
$\xi^{G'}_{ij}$ is $\frac{L}{\min\{\frac{\pi}{32}, \frac{\d}{1-\sqrt\eta}\}}$-Lipschitz, i.e., $\max\{\frac{32}{\pi}, \frac{1-\sqrt\eta}{\d}\}L$-Lipschitz. Note that $\max\{\frac{32}{\pi}, \frac{1-\sqrt\eta}{\d}\}\le \frac {1-\sqrt\eta}{\d}$. Set $L'=\frac{(1-\sqrt\eta)L}{\d}$, then $\xi^{G'}_{ij}$ is $L'$-Lipschitz.

By definition, $\xi_{ij}^G=\xi^{G'}_{ij}$ on $L_{ij}^{\frac{\pi}{16}, \frac{\pi}{8}}\cup L_{ij}^{\frac\pi 2-\theta(R_3-4\d), \frac\pi 2-\theta(\nu)}$ forall $1\le i,j\le 3$.

Now we define the new sliding competitor $G'$: we let $G'=G$ on $\pa\cU\bs [\cup_{1\le i,j\le 3}L^\nu_{ij}]$, and $G'=G(\xi^{G'}_{ij})$ on $L^\nu_{ij}$, $1\le i,j\le 3$.

Then we can see from definition that $G'\cap L^{\frac{5\pi}{32},\frac\pi 2-\theta(R_3-3\d)}_{ij}$ is the image of the $L$-Lipschitz map $\pi_{x_i}\circ\xi_{ij}^G$ from $L^{\frac{5\pi}{32},\theta(R_3-3\d)}_{ij}$ to $Q_{x_i}(0,\d)=[-\d x_i, \d x_i]$. The fact that $G'$ is a $(\eta,\d)$-sliding boundary is due to the following : for any $z\in Z\cap \pa\cU$, if $\varphi_G(z)\not\in [\cup_{1\le i,j\le 3}L^\nu_{ij}]$, set $\varphi_{G'}(z)=\varphi_G(z)$; otherwise, we have $\varphi_G(z)=\xi^G_{ij}(y)$ for some $y\in l^\nu_{ij}$, and we set $\varphi_{G'}(z)=\xi^{G'}_{ij}(y)$. Note that in this case, since $\varphi_G(z)\in P_{ij}\times Q_{x_i}\times Q_{y_j}=Q_{a_i\wg b_j\wg x_i}\times Q_{y_j}$, we have
\be \begin{split}||\varphi_G(z)-z||^2&=||\pi_{a_i\wg b_j\wg x_i}(\varphi_G(z)-z)||^2+||\pi_{y_j}(\varphi_G(z)-z)||^2\\
&=||\pi_{a_i\wg b_j\wg x_i}(\varphi_G(z)-z)||^2+||\pi_{y_i}\varphi_G(z)||^2.
\end{split}\ee
and similarly we have
\be ||\varphi_{G'}(z)-z||^2=||\pi_{a_i\wg b_j\wg x_i}(\varphi_{G'}(z)-z)||^2+||\pi_{y_i}\varphi_{G'}(z)||^2.\ee
But note that by definition, for all $z\in l^\nu_{ij}$,
\be ||\pi_{y_i}\varphi_{G'}(z)||\le ||\pi_{y_i}\varphi_G(z)||,\ee
and 
\be ||\pi_{a_i\wg b_j\wg x_i}(\varphi_G(z)-z)||=||\pi_{a_i\wg b_j\wg x_i}(\varphi_{G'}(z)-z)||,\ee
hence by \tb{(6.17) and (6.18)} we have
\be ||\varphi_{G'}(z)-z||\le ||\varphi_G(z)-z||\le \d.\ee

Now the rest is to prove \tb{(6.11)}. Let us first compare $q_{ij*}(\s_{ij}^{G'})$ and $q_{ij*}(\s_{ij}^G)$. To give an idea, take $i=j=1$ for example. By definition, we know that both $\s^G_{11}$ and $\s^{G'}_{11}$ are contained in $B(z_{11},\d)$, hence $|q_{11*}(\s^G_{11})|$ and $|q_{11*}(\s^{G'}_{11})|$ are contained in $B(q_{11}(z_{11}),\d)$.

A simple calculate gives that 
$q_{11}(z_{11})$ is the boundary of the region $\{x\in Q_{11}:||x||\le \frac{\sqrt 3}{2}(1-\eta)$ and $|\lg x,x_1\rg |\le \frac{\sqrt 3}{2}(1-\sqrt\eta)\}.$ And hence the region $B(q_{11}(z_{11}),\d)$ is as in the following picture. 

\medskip

\centerline{\includegraphics[width=0.6\textwidth]{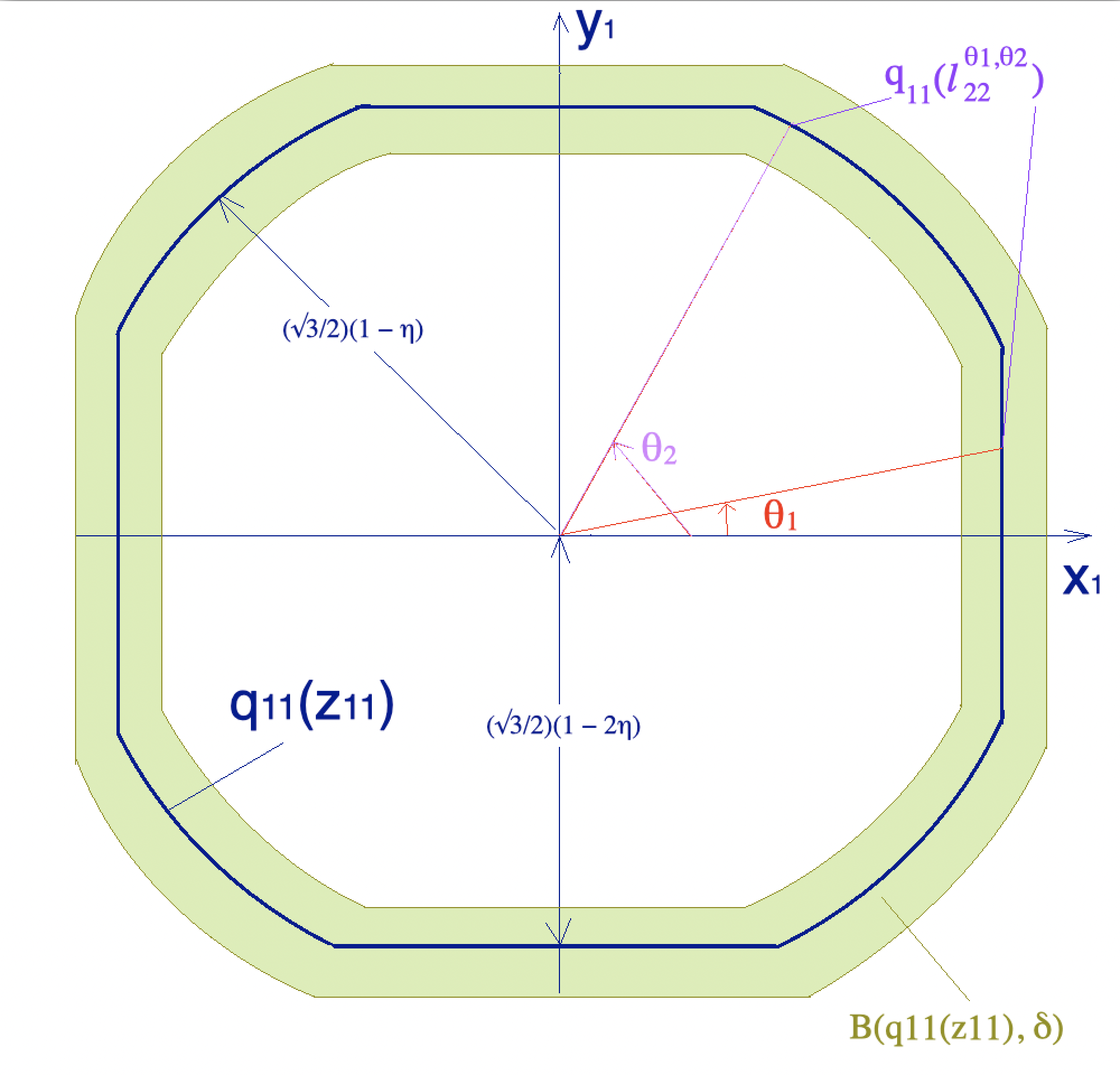}}

We also have that
\be q_{11}(l_{22}^{\theta_1,\theta_2})=\{z\in q_{11}(z_{11}):\arg z\in [\theta_1,\theta_2]\}.\ee
\be q_{11}(l_{23}^{\theta_1,\theta_2})=\{z\in q_{11}(z_{11}):\pi-\arg z\in [\theta_1,\theta_2]\}.\ee
\be q_{11}(l_{32}^{\theta_1,\theta_2})=\{z\in q_{11}(z_{11}):2\pi-\arg z\in [\theta_1,\theta_2]\}.\ee
\be q_{11}(l_{33}^{\theta_1,\theta_2})=\{z\in q_{11}(z_{11}):\arg z-\pi\in [\theta_1,\theta_2]\}.\ee
In all 
\be q_{11}(l_{kl}^{\theta_1,\theta_2})=\{z\in q_{11}(z_{11}):(-1)^{k+l}\arg z+l\pi\in [\theta_1,\theta_2]\}.\ee

See the figure above.

By definition, the difference of $\s^G_{11}$ and $\s^{G'}_{11}$ is contained in the disjoint union of $\cup_{k,l=2,3}l_{kl}^{\frac{\pi}{8}, \frac\pi 2-\theta(R_3-4\d)}\times B_{Q_{kl}}(0,\d)$, while
\be \begin{split}q_{11}(l_{kl}^{\frac{\pi}{8}, \frac\pi 2-\theta(R_3-4\d)}\times B_{Q_{kl}}(0,\d))=q_{11}(l_{kl}^{\frac{\pi}{8}, \frac\pi 2-\theta(R_3-4\d)})+q_{11}(B_{Q_{kl}}(0,\d))\\
=q_{11}(l_{kl}^{\frac{\pi}{8}, \frac\pi 2-\theta(R_3-4\d)})+B_{Q_{11}}(0,\frac\d 2)\subset 
q_{11}(l_{kl}^{\frac{3\pi}{32}, \frac\pi 2-\theta(R_3-5\d)})+[-\d y_1, \d y_1].
\end{split}\ee

As a result, we know that 
\be |q_{11*}(\s_{11}^G)|\Delta|q_{11*}(\s_{11}^{G'})|\subset \cup_{k,l=2,3}q_{11}(l_{kl}^{\frac{3\pi}{32}, \frac\pi 2-\theta(R_3-5\d)})+[-\d y_1, \d y_1],\ee
and the above union is disjoint.

On the other hand, for $k,l=2,3$, by \tb{(6.7)}, we know that 
\be \begin{split}&q_{11}(|\s_{11}^G|\bs [l_{kl}^{\frac{\pi}{16}, \frac\pi 2-\theta(\nu)}\times B_{Q_{kl}}(0,\d)])\subset \cup_{k',l'=2,3}
q_{11}(B(l_{k'l'}, \d)\bs [l_{k'l'}^{\frac{\pi}{16}, \frac\pi 2-\theta(\nu)}\times B_{Q_{kl}}(0,\d)])\\
&\subset[\cup_{k',l'=2,3, (k',l')\ne (k,l)}q_{11}(B(l_{k'l'}, \d))]\cup [q_{11}(B(l_{kl}, \d)\bs [l_{kl}^{\frac{\pi}{16}, \frac\pi 2-\theta(\nu)}\times B_{Q_{kl}}(0,\d)])]\\
&\subset [\cup_{k',l'=2,3, (k',l')\ne (k,l)}B(q_{11}(l_{k'l'}),\d)]\cup [B(q_{11}(l_{kl}\bs l_{kl}^{\frac{\pi}{16}, \frac\pi 2-\theta(\nu)}), \d)]\\
&\subset [\cup_{k',l'=2,3, (k',l')\ne (k,l)}B(q_{11}(l_{k'l'}),\d)]\cup [B(q_{11}(l_{kl}),\d)]\bs \{q_{11}(l_{kl}^{\frac{3\pi}{32}, \frac\pi 2-\theta(R_3-5\d)})+[-\d y_1, \d y_1]\}\\
&=B(q_{11}(z_{11}),\d)\bs \{q_{11}(l_{kl}^{\frac{3\pi}{32}, \frac\pi 2-\theta(R_3-5\d)})+[-\d y_1, \d y_1]\},
\end{split}\ee
hence
\be q_{11}(|\s_{11}^G|\bs [l_{kl}^{\frac{\pi}{16}, \frac\pi 2-\theta(\nu)}\times B_{Q_{kl}}(0,\d)])\cap \{q_{11}(l_{kl}^{\frac{3\pi}{32}, \frac\pi 2-\theta(R_3-5\d)})+[-\d y_1, \d y_1]\}=\emptyset.\ee

As a result, we have
\be q_{11}(|\s^G_{11}|)\cap \{q_{11}(l_{kl}^{\frac{3\pi}{32}, \frac\pi 2-\theta(R_3-5\d)})+[-\d y_1, \d y_1]\}\subset q_{11}[l_{kl}^{\frac{\pi}{16}, \frac\pi 2-\theta(\nu)}\times B_{Q_{kl}}(0,\d)],\ee
that is,
\be \begin{split}&q_{11}(|\s^G_{11}|)\cap \{q_{11}(l_{kl}^{\frac{3\pi}{32}, \frac\pi 2-\theta(R_3-5\d)})+[-\d y_1, \d y_1]\}\\
&\subset q_{11}(|\s^G_{11}|)\cap q_{11}[l_{kl}^{\frac{\pi}{16}, \frac\pi 2-\theta(\nu)}\times B_{Q_{kl}}(0,\d)]\\
&\subset q_{11}(|\s^G_{11}|\cap [l_{kl}^{\frac{\pi}{16}, \frac\pi 2-\theta(\nu)}\times B_{Q_{kl}}(0,\d)])\\
&=q_{11}(G\cap [l_{kl}^{\frac{\pi}{16}, \frac\pi 2-\theta(\nu)}\times B_{Q_{kl}}(0,\d)]),\end{split}\ee
the last equality is due to \tb{(6.10)}.

Similarly we have
\be q_{11}(|\s^{G'}_{11}|)\cap \{q_{11}(l_{kl}^{\frac{3\pi}{32}, \frac\pi 2-\theta(R_3-5\d)})+[-\d y_1,\d y_1]\}\subset q_{11}(G'\cap [l_{kl}^{\frac{\pi}{16}, \frac\pi 2-\theta(\nu)}\times B_{Q_{kl}}(0,\d)]).\ee

We claim that 
\be\mbox{the projection }\pi_{x_1}\mbox{ is injective on }q_{11}(G'\cap [l_{kl}^{\frac{\pi}{16}, \frac\pi 2-\theta(\nu)}\times B_{Q_{kl}}(0,\d)]).\ee

By definition, $G'\cap [l_{kl}^{\frac{\pi}{16}, \frac\pi 2-\theta(\nu)}\times B_{Q_{kl}}(0,\d)]=G(\xi^{G'}_{kl}, \frac{\pi}{16}, \frac\pi 2-\theta(\nu))$. Suppose that $z_1,z_2\in l_{kl}^{\frac{\pi}{16}, \frac\pi 2-\theta(\nu)}$, $\arg z_1<\arg z_2$. Then by definition, it is easy to that
\be \begin{split}||\pi_{x_1}(z_1)-\pi_{x_1}(z_2)||&\ge ||(1-\sqrt\eta)(\cos\arg z_2-\cos\arg z_1)||\\
&\ge (1-\sqrt\eta)(\arg z_2-\arg z_1)\sin \arg z_1\\
&\ge(1-\sqrt\eta)(\arg z_2-\arg z_1)\sin \theta(R_3)\\
&=\frac{(1-\sqrt\eta)R_3}{1-\eta}(\arg z_2-\arg z_1)
\end{split}\ee

On the other hand, since $\xi^{G'}$ is $L'$-Lipschitz, and $\pi_{x_1}$ is 1-Lipschitz, we know that
\be||\pi_{x_1}(\xi^{G'}(z_1))-\pi_{x_1}(\xi^{G'}(z_2))||\le L'||z_1-z_2||\le L'(1-\eta)|\arg z_1-\arg z_2|.\ee

As a result, we have
\be\begin{split} &||\pi_{x_1}(z_1, \xi^{G'}(z_1))-\pi_{x_1}(z_2, \xi^{G'}(z_2))||\\
\ge& ||\pi_{x_1}(z_1)-\pi_{x_1}(z_2)||-||\pi_{x_1}(\xi^{G'}(z_1)-\pi_{x_1}(\xi^{G'}(z_2)||\\
\ge& (\arg z_2-\arg z_1)\frac{(1-\sqrt\eta)R_3}{1-\eta}-L'(1-\eta)|\arg z_1-\arg z_2|\\
=&(\frac{(1-\sqrt\eta)R_3}{1-\eta}-L'(1-\eta))(\arg z_2-\arg z_1)>0
\end{split}\ee
because $(1-\eta)L'=\frac{(1-\eta)(1-\sqrt\eta)L}{\d}< \frac{(1-\sqrt\eta)R_3}{1-\eta}$. ($L<\frac{\d^2R_3}{(1-\eta)^2(1-\sqrt\eta)^2}$).

As a result, we get claim \tb{(6.34)}.

Similarly we have
\be\mbox{the projection }\pi_{x_1}\mbox{ is injective on }q_{11}(G\cap [l_{kl}^{\frac{\pi}{16}, \frac\pi 2-\theta(\nu)}\times B_{Q_{kl}}(0,\d)]).\ee

Now let us summerize:

We have two 1-chains $q_{11*}(\s^G_{11})$ and $q_{11*}(\s^{G'}_{11})$ in the plane $Q_{11}$, both are contained in $B(q_{11}(z_{11}), \d)$. And by \tb{(6.28)}, their difference are contained in the disjoint union
\be\begin{split} &\cup_{k,l=2,3}q_{11}(l_{kl}^{\frac{3\pi}{32}, \frac\pi 2-\theta(R_3-5\d)})+[-\d y_1, \d y_1]\\
&=\{z\in B(q_{11}(z_{11}), \d): |\lg z,x_1\rg |\in [\frac{\sqrt 3}{2}(R_3-5\d),\frac{\sqrt 3}{2}(1-\eta)\cos\frac{3\pi}{32}]\}.
\end{split}\ee

Moreover, by \tb{(6.31), (6.33), (6.34) and (6.38)}, we know that $\pi_{x_1}$ is injective on $|q_{11}(\s^G_{11})|\cap \{q_{11}(l_{kl}^{\frac{3\pi}{32}, \frac\pi 2-\theta(R_3-5\d)})+[-\d y_1, \d y_1]\}$, and on $|q_{11}(\s^{G'}_{11})|\cap \{q_{11}(l_{kl}^{\frac{3\pi}{32}, \frac\pi 2-\theta(R_3-5\d)})+[-\d y_1, \d y_1]\}$, $k,l=2,3$.


Let $\z_{11}^{kl,G}$ be the map from $ (-1)^k[\frac{\sqrt 3}{2}(R_3-5\d)x_1,\frac{\sqrt 3}{2}(1-\eta)\cos\frac{3\pi}{32}x_1]$ to $Q_{y_1}$, so that $|q_{11}(\s^G_{11})|\cap q_{11}(l_{kl}^{\frac{3\pi}{32}, \frac\pi 2-\theta(R_3-5\d)})+[-\d y_1, \d y_1]$ coincides with the graph of $\z_{11}^{kl,G}$. Define $\z_{11}^{kl,G'}$ similarly.

As a result, we know that $\S^{G'}_{ij}$ and $\S^{G'}_{ij}$ are the same outside the set
\be \{z\in Q_{11}: |\lg z,x_1\rg |\in [\frac{\sqrt 3}{2}(R_3-5\d),\frac{\sqrt 3}{2}(1-\eta)\cos\frac{3\pi}{32}]\}.\ee
which gives that
\be \begin{split}&|q_{11*}(\s^{G'}_{11})|_2-|q_{11*}(\s^G_{11})|_2=|\S^{G'}_{11}|-|\S^G_{11}|\\
=&|\S^{G'}_{11}\cap\{z\in Q_{11}: |\lg z,x_1\rg |\in [\frac{\sqrt 3}{2}(R_3-5\d),\frac{\sqrt 3}{2}(1-\eta)\cos\frac{3\pi}{32}]\}|\\
&-|\S^G_{11}\cap \{z\in Q_{11}: |\lg z,x_1\rg |\in [\frac{\sqrt 3}{2}(R_3-5\d),\frac{\sqrt 3}{2}(1-\eta)\cos\frac{3\pi}{32}]\}|
\end{split}\ee
 
 Note that the boundary of $\S^G_{11}\cap \{z\in Q_{11}: \lg z,x_1\rg \in [\frac{\sqrt 3}{2}(R_3-5\d),\frac{\sqrt 3}{2}(1-\eta)\cos\frac{3\pi}{32}]\}$ is the disjoint union of $q_{11*}(\s^G_{22})\cap [q_{11}(l_{22}^{\frac{3\pi}{32}, \frac\pi 2-\theta(R_3-5\d)})+[-\d y_1, \d y_1]]$, $q_{11*}(\s^G_{23})\cap [q_{11}(l_{23}^{\frac{3\pi}{32}, \frac\pi 2-\theta(R_3-5\d)})+[-\d y_1, \d y_1]]$, and some segments parallel to $y_1$. The parts $q_{11*}(\s^G_{22})\cap [q_{11}(l_{22}^{\frac{3\pi}{32}, \frac\pi 2-\theta(R_3-5\d)})+[-\d y_1, \d y_1]]$ and $q_{11*}(\s^G_{23})\cap [q_{11}(l_{23}^{\frac{3\pi}{32}, \frac\pi 2-\theta(R_3-5\d)})+[-\d y_1, \d y_1]]$ are two disjoint graphs $\z_{11}^{22,G}$ and $\z_{11}^{23,G}$ from $[\frac{\sqrt 3}{2}(R_3-5\d)x_1,\frac{\sqrt 3}{2}(1-\eta)\cos\frac{3\pi}{32}x_1]$, with $\z_{11}^{22,G}(z)>\z_{11}^{23,G}(z)$, for any $z\in [\frac{\sqrt 3}{2}(R_3-5\d)x_1,\frac{\sqrt 3}{2}(1-\eta)\cos\frac{3\pi}{32}x_1]$. Thus we have
 \be \begin{split}|\S^G_{ij}\cap \{z\in Q_{11}: \lg z,x_1\rg \in [\frac{\sqrt 3}{2}(R_3-5\d),\frac{\sqrt 3}{2}(1-\eta)\cos\frac{3\pi}{32}]\}|\\
 =\int_{\frac{\sqrt 3}{2}(R_3-5\d)}^{\frac{\sqrt 3}{2}(1-\eta)\cos\frac{3\pi}{32}}(\z_{11}^{22,G}(tx_1)-\z_{11}^{23,G}(tx_1))dt.
 \end{split}\ee
 Similarly we have
  \be\begin{split} |\S^G_{ij}\cap \{z\in Q_{11}: \lg z,x_1\rg \in -[\frac{\sqrt 3}{2}(R_3-5\d),\frac{\sqrt 3}{2}(1-\eta)\frac{3\pi}{32}]\}|\\
  =\int_{-\frac{\sqrt 3}{2}(1-\eta)\cos\frac{3\pi}{32}}^{-\frac{\sqrt 3}{2}(R_3-5\d)}(\z_{11}^{32,G}(tx_1)-\z_{11}^{33,G}(tx_1))dt,
 \end{split}\ee
\be \begin{split}|\S^{G'}_{ij}\cap \{z\in Q_{11}: \lg z,x_1\rg \in [\frac{\sqrt 3}{2}(R_3-5\d),\frac{\sqrt 3}{2}(1-\eta)\frac{3\pi}{32}]\}|\\
 =\int_{\frac{\sqrt 3}{2}(R_3-5\d)}^{\frac{\sqrt 3}{2}(1-\eta)\cos\frac{3\pi}{32}}(\z_{11}^{22,G'}(tx_1)-\z_{11}^{23,G'}(tx_1))dt,
 \end{split}\ee
 and 
 \be\begin{split} |\S^{G'}_{ij}\cap \{z\in Q_{11}: \lg z,x_1\rg \in -[\frac{\sqrt 3}{2}(R_3-5\d),\frac{\sqrt 3}{2}(1-\eta)\frac{3\pi}{32}]\}|\\
  =\int_{-\frac{\sqrt 3}{2}(1-\eta)\cos\frac{3\pi}{32}}^{-\frac{\sqrt 3}{2}(R_3-5\d)}(\z_{11}^{32,G'}(tx_1)-\z_{11}^{33,G'}(tx_1))dt.
 \end{split}\ee

 Then by \tb{(6.41)} we have
 \be \begin{split}&|q_{11*}(\s^{G'}_{11})|_2-|q_{11*}(\s^G_{11})|_2\\
 =&
\int_{\frac{\sqrt 3}{2}(R_3-5\d)}^{\frac{\sqrt 3}{2}(1-\eta)\cos\frac{3\pi}{32}}[\z_{11}^{22,G'}(tx_1)-\z_{11}^{22,G}(tx_1)]-[\z_{11}^{23,G'}(tx_1)-\z_{11}^{23,G}(tx_1)]dt\\
&+\int_{-\frac{\sqrt 3}{2}(1-\eta)\cos\frac{3\pi}{32}}^{-\frac{\sqrt 3}{2}(R_3-5\d)}[\z_{11}^{32,G'}(tx_1)-\z_{11}^{32,G}(tx_1)]-[\z_{11}^{33,G'}(tx_1)-\z_{11}^{33,G}(tx_1))]dt.\end{split}\ee

Similar discuss gives
 \be \begin{split}&|q_{12*}(\s^{G'}_{12})|_2-|q_{12*}(\s^G_{12})|_2\\
 =&
\int_{\frac{\sqrt 3}{2}(R_3-5\d)}^{\frac{\sqrt 3}{2}(1-\eta)\frac{3\pi}{32}}[\z_{12}^{23,G'}(tx_1)-\z_{12}^{23,G}(tx_1)]-[\z_{12}^{21,G'}(tx_1)-\z_{12}^{21,G}(tx_1)]dt\\
&+\int_{-\frac{\sqrt 3}{2}(1-\eta)\cos\frac{3\pi}{32}}^{-\frac{\sqrt 3}{2}(R_3-5\d)}[\z_{12}^{33,G'}(tx_1)-\z_{12}^{33,G}(tx_1)]-[\z_{12}^{31,G'}(tx_1)-\z_{12}^{31,G}(tx_1))]dt,\end{split}\ee
and 
 \be \begin{split}&|q_{13*}(\s^{G'}_{13})|_2-|q_{13*}(\s^G_{13})|_2\\
 =&
\int_{\frac{\sqrt 3}{2}(R_3-5\d)}^{\frac{\sqrt 3}{2}(1-\eta)\cos\frac{3\pi}{32}}[\z_{13}^{21,G'}(tx_1)-\z_{13}^{21,G}(tx_1)]-[\z_{13}^{22,G'}(tx_1)-\z_{13}^{22,G}(tx_1)]dt\\
&+\int_{-\frac{\sqrt 3}{2}(1-\eta)\cos\frac{3\pi}{32}}^{-\frac{\sqrt 3}{2}(R_3-5\d)}[\z_{13}^{31,G'}(tx_1)-\z_{13}^{31,G}(tx_1)]-[\z_{13}^{32,G'}(tx_1)-\z_{13}^{32,G}(tx_1))]dt,\end{split}\ee

Let us now look at the terms $\z_{1j}^{kl, G}$, for $k\ne 1, l\ne j$. For each pair $(k,l)$ with $k\ne 1$, take $k=l=2$ for example, it concerns of two terms on $G$ in \tb{(6.46)-(6.48)}: $\z_{11}^{22, G}$, and $\z_{13}^{22, G}$. Note that they are part of projections of $G(\xi^G_{22}, \frac{\pi}{16}, \frac\pi2-\theta(\nu))$ under $q_{11}$ and $q_{13}$ respectively. So take any $t\in [\frac{\sqrt 3}{2}(R_3-5\d),\frac{\sqrt 3}{2}(1-\eta)\frac{3\pi}{32}]$, we know that $(tx_1, \z_{11}^{22, G}(t)y_1)$ is the image of a point $z_t$ in $G(\xi^G_{22}, \frac{\pi}{16}, \frac\pi2-\theta(\nu))$ under $q_{11}$. And hence 
\be\begin{split} tx_1=\pi_{x_1}((tx_1, \z_{11}^{22, G}(t)y_1))=\pi_{x_1}(q_{11}(z_t))=\pi_{x_1}(z_t)\\
=\pi_{x_1}(q_{13}(z_t))\in \pi_{x_1}(q_{13}(G(\xi^G_{22}, \frac{\pi}{16}, \frac\pi2-\theta(\nu)))).\end{split}\ee
Note that $\pi_{x_1}$ is injective on $q_{13}(G(\xi^G_{22}, \frac{\pi}{16}, \frac\pi2-\theta(\nu)))$, hence we know that
$(tx_1, \z_{11}^{23, G}(t)y_1)$ is the image of the same point $z_t$ under $q_{13}$.

Since $z_t\in G(\xi^G_{22}, \frac{\pi}{16}, \frac\pi2-\theta(\nu))$, there exists a unique $w_t\in l_{22}^{\frac{\pi}{16}, \frac\pi2-\theta(\nu)}$, so that $z_t=(w_t, \xi^G_{22}(w_t))\in P_{22}\times Q_{22}$.  

Then by definition, 
\be \z_{11}^{22, G}(tx_1)=\lg z_t, y_1\rg =\lg w_t,y_1\rg +\lg \xi^G_{22}(w_t),y_1\rg ,\ee
and 
\be \z_{13}^{22, G}(tx_1)=\lg z_t, y_3\rg =\lg w_t,y_3\rg +\lg \xi^G_{22}(w_t),y_3\rg .\ee
Since $w_t\in P_{22}=Q_{a_2, b_2}$, we know that $\lg w_t, y_1\rg +\lg w_t, y_3\rg =0$. Since $\xi^G_{22}(w_t)\in Q_{x_2, y_2}$, we know that $\lg \xi^G_{22}(w_t),y_1\rg =\lg \xi^G_{22}(w_t),y_3\rg $. As a result, we have
\be \z_{11}^{22, G}(tx_1)-\z_{13}^{22, G}(tx_1)=\lg 2w_t, y_1\rg .\ee

On the other hand, for the point $w_t$, we know that $\xi^G_{22}(w_t)-\xi^{G'}_{22}(w_t)\in Q_{y_2}$, hence by the same calculation as in \tb{(6.49)}, we know that
\be tx_1=\pi_{x_1}(w_t, \xi^{G'}_{22}(w_t)),\ee
that is, $\z_{11}^{22, G}(tx_1)$ and $\z_{11}^{22, G'}(tx_1)$ corresponds to the same point $w_t\in l_{22}^{\frac{\pi}{16}, \frac\pi2-\theta(\nu)}$. Hence by the same argument as above, we know that
\be \z_{11}^{22, G'}(tx_1)-\z_{13}^{22, G'}(tx_1)=\lg 2w_t, y_1\rg .\ee

Thus, \tb{(6.52) and (6.54)} gives
\be [\z_{11}^{22, G'}(tx_1)-\z_{11}^{22, G}(tx_1)]-[\z_{13}^{22, G'}(tx_1)-\z_{13}^{22, G}(tx_1)]=0.\ee

Same arguments gives that, for each $k\ne 1, 1\le l\le 3$, and for $i,j\ne l$, we have
\be [\z_{1i}^{kl, G'}(tx_1)-\z_{1i}^{kl, G}(tx_1)]-[\z_{1j}^{kl, G'}(tx_1)-\z_{1j}^{kl, G}(tx_1)]=0.\ee 

We sum over \tb{(6.46)-(6.48)}, taking the relation \tb{(6.56)} into account, and get that
\be \sum_{1\le j\le 3}|q_{1j*}(\s^{G'}_{1j})|_2-|q_{1j*}(\s^G_{1j})|_2=0.\ee

Similarly argument gives, for $1\le i\le 3$, 
\be \sum_{1\le j\le 3}|q_{ij*}(\s^{G'}_{ij})|_2-|q_{1j*}(\s^G_{ij})|_2=0.\ee

This gives the relation \tb{(6.11)}.\qed

\begin{rem}Note that here the fastest way is to project the whole $\xi_{ij}^G$ to $Q_{x_i}$ to get a new sliding boundary $G'$ which is a graph from the whole $l^\nu_{ij}$ to $Q_{x_i}$. But note that this cannot help, even if we only do it on the arc $(1-\eta)\gamma^0_{ij}$, because for any point $z\in l^\nu_{ij}$ with $\arg z<\theta(R_3)$, the injectivity condition as \tb{(6.34)} cannot hold, no matter how small $L$ is, because $\pi_{x_1}$ is not even injective on $l^{\theta(\nu), \theta(R_3)}_{ij}$, hence the graph cannot be injective on a neighborhood of this segment. 

There might be problem if we loose injectivity. For example 

That is why we have to do the next step, to get a sliding boundary which coincide with $z_{ij}$ on $L^{R_3}_{ij}$. See the following propositions.\end{rem}

\begin{pro}Let $G$ be an $(\eta,\d,\nu, L)$-sliding boundary for $Z$ for some $\eta<\eta_1$, $\d<R_1$, $\nu\in (0,R_3)$, and $L>0$. Suppose that $\arctan\frac{R_3}{1-\sqrt\eta}<\frac{\pi}{16}$, $\d<\min\{\nu,\frac{R_3-\nu}{6}\}$, and $L<\frac{R_3\d^2}{(1-\sqrt\eta)(1-\eta)^2}$. Then there exists an $(\eta,\d,\nu,(1-\sqrt\eta)^2L/\d^2)$-sliding boundary $G''$ of $Z$, such that  

$1^\circ$ $G''\cap L^{\theta(R_3-3\d),\frac\pi 2-\theta(R_3)}_{ij}$ is the image of an $(1-\sqrt\eta)^2L/\d^2$-Lipschitz graph from $l^{\theta(R_3-3\d),\frac\pi 2-\theta(R_3)}_{ij}$ to $Q_{y_j}(0,\d)=[-\d y_j, \d y_j]$;

$2^\circ$ $G''\cap L^{\frac{\pi}{4},\frac\pi 2-\theta(R_3)}_{ij}$ coincides with $Z\cap L^{\frac{\pi}{4},\frac\pi 2-\theta(R_3)}_{ij}$;

$3^\circ$
\be \sum_{1\le i\le 3}\sum_{1\le j\le 3}|q_{ij*}(\s^{G''}_{ij})|_2=\sum_{1\le i\le 3}\sum_{1\le j\le 3}|q_{ij*}(\s^G_{ij})|_2.\ee
\end{pro}

\nd Fix any $G$. Take the corresponding $G'$ as in the last proposition. We will do the same projection for $G'\cap L^\nu_{ij}$ as for $G$, but this time to the direction $Q_{y_j}$.

Recall that $G'$ is a $(\eta,\d,\nu, L')$-sliding boundary for $Z$, with $\arctan\frac{R_3}{1-\sqrt\eta}<\frac{\pi}{16}$, $\d<\min\{\nu,\frac{R_3-\nu}{6}\}$, $L'=\frac{(1-\sqrt\eta)L}{\d}$, and $L<\frac{R_3\d^2}{(1-\sqrt\eta)(1-\eta)^2}$. Moreover, 
\be \begin{split}G'\cap L^{\frac{5\pi}{32},\frac\pi 2-\theta(R_3-3\d)}_{ij}\mbox{ is the image of an }L\mbox{-Lipschitz graph from }\\
l^{\frac{5\pi}{32},\frac\pi 2-\theta(R_3-3\d)}_{ij}\mbox{ to }Q_{x_i}(0,\d)=[-\d x_i, \d x_i].\end{split}\ee

Take all the notations as in the last proposition. For $\theta\in [\frac\pi 2-\theta(R_3),\frac\pi 2-\theta(R_3-\d)]$, let $t_\theta=\frac{\theta-(\frac\pi 2-\theta(R_3))}{\theta(R_3)-\theta(R_3-\d)}$; and for $\theta\in [\theta(R_3-4\d), \theta(R_3-3\d)]$, let $t_\theta=\frac{\theta-\theta(R_3-4\d)}{\theta(R_3-3\d)-\theta(R_3-4\d)}$. 
 We define, for $1\le i,j\le 3$, 
$\xi^{G''}_{ij}: l^\nu_{ij}\to B_{Q_{ij}}(0,\d)$: 
\be \xi^{G''}_{ij}(z)=\left\{\begin{array}{rcl}
\xi^{G'}_{ij}(z)&,\ &\arg z\in [\theta(\nu), \theta(R_3-4\d)];\\
(1-t_{\arg z})\xi^{G'}_{ij}(z)+t_{\arg z}\pi_{y_j}\circ \xi^{G'}_{ij}(z)&,\ &\arg z\in [\theta(R_3-4\d), \theta(R_3-3\d)];\\
\pi_{y_j}\circ \xi^{G'}_{ij}(z)&,\ &\arg z\in [\theta(R_3-4\d), \frac\pi 2-\theta(R_3)];\\
(1-t_{\arg z})\pi_{y_j}\circ \xi^{G'}_{ij}(z)+t_{\arg z}\xi^{G'}_{ij}(z)&,\ &\arg z\in [\frac\pi 2-\theta(R_3), \frac\pi 2-\theta(R_3-\d)];\\
\xi^{G'}_{ij}(z)&,\ &\arg z\in [\frac\pi 2-\theta(R_3-\d), \frac\pi 2-\theta(\nu)].
\end{array}\right.\ee

Then similar as in Proposition \tb{6.2}, $1^\circ$ and $2^\circ$ follows directly from the definition of $G''$. The rest is to prove \tb {(6.59)}.

Note that this time we do the projection on $y_j$. The main difference here is that for all points $z$ in $l_{ij}$ with $\arg z\in [\frac \pi2-\theta(R_3), \frac\pi 2-\nu]$, we know that $\pi_{y_j}(z)$ are the same, and hence we cannot expect that $\pi_{y_i}$ is injective on either $G(\xi^{G'}_{ij}, \frac \pi2-\theta(R_3), \frac\pi 2-\nu)$, or $G(\xi^{G''}_{ij}, \frac \pi2-\theta(R_3), \frac\pi 2-\nu)$. That is, the injectivity that we used cannot apply here. But here, compare to an arbitrary sliding boundary, the set $G'$ has a special property \tb{(6.60)}. This property guarantees that $\pi_{y_j}(z, \xi^{G'}_{ij}(z))=\pi_{y_j}(z)$ for all $z\in l_{ij}^{\frac{5\pi}{32}, \frac\pi2-\theta(R_3-3\d)}$.

Again let us look at $q_{11*}$ for example. By \tb{(6.60)}, we know that, for any $\frac{5\pi}{32}\le\theta_1<\theta_2\le \frac\pi 2-\theta(R_3-3\d)$, 
\be G(\xi^{G'}_{kl},\theta_1, \theta_2)\subset l_{kl}^{\theta_1, \theta_2}+[-\d x_1,\d x_1],\ee and 
\be G(\xi^{G''}_{kl},\theta_1, \theta_2)\subset l_{kl}^{\theta_1, \theta_2}+[-\d x_1,\d x_1].\ee

Since $\xi^{G'}_{kl}$ and $\xi^{G''}_{kl}$ are different only in $l_{kl}^{\theta(R_3-4\d), \frac\pi 2-\theta(R_3-\d)}$, $k,l=2,3$, by \tb{(6.62) and (6.63)}, hence we know that
\be \begin{split}|\xi^{G'}_{kl}|\Delta|\xi^{G''}_{kl}|\subset &\cup_{k,l=2,3}[G(\xi^{G'}_{kl},\theta(R_3-4\d), \frac\pi 2-\theta(R_3-\d))\cup G(\xi^{G''}_{kl},\theta(R_3-4\d), \frac\pi 2-\theta(R_3-\d))]\\
\subset &\cup_{k,l=2,3}[G(\xi^{G'}_{kl},\theta(R_3-4\d), \frac{7\pi}{32})\cup G(\xi^{G'}_{kl},\frac{7\pi}{32},\frac\pi 2-\theta(R_3-\d))\\
&\cup G(\xi^{G''}_{kl},\theta(R_3-4\d), \frac{7\pi}{32})\cup G(\xi^{G''}_{kl},\frac{7\pi}{32},\frac\pi 2-\theta(R_3-\d))]\\
\subset &\cup_{k,l=2,3}\{[l_{kl}^{\frac{7\pi}{32}, \frac\pi 2-\theta(R_3-\d)}+[-\d x_1,\d x_1]]\cup [l_{kl}^{\theta(R_3-4\d), \frac{7\pi}{32}}\times B_{Q_{11}}(0,\d)]\}
\end{split}\ee

As a result,
\be \begin{split}&|q_{11*}(\s_{11}^{G'})|\Delta |q_{11*}(\s_{11}^{G''})|\\
&\subset \cup_{k,l=2,3}q_{11}(\{[l_{kl}^{\frac{7\pi}{32}, \frac\pi 2-\theta(R_3-\d)}+[-\d x_1,\d x_1]]\cup [l_{kl}^{\theta(R_3-4\d), \frac{7\pi}{32}}\times B_{Q_{11}}(0,\d)]\})\\
&\subset\cup_{k,l=2,3}[q_{11}(l_{kl}^{\theta(R_3-5\d), \frac\pi 2-\theta(R_3-\d)})+[-\d x_1,\d x_1]],
\end{split}\ee
where the above union is disjoint.

%

On the other hand, still by \tb{(6.60)} (for the part $[\frac{6\pi}{32},\frac\pi 2-\theta(R_3)]$ and the same argument as in Proposition \tb{6.2} (for the part $[\theta(R_3-5\d), \frac{6\pi}{32}]$), we know that for $k\ne 2, l\ne 2$, $\pi_{y_1}$ is injective on
\be |q_{11*}(\s^{G''}_{11})|\cap \{q_{11}(l_{kl}^{\theta(R_3-5\d), \frac\pi 2-\theta(R_3)})+[-\d x_1,\d x_1]\},\ee
and 
\be\pi_{y_1}(|q_{11*}(\s^{G''}_{11})|\cap \{q_{11}(l_{kl}^{\theta(R_3-5\d), \frac\pi 2-\theta(R_3)})+[-\d x_1,\d x_1])=(-1)^k[\frac{\sqrt 3}{2}(R_3-5\d),\frac{\sqrt 3}{2}(1-\sqrt\eta)].\ee

And also note that for $z$ so that $\arg z\in [\frac\pi2-\theta(R_3), \frac\pi2-\theta(R_3-3\d)]$, their projections under $\pi_{y_1}$ are $1-\sqrt\eta$, and both $q_{11}(G(\xi^{G'}_{kl},\frac\pi2-\theta(R_3), \frac\pi2-\theta(R_3-3\d)))$ and $q_{11}(G(\xi^{G''}_{kl},\frac\pi2-\theta(R_3), \frac\pi2-\theta(R_3-3\d)))$ are contained in the line $\pi_{y_1}^{-1}\{1-\sqrt\eta)$, hence their difference makes no contribution to $|\S_{11}^{G'}|\Delta |\S_{11}^{G''}|$. 

As before we let $\lambda_{11}^{kl,G'}$ be the map from $(-1)^k[\frac{\sqrt 3}{2}(R_3-5\d),\frac{\sqrt 3}{2}(1-\sqrt\eta)]$ to $Q_{x_1}$, so that $|q_{11*}(\s^{G'}_{11})|\cap \{q_{11}(l_{kl}^{\frac{6\pi}{32}, \frac\pi 2-\theta(R_3)})+[-\d x_1,\d x_1]\}$ coincides with the graph of $\lambda_{11}^{kl,G}$. Define $\lambda_{11}^{kl,G''}$ similarly.

Then the same argument as in the proof of Proposition \tb{6.2} gives
\be \begin{split}&|q_{11*}(\s_{11}^{G'})|_2-|q_{11*}(\S_{11}^{G''})|_2\\
=&\sum_{k,l=2,3}\int_{(-1)^k[\frac{\sqrt 3}{2}(1-\eta)\sin\frac{6\pi}{32},\frac{\sqrt 3}{2}(1-\sqrt\eta)]}(-1)^l(\lambda_{11}^{kl,G'}(ty_1)-\lambda_{11}^{kl,G''}(ty_1))dt\\
=&\int_{\frac{\sqrt 3}{2}(1-\eta)\sin\frac{6\pi}{32}}^{\frac{\sqrt 3}{2}(1-\sqrt\eta)}[\lambda_{11}^{22,G'}(ty_1)-\lambda_{11}^{22,G''}(ty_1)]-[\lambda_{11}^{23,G'}(ty_1)-\lambda_{11}^{23,G''}(ty_1)]dt\\
&+\int_{-\frac{\sqrt 3}{2}(1-\sqrt\eta)}^{-\frac{\sqrt 3}{2}(1-\eta)\sin\frac{6\pi}{32}}[\lambda_{11}^{32,G'}(ty_1)-\lambda_{11}^{32,G''}(ty_1)]-[\lambda_{11}^{33,G'}(ty_1)-\lambda_{11}^{33,G''}(ty_1))]dt.\end{split}\ee

and similarly, we have
\be \begin{split}&|q_{21*}(\s_{21}^{G'})|_2-|q_{21*}(\S_{21}^{G''})|_2\\
=&\int_{\frac{\sqrt 3}{2}(1-\eta)\sin\frac{6\pi}{32}}^{\frac{\sqrt 3}{2}(1-\sqrt\eta)}[\lambda_{21}^{32,G'}(ty_1)-\lambda_{21}^{32,G''}(ty_1)]-[\lambda_{21}^{12,G'}(ty_1)-\lambda_{21}^{12,G''}(ty_1)]dt\\
&+\int_{-\frac{\sqrt 3}{2}(1-\sqrt\eta)}^{-\frac{\sqrt 3}{2}(1-\eta)\sin\frac{6\pi}{32}}[\lambda_{21}^{33,G'}(ty_1)-\lambda_{21}^{33,G''}(ty_1)]-[\lambda_{21}^{13,G'}(ty_1)-\lambda_{21}^{13,G''}(ty_1))]dt,\end{split}\ee
and
\be \begin{split}&|q_{31*}(\s_{31}^{G'})|_2-|q_{31*}(\S_{31}^{G''})|_2\\
=&\int_{\frac{\sqrt 3}{2}(1-\eta)\sin\frac{6\pi}{32}}^{\frac{\sqrt 3}{2}(1-\sqrt\eta)}[\lambda_{31}^{12,G'}(ty_1)-\lambda_{31}^{12,G''}(ty_1)]-[\lambda_{31}^{22,G'}(ty_1)-\lambda_{31}^{22,G''}(ty_1)]dt\\
&+\int_{-\frac{\sqrt 3}{2}(1-\sqrt\eta)}^{-\frac{\sqrt 3}{2}(1-\eta)\sin\frac{6\pi}{32}}[\lambda_{31}^{13,G'}(ty_1)-\lambda_{31}^{13,G''}(ty_1)]-[\lambda_{31}^{23,G'}(ty_1)-\lambda_{31}^{23,G''}(ty_1))]dt.\end{split}\ee

Again, the same argument as in the proof of Proposition \tb{6.2} gives
\be \sum_{1\le j\le 3}|q_{j1*}(\s^{G'}_{j1})|_2-|q_{j1*}(\s^{G''}_{j1*})|_2=0.\ee

Similarly argument gives, for $1\le i\le 3$, 
\be \sum_{1\le j\le 3}|q_{ji*}(\s^{G'}_{ji})|_2-|q_{ji*}(\s^{G''}_{ji})|_2=0.\ee

We sum over $j$, and get
\be \sum_{1\le i\le 3}\sum_{1\le j\le 3}|q_{ij*}(\s^{G''}_{ij})|_2=\sum_{1\le i\le 3}\sum_{1\le j\le 3}|q_{ij*}(\s^{G'}_{ij})|_2.\ee
Combine with \tb{(6.11)}, we get \tb{(6.59)}.\qed

\begin{pro}Let $G$ be an $(\eta,\d,\nu, L)$-sliding boundary for $Z$ for some $\eta<\eta_1$, $\d<R_1$, $\nu\in (0,R_3)$, and $L>0$. Suppose that $\arctan\frac{R_3}{1-\sqrt\eta}<\frac{\pi}{16}$, $\d<\min{\nu,\frac{R_3-\nu}{6}}$, and $L<\frac{R_3\d^2}{(1-\sqrt\eta)(1-\eta)^2}$. Then there exists an $(\eta,\d,\nu,(1-\sqrt\eta)^2L/\d^2)$-sliding boundary $G_0$ of $Z$, such that  

$1^\circ$ $G_0\cap L^{\theta(R_3-3\d),\frac\pi 4}_{ij}$ is the image of an $(1-\sqrt\eta)^2L/\d^2)$-Lipschitz graph from $l^{\theta(R_3-3\d),\frac\pi 4}_{ij}$ to $Q_{y_j}(0,\d)=[-\d y_j, \d y_j]$;

$2^\circ$ $G_0\cap L^{\frac\pi 4,\frac\pi 2-\theta(R_3-3\d)}_{ij}$ is the image of an $(1-\sqrt\eta)^2L/\d^2)$-Lipschitz graph from $l^{\frac\pi 4,\frac\pi 2-\theta(R_3-3\d)}_{ij}$ to $Q_{x_i}(0,\d)=[-\d x_i, \d x_i]$;

$3^\circ$ $G_0\cap L^{\theta(R_3),\frac\pi 2-\theta(R_3)}_{ij}$ coincides with $\pa Z\cap L^{\theta(R_3),\frac\pi 2-\theta(R_3)}_{ij}$;

$4^\circ$
\be \sum_{1\le i\le 3}\sum_{1\le j\le 3}|q_{ij*}(\s^{G_0}_{ij})|_2=\sum_{1\le i\le 3}\sum_{1\le j\le 3}|q_{ij*}(\s^G_{ij})|_2\ee
\end{pro}

\nd It is enough to take the $G''$ as in Proposition \tb{6.4}, and project $G''\cap L_{ij}^{\theta(R_3), \frac\pi 2-\theta(R_3)}$ again to $l_{ij}^{\theta(R_3), \frac\pi 4}\times Q_{x_i}$, via a homotopy in $[\theta(R_3-\d), \theta(R_3)]$ (note that in $[\frac\pi 4,\frac\pi 2-\theta(R_3)]$, $G''\cap L_{ij}^{\frac\pi 4,\frac\pi 2-\theta(R_3)}$ coincides alreay with $Z$). The same argument in Proposition \tb{6.4} gives the conclusion of Proposition \tb{6.5}. \qed

\begin{pro}Let $G$ be an $(\eta,\d,\nu, L)$-sliding boundary for $Z$ for some $\eta<\eta_1$, $\d<R_1$, $\nu\in (0,R_3)$, and $L>0$. Suppose that $\d<\min\{\nu,\frac{R_3-\nu}{6}\}$, $L<\frac 19$, and $1^\circ$-$3^\circ$ in Proposition \tb{6.5} holds for when we replace $G_0$ by $G$ and $(1-\sqrt\eta)^2L/\d^2$ by $L$. Then
\be \sum_{1\le i\le 3}\sum_{1\le j\le 3}|q_{ij*}(\s^G_{ij})|_2=\sum_{1\le i\le 3}\sum_{1\le j\le 3}|q_{ij*}(z_{ij})|_2\ee
\end{pro}

\nd Take $i,j=1$ for example again. We already know that 
\be \s^G_{11}\cap (l_{kl}^{\theta(R_3),\frac\pi 2-\theta(R_3)}\times B_{Q_{kl}}(0,\d))=l_{kl}^{\theta(R_3),\frac\pi 2-\theta(R_3)},\ee
hence by \tb{(6.7)}, it is enough to study the structure of $\s^G_{11}$ in 
\be \begin{split}
B(z_{11}, \d)\bs [\cup_{k,l=2.3} &l_{kl}^{\theta(R_3),\frac\pi 2-\theta(R_3)}\times B_{Q_{kl}}(0,\d)]\\
=\cup_{k,l=2.3}&[B(l_{kl},\d)\bs (l_{kl}^{\theta(R_3),\frac\pi 2-\theta(R_3)}\times B_{Q_{kl}}(0,\d))]\\
=\cup_{k,l=2,3}&\{B([l_{kl}^{0,\theta(R_3-3\d)}\cup l_{kl}^{\frac\pi 2-\theta(R_3-3\d),\frac\pi 2}],\d)\\
&\cup [[l_{kl}^{\theta(R_3-3\d),\theta(R_3)}\cup l_{kl}^{\frac\pi 2-\theta(R_3),\frac\pi 2-\theta(R_3-3\d)}] \times B_{Q_{kl}}(0,\d)]\}
\end{split}\ee

Let us first look at the situation in the second component $[l_{kl}^{\theta(R_3-3\d),\theta(R_3)}\cup l_{kl}^{\frac\pi 2-\theta(R_3),\frac\pi 2-\theta(R_3-3\d)}] \times B_{Q_{kl}}(0,\d)$ of the last union.

Now by $1^\circ$ and $2^\circ$, for each $k,l=2,3$, there exists an $L$-Lipschitz map $\xi^G_{kl}:l^\nu_{kl}\to B_{Q_{kl}}(0,R)$, so that 
\be G\cap L^\nu_{kl}=G(\xi^G_{kl}, \theta(\nu), \frac\pi2-\theta(\nu)),\ee
\be\xi^G_{kl}(l_{kl}^{\theta(R_3-3\d), \frac\pi 4})\subset [-\d y_l,\d y_l],\ \xi^G_{kl}(l_{kl}^{\frac\pi 4,\frac\pi 2-\theta(R_3-3\d)})\subset [-\d x_k,\d x_k]\ee
and
\be \xi^G_{kl}(z)=0\mbox{ for any }z\in l_{kl}^{\theta(R_3), \frac\pi2-\theta(R_3)}.\ee 

As a result, we have

\be \pi_{x_1}(z,\xi^G_{kl}(z))=\pi_{x_1}(z)=\frac{\sqrt 3}{2}(1-\sqrt\eta), \forall z\in l_{kl}^{\theta(R_3-3\d), \theta(R_3)},\ee  and 
\be \pi_{y_1}(z, \xi^G_{kl}(z))=\pi_{y_1}(z)=\frac{\sqrt 3}{2}(1-\sqrt\eta), \forall z\in l_{kl}^{ \frac\pi 2-\theta(R_3), \frac\pi 2-\theta(R_3-3\d)},\ee 
and for any $z, w\in l_{kl}^{\theta(R_3-3\d), \theta(R_3)}$, suppose that $\arg z=\theta(R_3-t), t\in [0,3\d]$, $\arg w=\theta(R_3-s), s\in [0,3\d]$, then
\be\begin{split} \pi_{y_1}(z, \xi^G_{kl}(z))-\pi_{y_1}(w, \xi^G_{kl}(w))
&=[\pi_{y_1}(z)-\pi_{y_1}(w)]+[\pi_{y_1}(\xi^G_{kl}(z)-\xi^G_{kl}(w)]\\
&=\frac{\sqrt 3}{2}(t-s)+[\pi_{y_1}\circ\xi^G_{kl}(z)-\pi_{y_1}\circ\xi^G_{kl}(w)]\\
&\in [\frac{\sqrt 3}{2}(t-s)-\frac L2||t-s||, -\frac{\sqrt 3}{2}t+\frac L2||t-s||]\\
&= (t-s)[\frac{\sqrt 3}{2}-\frac L2, \frac{\sqrt 3}{2}+\frac L2].
\end{split}\ee

Since $L$ is very small, we know that $\frac{\sqrt 3}{2}-\frac L2>\frac{\sqrt 3}{2}(1-\frac 19)$, and $\frac{\sqrt 3}{2}+\frac L2<\frac{\sqrt 3}{2}(1+\frac 19)$, hence
\be (-1)^l[\frac{\sqrt 3}{2}(R_3-\frac83\d), \frac{\sqrt 3}{2}R_3]
\subset \{\pi_{y_1}(z, \xi^G_{kl}(z)):z\in l_{kl}^{\theta(R_3-3\d), \theta(R_3)}\}
\subset (-1)^l[\frac{\sqrt 3}{2}( R_3-\frac{10}{3}\d),\frac{\sqrt 3}{2}R_3],\ee

That is, 
\be\begin{split}& (-1)^k\frac{\sqrt 3}{2}(1-\sqrt\eta)x_1+(-1)^l[\frac{\sqrt 3}{2}(R_3-\frac83\d), \frac{\sqrt 3}{2}R_3]y_1\\
\subset &\pi_{y_1}(\s^G_{11}\cap [l_{kl}^{\theta(R_3-3\d), \theta(R_3)}\times B_{Q_{kl}}(0,\d)])\\
\subset &(-1)^k\frac{\sqrt3}{2}(1-\sqrt\eta)x_1+(-1)^l[\frac{\sqrt 3}{2}(R_3-\frac{10}{3}\d), \frac{\sqrt 3}{2}R_3]y_1.
\end{split}\ee

Similarly we have
\be\begin{split}& (-1)^l\frac{\sqrt 3}{2}(1-\sqrt\eta)y_1+(-1)^k[\frac{\sqrt 3}{2}(R_3-\frac83\d), \frac{\sqrt 3}{2}R_3]x_1\\
\subset &\pi_{y_1}(\s^G_{11}\cap [l_{kl}^{ \frac\pi 2-\theta(R_3),\frac\pi 2-\theta(R_3-3\d)}\times B_{Q_{kl}}(0,\d)])\\
\subset &(-1)^l\frac{\sqrt3}{2}(1-\sqrt\eta)y_1+(-1)^k[\frac{\sqrt 3}{2}(R_3-\frac{10}{3}\d), \frac{\sqrt 3}{2}R_3]x_1.
\end{split}\ee

Next let us look at the part of $\s^G_{11}$ in the first component $B([l_{kl}^{0,\theta(R_3-3\d)}\cup l_{kl}^{\frac\pi 2-\theta(R_3-3\d),\frac\pi 2}] ,\d)$ of the last union in \tb{(6.77)}. For this part, we simply have, for $z=a_za_k+b_zb_l+x_zx_k+y_zy_l\in \R^4$, 
\be q_{11}(z)=[(-1)^k\frac{\sqrt 3}{2}a_z-\frac12x_z]x_1+[(-1)^l\frac{\sqrt 3}{2}b_z-\frac 12 y_z]y_1,\ee
and hence for any $z\in B(l_{kl}^{0,\theta(R_3-3\d)},\d)$, since
\be z-(1-\sqrt\eta)c_{k2}\in (-\d,(R_3-3\d)b_l)\times B_{Q_{kl}}(0,\d),\ee
\be \begin{split}||q_{11}(z)-q_{11}((1-\sqrt\eta)c_{k2})||=||q_{11}(z-(1-\sqrt\eta)c_{k2})||\\
\in [-\frac{\sqrt 3}{2}(R_3-3\d)-\frac\d 2, \frac{\sqrt 3}{2}(R_3-3\d)+\frac\d 2]y_1\times [-\frac\d 2, \frac\d 2]x_1,
\end{split}\ee
hence
\be q_{11}(z)\in (-1)^k\frac{\sqrt 3}{2}(1-\sqrt\eta)x_1+[-\frac{\sqrt 3}{2}(R_3-3\d)-\frac\d 2, \frac{\sqrt 3}{2}(R_3-3\d)+\frac\d 2]y_1\times [-\frac\d 2, \frac\d 2]x_1.\ee
That is, 
\be \begin{split}&q_{11}[B(l_{kl}^{0,\theta(R_3-3\d)} ,\d)]\\
\subset& (-1)^k\frac{\sqrt 3}{2}(1-\sqrt\eta)x_1
+[-\frac{\sqrt 3}{2}(R_3-3\d)-\frac\d 2, \frac{\sqrt 3}{2}(R_3-3\d)+\frac\d 2]y_1\times [-\frac\d 2, \frac\d 2]x_1\\
\subset& (-1)^k\frac{\sqrt 3}{2}(1-\sqrt\eta)x_1
+[-\frac{\sqrt 3}{2}R_3, \frac{\sqrt 3}{2}R_3]y_1\times [-\frac\d 2, \frac\d 2]x_1.\end{split}\ee
Similarly we know that
\be \begin{split}&q_{11}[B(l_{kl}^{\frac\pi2-\theta(R_3-3\d),\frac\pi2},\d)]\\
\subset &(-1)^l\frac{\sqrt 3}{2}(1-\sqrt\eta)y_1
+[-\frac{\sqrt 3}{2}(R_3-3\d)-\frac\d 2, \frac{\sqrt 3}{2}(R_3-3\d)+\frac\d 2]x_1\times [-\frac\d 2, \frac\d 2]y_1\\
\subset& (-1)^l\frac{\sqrt 3}{2}(1-\sqrt\eta)y_1
+[-\frac{\sqrt 3}{2}R_3, \frac{\sqrt 3}{2}R_3]x_1\times [-\frac\d 2, \frac\d 2]y_1.\end{split}\ee

 
 To summerize, by \tb{(6.76), (6.85), (6.86), (6.91) and (6.92)} we have
 \be \begin{split}q_{11*}(\s^G_{11})&\subset q_{11}(z_{11})\\
 &\cup\{\cup_{k=2,3}(-1)^k\frac{\sqrt 3}{2}(1-\sqrt\eta)x_1
+[-\frac{\sqrt 3}{2}R_3, \frac{\sqrt 3}{2}R_3]y_1\times [-\frac\d 2, \frac\d 2]x_1\}\\
&\cup \{\cup_{l=2,3} (-1)^l\frac{\sqrt 3}{2}(1-\sqrt\eta)y_1
+[-\frac{\sqrt 3}{2}R_3, \frac{\sqrt 3}{2}R_3]x_1\times [-\frac\d 2, \frac\d 2]y_1\},
\end{split}\ee
and 
\be \begin{split}&q_{11*}(\s^G_{11}\cap \mI_k)\\
&= q_{11*}(\s^G_{11})\cap \{(-1)^k\frac{\sqrt 3}{2}(1-\sqrt\eta)x_1+[-\frac{\sqrt 3}{2}R_3, \frac{\sqrt 3}{2}R_3]y_1\times [-\frac\d 2, \frac\d 2]x_1\}\mbox{ for }k=1,2,\end{split}\ee
\be \begin{split}&q_{11*}(\s^G_{11}\cap \mK_l])\\
&=  q_{11*}(\s^G_{11})\cap \{(-1)^l\frac{\sqrt 3}{2}(1-\sqrt\eta)y_1
+[-\frac{\sqrt 3}{2}R_3, \frac{\sqrt 3}{2}R_3]x_1\times [-\frac\d 2, \frac\d 2]y_1\}\mbox{ for }l=1,2,\end{split}\ee

where \be \mI_k=\cup_{l=2,3}[l_{kl}^{0,\theta(R_3)}\times B_{Q_{kl}}(0,\d)\cup B(l_{kl}^{0,\theta(R_3-3\d)},\d)], k=2,3,\ee and
\be \mK_l=\cup_{k=2,3}[l_{kl}^{\frac\pi2-\theta(R_3),\frac \pi2}\times B_{Q_{kl}}(0,\d)\cup B(l_{kl}^{\frac\pi2-\theta(R_3-3\d),\frac\pi 2},\d)], l=2,3.\ee
Also, for $k=2,3$, $l=2,3$, set
\be\Xi_{11}^k=(-1)^k\frac{\sqrt 3}{2}(1-\sqrt\eta)x_1+[-\frac{\sqrt 3}{2}R_3, \frac{\sqrt 3}{2}R_3]y_1\times [-\frac\d 2, \frac\d 2]x_1=q_{11}([o_{I_{k2}}, o_{I_{k3}}])\times [-\frac\d 2, \frac\d 2]x_1,\ee
\be \Xi_{11,l}=(-1)^l\frac{\sqrt 3}{2}(1-\sqrt\eta)y_1
+[-\frac{\sqrt 3}{2}R_3, \frac{\sqrt 3}{2}R_3]x_1\times [-\frac\d 2, \frac\d 2]y_1=q_{11}([o_{J_{2l}}, o_{J_{3l}}])\times [-\frac\d 2, \frac\d 2]y_1.\ee

Then we have
\be q_{11*}(\s^G_{11})\subset q_{11}(z_{11})\cup[\cup_{k=2,3}\Xi_{11}^k]\cup[\cup_{l=2,3}\Xi_{11,l}],\ee

\be q_{11*}(\s^G_{11}\cap\mI_k)= q_{11*}(\s^G_{11})\cap  \Xi_{11}^k, q_{11*}(\s^G_{11}\cap\mK_l)= q_{11*}(\s^G_{11})\cap\Xi_{11,l}. \ee

See the following picture. 

\centerline{\includegraphics[width=0.5\textwidth]{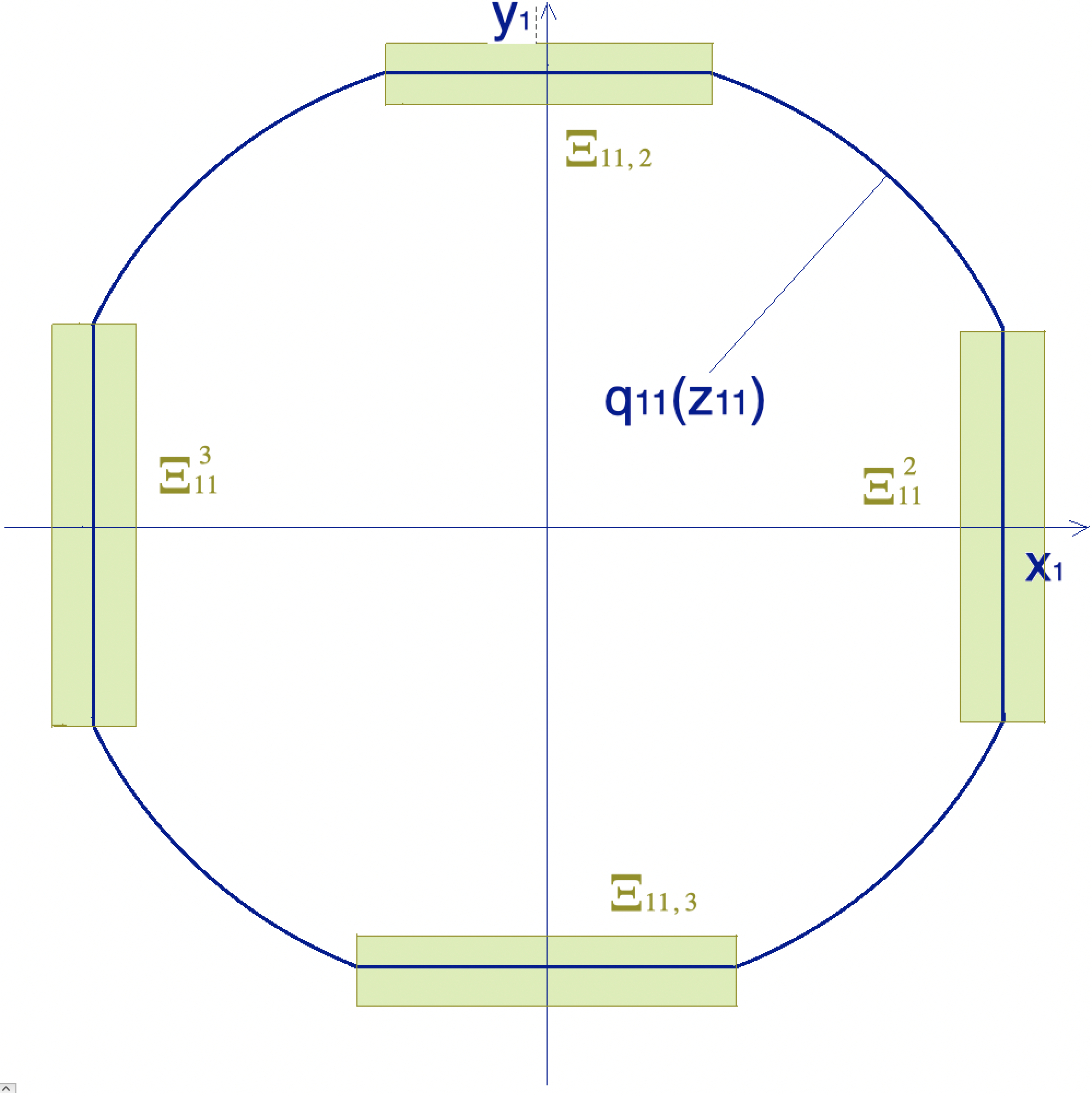}}

%

Now let us work in the plane $Q_{11}$ with coordinate under the basis $\{x_1,y_1\}$. Set $\frac{\sqrt 3}{2}(1-\sqrt\eta)=R_4$, $\frac{\sqrt 3}{2}R_3=R_5$. Set $\cG_{11}=q_{11*}(Z_{11})\bs [(\cup_{k=2,3}\Xi_{11}^k)\cup(\cup_{l=2,3}\Xi_{11,l})]$.

Since $\s^G_{11}$ is closed, so does $q_{11*}(\s^G_{11})$. By \tb{(6.100)}, it is easy to see that
\be \S^G_{11}=\cG_{11}+(\sum_{k=2,3}\S^G_{11}\cap \Xi_{11}^k)+(\sum_{l=2,3}\S^G_{11}\cap \Xi_{11,l}),\ee
 and hence
 \be |q_{11*}(\s^G_{11})|_2=\H^2(\cG_{11})+(\sum_{k=2,3}|\S^G_{11}\cap \Xi_{11}^k|)+(\sum_{l=2,3}|\S^G_{11}\cap \Xi_{11,l}|).\ee

Now we introduce the following notation: for finitely many points $z_1,z_2,\cdots, z_m$, let $[z_1,z_2,\cdots z_m]$ denote the piecewise linear 1-chain $\sum_{i=1}^{m-1}[z_i,z_{i+1}]$. 

Note that $\pa(\S^G_{11}\cap \Xi_{11}^k)$ is the sum of $q_{11*}(\s^G_{11})\cap \Xi_{11}^k$ with $[((-1)^kR_4, -R_5), ((-1)^k(R_4-\frac\d 2), -R_5), ((-1)^k(R_4-\frac\d 2), R_5), ((-1)^kR_4, R_5)]$. On the other hand, we know that
$ [((-1)^kR_4, -R_5), ((-1)^k(R_4-\frac\d 2), -R_5), ((-1)^k(R_4-\frac\d 2), R_5), ((-1)^kR_4, R_5)]$ is the projection under $q_{11}$ of 
\be [o_{I{k2}}, o_{I{k2}}-(-1)^k\d x_k, o_{I{k3}}-(-1)^k\d x_k, o_{I{k3}}].\ee
By \tb{(6.101)}, we have, for $k=2,3$, 
\be \begin{split}&\pa(\S^G_{11}\cap \Xi_{11}^k)=[q_{11*}(\s^G_{11})\cap \Xi_{11}^k]\\
&+ [((-1)^kR_4, -R_5), ((-1)^k(R_4-\frac\d 2), -R_5), ((-1)^k(R_4-\frac\d 2), R_5), ((-1)^kR_4, R_5)]\\
=&q_{11*}(\s^G_{11}\cap \mI_{k}+[o_{I{k2}}, o_{I{k2}}-(-1)^k\d x_k, o_{I{k3}}-(-1)^k\d x_k, o_{I{k3}}]).
\end{split}\ee

Similarly we have
\be \pa(\S^G_{11}\cap \Xi_{11,l})=q_{11*}(\s^G_{11}\cap \mK_{k\l}+[o_{J{2l}}, o_{J{2l}}-(-1)^k\d y_l, o_{J{2l}}-(-1)^k\d y_l, o_{J{3l}}]).\ee

As a result, combine with \tb{(6.103)}, we have
\be \begin{split} &|q_{11*}(\s^G_{11})|_2=\H^2(\cG_{11})\\
&+(\sum_{k=2,3}|q_{11*}(\s^G_{11}\cap \mI_{k}+[o_{I{k2}}, o_{I{k2}}-(-1)^k\d x_k, o_{I{k3}}-(-1)^k\d x_k, o_{I{k3}}])|_2)\\
&+(\sum_{l=2,3}|q_{11*}(\s^G_{11}\cap \mK_{l}+[o_{J{2l}}, o_{J{2l}}-(-1)^k\d y_l, o_{J{3l}}-(-1)^k\d y_l, o_{J{3l}}])|_2).
\end{split}\ee

The above discuss is for $i=j=1$. Now we define similarly, for $1\le i,j\le 3$, for $k\ne i$, 
\be \mI^j_k=\cup_{l\ne j}[l_{kl}^{0,\theta(R_3)}\times B_{Q_{kl}}(0,\d)\cup B(l_{kl}^{0,\theta(R_3-3\d)},\d)], \ee and set, for $l\ne j$,
\be \mK^i_l=\cup_{k\ne i}[l_{kl}^{\frac\pi2-\theta(R_3),\frac \pi2}\times B_{Q_{kl}}(0,\d)\cup B(l_{kl}^{\frac\pi2-\theta(R_3-3\d),\frac\pi 2},\d)].\ee

Also set, for $k\ne i$, 
\be\Xi_{ij}^k=\cup_{l\ne j}q_{ij}([o_{I_{kl}}, c_{k2}])\times [-\frac\d 2, \frac\d 2]x_i,\ee
\be \Xi_{ij,l}=\cup_{k\ne i}q_{ij}([o_{J_{kl}}, c_{1l}])\times [-\frac\d 2, \frac\d 2]y_i.\ee
Set $\cG_{ij}=q_{ij*}(Z_{ij})\bs [(\cup_{k\ne i}\Xi_{ij}^k)\cup(\cup_{l\ne j}\Xi_{ij,l})]$. Then as above, the exact same argument gives
\be \begin{split} &|q_{ij*}(\s^G_{ij})|_2=\H^2(\cG_{ij})\\
&+(\sum_{k\ne i}|q_{ij*}(\s^G_{ij}\cap \mI^j_{k}+[o_{I{k\lg j+1\rg }}, o_{I{k\lg j+1\rg }}+(-1)^{\lg k-i\rg }\d x_k, o_{I{k\lg j+2\rg }}+(-1)^{\lg k-i\rg }\d x_k, o_{I{k\lg j+2\rg }}])|_2)\\
&+(\sum_{l\ne j}|q_{11*}(\s^G_{ij}\cap \mK^i_{l}+[o_{J{\lg i+1\rg l}}, o_{J{\lg i+1\rg l}}+(-1)^{\lg l-j\rg }\d y_l, o_{J{\lg i+2\rg l}}+(-1)^{\lg l-j\rg }\d y_l, o_{J{\lg i+2\rg l}}])|_2),
\end{split}\ee
where $\lg m\rg \in 1,2,3$ is congruent to $m$ modulo 3.

We sum over $1\le i,j\le 3$, and get
\be \begin{split}&\sum_{1\le i,j\le 3}|q_{ij*}(\s^G_{ij})|_2=\sum_{1\le i,j\le 3}\{\H^2(\cG_{ij})\\
&+(\sum_{k\ne i}|q_{ij*}(\s^G_{ij}\cap \mI^j_{k}+[o_{I{k\lg j+1\rg }}, o_{I{k\lg j+1\rg }}+(-1)^{\lg k-i\rg }\d x_k, o_{I{k\lg j+2\rg }}+(-1)^{\lg k-i\rg }\d x_k, o_{I{k\lg j+2\rg }}])|_2)\\
&+(\sum_{l\ne j}|q_{11*}(\s^G_{ij}\cap \mK^i_{l}+[o_{J{\lg i+1\rg l}}, o_{J{\lg i+1\rg l}}+(-1)^{\lg l-j\rg }\d y_l, o_{J{\lg i+2\rg l}}+(-1)^{\lg l-j\rg }\d y_l, o_{J{\lg i+2\rg l}}])|_2)\}\\
=&[\sum_{1\le i,j\le 3}\H^2(\cG_{ij})]+\\
&+[\sum_{1\le k,j\le 3}(\sum_{i\ne k}|q_{ij*}(\s^G_{ij}\cap \mI^j_{k}+[o_{I{k\lg j+1\rg }}, o_{I{k\lg j+1\rg }}+(-1)^{\lg k-i\rg }\d x_k, o_{I{k\lg j+2\rg }}+(-1)^{\lg k-i\rg }\d x_k, o_{I{k\lg j+2\rg }}])|_2)]\\
&+[\sum_{1\le i,l\le 3}(\sum_{j\ne l}|q_{ij*}(\s^G_{ij}\cap \mK^i_{l}+[o_{J{\lg i+1\rg l}}, o_{J{\lg i+1\rg l}}+(-1)^{\lg l-j\rg }\d y_l, o_{J{\lg i+2\rg l}}+(-1)^{\lg l-j\rg }\d y_l, o_{J{\lg i+2\rg l}}])|_2)].
\end{split}\ee

Let us now fix any pair $k, j$, and look at what is the term 
\be (\sum_{i\ne k}|q_{ij*}(\s^G_{ij}\cap \mI^j_{k}+[o_{I{k\lg j+1\rg }}, o_{I{k\lg j+1\rg }}+(-1)^{\lg k-i\rg }\d x_k, o_{I{k\lg j+2\rg }}+(-1)^{\lg k-i\rg }\d x_k, o_{I{k\lg j+2\rg }}])|_2).\ee
Again take $k=1,j=1$ for example, we have to calculate
\be \begin{split}&|q_{21*}(\s^G_{21}\cap \mI^1_1+[o_{I{12}}, o_{I{12}}-\d x_1, o_{I{13}}-\d x_1, o_{I{13}}])|_2\\
+&
|q_{31*}(\s^G_{31}\cap \mI^1_1+[o_{I{12}}, o_{I{12}}+\d x_1, o_{I{13}}+\d x_1, o_{I{13}}])|_2.
\end{split}\ee

By definition, $\s_{i1^G}\cap \mI^1_1=\varphi_{G*}(z_{i1})\cap \mI^1_1,i=2,3$. Note that $z_{21}=z_{31}=l_{12}\cup l_{13}$ in $B([l_{12}^{0,\frac\pi2-\theta(\d)}\cup l_{13}^{0,\frac\pi2-\theta(\d)}],\d)$, $\mI_1^1=B([l_{12}^{0,\frac\pi2-\theta(\d)}\cup l_{13}^{\theta(R_3),\frac\pi2-\theta(\d)}],\d)   \bs [(l_{12}^{\theta(R_3),\frac\pi2-\theta(\d)}\times B_{12}(0,\d))\cup (l_{13}^{0,\frac\pi2-\theta(\d)}\times B_{13}(0,\d)) ]$, and $|\varphi_G(z)-z|<\d$ on $Z$, hence $\s_{21^G}\cap \mI^1_1=\s_{31^G}\cap \mI^1_1$. 

Also, by definition, we have 
\be\s_{21^G}\cap \mI^1_1\subset (1-\sqrt\eta)a_1+a_1^\perp.\ee

Now for any $z\in (1-\sqrt\eta)a_1+a_1^\perp$, write $z=(1-\sqrt\eta)a_1+x_zx_1+b_zb_1+y_zy_1$, then
\be q_{21}(z)=[-R_4-\frac 12 x_z]x_2+y_zy_1, q_{31}(z)=[R_4-\frac 12 x_z]x_3+y_zy_1.
\ee

Let $f: Q_{21}\to Q_{31}$ be the affine map so that $f(xx_2+yy_1)=(x+2R_4)x_3+yy_1$. Then $f$ is an isomorphism, and we have
$f\circ q_{21}(z)=q_{31}(z)$ for $z\in (1-\sqrt\eta)a_1+a_1^\perp$. In particular, we have 
\be f_*\circ q_{21*}(s^G_{21}\cap \mI^1_1)=q_{31*}(s^G_{21}\cap \mI^1_1)=q_{31*}(s^G_{31}\cap \mI^1_1).\ee

 On the other hand, we observe that
\be \begin{split}& q_{21}([o_{I{12}}, o_{I{12}}-\d x_1, o_{I{13}}-\d x_1, o_{I{13}}])=\\
&[-R_4x_2+R_5y_1, 
(-R_4+\frac\d 2)x_2+R_5y_1,\\
 &(-R_4+\frac\d 2)x_2-R_5y_1, -R_4x_2-R_5y_1],
\end{split}\ee
and
\be \begin{split}&q_{31}([o_{I{12}}, o_{I{12}}+\d x_1, o_{I{13}}+\d x_1, o_{I{13}}])=\\
&[R_4x_2+R_5y_1, 
(R_4-\frac\d 2)x_2+R_5y_1,\\
 &(R_4-\frac\d 2)x_2-R_5y_1, R_4 x_2-R_5y_1].
\end{split}\ee
Thus we have
\be \begin{split}f\circ q_{21}&([o_{I{12}}, o_{I{12}}-\d x_1, o_{I{13}}-\d x_1, o_{I{13}}])+q_{31}([o_{I{12}}, o_{I{12}}+\d x_1, o_{I{13}}+\d x_1, o_{I{13}}])\\
&\mbox{ is the square }\\
[(R_4+&\frac\d 2)x_2+R_5y_1, (R_4+\frac\d 2)x_2-R_5y_1,(R_4-\frac\d 2)x_2-R_5y_1, (R_4-\frac\d 2)x_2-R_5y_1].\end{split}\ee

Combine with \tb{(6.118)}, we know that the $\Z_2$-1-chain
\be\begin{split} f_*\circ &q_{21*}(\s^G_{21}\cap \mI^1_1+[o_{I{12}}, o_{I{12}}-\d x_1, o_{I{13}}-\d x_1, o_{I{13}}])\\
&+
q_{31*}(\s^G_{31}\cap \mI^1_1+[o_{I{12}}, o_{I{12}}+\d x_1, o_{I{13}}+\d x_1, o_{I{13}}])\\
=&[(R_4+\frac\d 2)x_2+R_5y_1, (R_4+\frac\d 2)x_2-R_5y_1,(R_4-\frac\d 2)x_2-R_5y_1, (R_4-\frac\d 2)x_2-R_5y_1].\end{split}\ee
As consequence, since the support of $f_*\circ q_{21*}(\s^G_{21}\cap \mI^1_1)=q_{31*}(\s^G_{31}\cap \mI^1_1)$ are contained in the square $[(R_4+\frac\d 2)x_2+R_5y_1, (R_4+\frac\d 2)x_2-R_5y_1,(R_4-\frac\d 2)x_2-R_5y_1, (R_4-\frac\d 2)x_2-R_5y_1]$, 
\be\begin{split} |f_*\circ &q_{21*}(\s^G_{21}\cap \mI^1_1+[o_{I{12}}, o_{I{12}}-\d x_1, o_{I{13}}-\d x_1, o_{I{13}}])|_2\\
&+
|q_{31*}(\s^G_{31}\cap \mI^1_1+[o_{I{12}}, o_{I{12}}+\d x_1, o_{I{13}}+\d x_1, o_{I{13}}])|_2\\
\ge &|[(R_4+\frac\d 2)x_2+R_5y_1, (R_4+\frac\d 2)x_2-R_5y_1,(R_4-\frac\d 2)x_2-R_5y_1, (R_4-\frac\d 2)x_2-R_5y_1]|_2.
\end{split}\ee
 
Now if we let $G=\pa Z$, (thus $\varphi^{\pa Z}=id$, then $s^{\pa Z}_{ij}=z_{ij}$, and $\s^{\pa Z}_{21}\cap \mI^1_1=\s^{\pa Z}_{31}\cap \mI^1_1=[o_{I{12}}, o_{I{13}}]$, hence it is easy to see that in this case equality holds in \tb{(6.123)}, this yields
\be\begin{split}  |f_*\circ &q_{21*}(\s^G_{21}\cap \mI^1_1+[o_{I{12}}, o_{I{12}}-\d x_1, o_{I{13}}-\d x_1, o_{I{13}}])|_2\\
&+
|q_{31*}(\s^G_{31}\cap \mI^1_1+[o_{I{12}}, o_{I{12}}+\d x_1, o_{I{13}}+\d x_1, o_{I{13}}])|_2\\
\ge |f_*\circ &q_{21*}(\s^{\pa Z}_{21}\cap \mI^1_1+[o_{I{12}}, o_{I{12}}-\d x_1, o_{I{13}}-\d x_1, o_{I{13}}])|_2\\
&+
|q_{31*}(\s^{\pa Z}_{31}\cap \mI^1_1+[o_{I{12}}, o_{I{12}}+\d x_1, o_{I{13}}+\d x_1, o_{I{13}}])|_2.
\end{split}\ee

Similarly we have, for any pair $k,j$, 
\be \begin{split}(\sum_{i\ne k}|q_{ij*}(\s^G_{ij}\cap \mI^j_{k}+[o_{I{k\lg j+1\rg }}, o_{I{k\lg j+1\rg }}+(-1)^{\lg k-i\rg }\d x_k, o_{I{k\lg j+2\rg }}+(-1)^{\lg k-i\rg }\d x_k, o_{I{k\lg j+2\rg }}])|_2)\\
\ge (\sum_{i\ne k}|q_{ij*}(\s^{\pa Z}_{ij}\cap \mI^j_{k}+[o_{I{k\lg j+1\rg }}, o_{I{k\lg j+1\rg }}+(-1)^{\lg k-i\rg }\d x_k, o_{I{k\lg j+2\rg }}+(-1)^{\lg k-i\rg }\d x_k, o_{I{k\lg j+2\rg }}])|_2).\end{split}\ee

Same argument gives that, for any pair $i,l$, 
\be \begin{split}\sum_{1\le i,l\le 3}(\sum_{j\ne l}|q_{ij*}(\s^G_{ij}\cap \mK^i_{l}+[o_{J{\lg i+1\rg l}}, o_{J{\lg i+1\rg l}}+(-1)^{\lg l-j\rg }\d y_l, o_{J{\lg i+2\rg l}}+(-1)^{\lg l-j\rg }\d y_l, o_{J{\lg i+2\rg l}}])|_2)\\
\ge \sum_{1\le i,l\le 3}(\sum_{j\ne l}|q_{ij*}(\s^{\pa Z}\cap \mK^i_{l}+[o_{J{\lg i+1\rg l}}, o_{J{\lg i+1\rg l}}+(-1)^{\lg l-j\rg }\d y_l, o_{J{\lg i+2\rg l}}+(-1)^{\lg l-j\rg }\d y_l, o_{J{\lg i+2\rg l}}])|_2).
\end{split}\ee
Combine with \tb{(6.113)}, we get
\be \sum_{1\le i,j\le 3}|q_{ij*}(\s^G_{ij})|_2\ge \sum_{1\le i,j\le 3}|q_{ij*}(z_{ij})|_2.\ee\qed

\noindent \textbf{Proof of Theorem \tb{6.1}}.

It is enough to combine the results of Propositions \tb{6.2, 6.4, 6.5 and 6.6}.\qed

 \begin{thm}The set $Y\times Y$ is $\Z_2$-topological sliding stable and Almgren sliding stable.
 \end{thm}
 
 \nd It is enough to apply Corollary \tb{4.4}, Proposition \tb{5.3}, Corollary \tb{5.5} and Theorem \tb{6.1}.\qed

\section{Uniqueness for $Y\times Y$}

In this section we deal with the uniqueness property for $Y\times Y$. Let us first introduce the notions of uniqueness:

\begin{defn}
Let $C$ be a $d$-dimensional reduced Almgren minimal set in a bounded domain $U$, we say that 

$1^\circ$ $C$ is Almgren unique in $U$, if $\H^d(C)=\inf_{F\in \overline \F(C,U)}\H^d(F)$, and
\be \forall \mbox{ reduced set }E\in \overline\F(C,U), \H^d(E)=\inf_{F\in \overline \F(C,U)}\H^d(F)\mbox{ implies } E=C.\ee

$2^\circ$ $C$ is $G$-topological unique in $U$, if $C$ is $G$-topological minimal, and 
\be \begin{split}\mbox{for any reduced }d\mbox{-dimensional }G-\mbox{topological competitor }E\mbox{ of }C\mbox{ in }U,\\
\H^d(E\cap U)=\H^d(C\cap U) \mbox{ implies }C=E;\end{split}\ee

$3^\circ$ We say that a $d$-dimensional minimal set $C$ in $\R^n$ is Almgren (resp. $G$-topological) unique, if it is Almgren (resp. $G$-topologial) unique in every bounded domain $U\subset \R^n$.
\end{defn}

For minimal cones, we have immediately:

\begin{pro}[Unique minimal cones, cf. \cite{uniquePYT}, Proposition 3.2]Let $K$ be a $d$-dimensional Almgren minimal cone in $\R^n$. Then it is Almgren (resp. $G$-topological) unique, if and only if it is Almgren (resp. $G$-topological) unique in some bounded convex domain $U$ that contains the origin. 
\end{pro}

\begin{pro}[cf. \cite{uniquePYT}, Proposition 3.4 and Proposition 3.5]Let $K\subset \R^n$ be a $G$-topological unique minimal cone of dimension $d$. Then it is also $G$-topological unique and Almgren unique of dimension $d$ in $\R^m$ for all $m\ge n$.
\end{pro}

After the above preliminaries, we are now going to prove the following :

\begin{thm}The 2-dimensional minimal cone $Y\times Y$ is $\Z_2$ topological unique and Almgren unique of dimension 2 in $\R^n$ for all $n\ge 4$.\end{thm}

\nd Let us first prove that $Y\times Y$ is $\Z_2$-topological unique in $\R^4$. By Proposition \tb{7.2}, it is enough to prove that $Y\times Y$ is $\Z_2$-topological unique in $\cU=\cU(Y\times Y, \eta)$ for some $\eta<\eta_1$. So fix any $\eta<\eta_1$.

Let $Z$ be the set $Y_1\times Y_2$ as defined at the beginning of Subsection \tb{3.2}, and take all the notations there. 

Suppose that $E$ is a reduced $\Z_2$-topological competitor for $Z$ in $\cU$, so that 
\be \H^2(E\cap \cU)=\H^2(Z\cap \cU).\ee

Then since $E$ is a $\Z_2$-topological competitor for $Z$ in $\cU$, all $\Z_2$-topological competitors $F$ for $E$ in $\cU$ are also $\Z_2$-topological competitors for $Z$ in $\cU$. Since $Z$ is $\Z_2$-topological minimal in $\cU$, 
\be\begin{split} \H^2(E\cap \cU)&=\H^2(Z\cap \cU)=\inf\{\H^2(F\cap \cU): F\mbox{ is a }\Z_2\mbox{-topological competitor for }Z\mbox{ in }\cU\}\\
&\le 
 \inf\{\H^2(F\cap \cU): F\mbox{ is a }\Z_2\mbox{-topological competitor for }E\mbox{ in }\cU\}\le \H^2(E\cap\cU).
 \end{split}\ee
 Hence $E$ is also $\Z_2$-topological minimal in $\cU$. By regularity for minimal sets, $E$ is 2-regular.
In particular, for almost all $x\in E$, the tangent plane $T_xE$ exists.

By definition, since $E$ is a 2-regular $\Z_2$-topological competitor for $Z$ in $\cU$, $F:=E\cap\bar\cU$ is atomatically in the class $\F(\eta,\d,\nu,L)$ for any $\d\in (0,R_1)$ and $\nu\in (0,R_3)$, and the corresponding $\s_{ij}$ defined in \tb{(5.12)} is $z_{ij}$. Hence by Proposition \tb{5.2}, there exists subsets $F_{ij}, 1\le i,j\le 3$, so that $1^\circ$ and $2^\circ$ of Proposition \tb{5.2} holds. Thus, by Proposition \tb{5.3} and Corollary \tb{5.5}, we know that
\be \H^2(F)\ge \sum_{1\le i,j\le 3}|q_{ij*}(z_{ij})|_2=\H^2(Z\cap \cU)=\H^2(E\cap \cU)=\H^2(F),\ee
hence we have
\be \H^2(F)= \sum_{1\le i,j\le 3}|q_{ij*}(z_{ij})|_2.\ee

Thus, the inequalities in the proof of Proposition \tb{5.3} and Lemma \tb{5.4} are all equalities. Therefore we get, in particular, 
that

$1^\circ$ If we set, for $1\le i,j\le 3$, $K_{ij}=\cap_{(k,l)\ne (i,j)} F_{kl}$, then we have
\be F=\cup_{1\le i,j\le 3}K_{ij}\mbox{, and the union is disjoint modulo }\H^2\mbox{-hull sets;}\ee

$2^\circ$ For $\H^2$-a.e. $x\in K_{ij}$, $T_xK_{ij}\perp \sum_{(k,l)\ne (i,j)}(-1)^{k+l-i-j}v_{kl}$. As a result, 
\be T_xF=T_x K_{ij}=Q_{a_i\wg b_j}=P_{ij}\mbox{ for }\H^2-a.e. x\in K_{ij}.\ee

Now since $F=E\cap \bar\cU$ is minimal in $\cU$, if $x\in F\cap\cU$ is a regular point of $F$, then by Theorem \tb{2.20}, there exists $r=r(x)>0$ such that $B(x,r)\subset \cU$, and in $B(x,r)$, $F$ is the graph of a $C^1$ function from $T_xF$ to $T_xF^\perp$, hence for all $y\in F\cap B(x,r)$, the tangent plane $T_yF$ exists, and the map $f:F\cap B(x,r)\to G(3,2): y\mapsto T_yF$ is continuous. But by \tb{(7.7) and (7.8)}, we have only nine choices (which are isolated points in $G(3,2)$) for $T_yF$, hence $f$ is constant, and $T_yF=T_xF$ for all $y\in F\cap B(x,r)$. As a result, 
\be F\cap B(x,r)=(T_xF+x)\cap B(x,r)\ee
is a disk parallel to one of the $P_{ij}$.

Still by the $C^1$ regularity Theorem \tb{2.20}, the set $F_P\cap \cU$ is a $C^1$ manifold, and is open in $F$. Thus, we deduce that
\be \begin{split}\mbox{each connected component of }F_P\cap \cU\mbox{ is part of a plane}\\
\mbox{ that is parallel to one of the }P_{ij}.\end{split}\ee

Let us look at $F_Y$. First, $F_Y\ne\emptyset$: otherwise, by \tb{Corollary 2.23 $2^\circ$}, $F\cap \cU=F_P\cap \cU$, and hence is a union of planes. But $F\cap \pa\cU$ does not coincide with any union of planes.

Take any $x\in F_Y$, then by the $C^1$ regularity around $\Y$ points (Theorem \tb{2.20} and Remark \tb{2.21}), there exists $r=r(x)>0$ such that $B(x,r)\subset \cU$, and in $B(x,r)$, $F$ is the image of a $C^1$ diffeomorphism $\varphi$ of a $\Y$-set Y, and $Y$ is tangent to $F$ at $x$. Denote by $L_Y$ the spine of $Y$, and by $R_i,1\le i\le 3$ the three open half planes of $Y$. Then $\varphi(R_i),1\le i\le 3$ are connected subsets $F_P$, hence each of them is a part of a plane parallel to one of the $P_{ij},1\le i,j\le 3$. As consequence, $\varphi(L_Y)\cap B(x,r)$ is an open segment passing through $x$ and parallel to one of the spines $D_{1i}, 1\le i\le 3$ and $D_{2j},1\le j\le 3$, where $D_{1i}$ is the line generated by $c_{1i}$, and $D_{2j}$ is the line generated by $c_{2j}$. Note that $D_{1i}$ is the intersection of the three $P_{ij}, 1\le j\le 3$, and $D_{2j}$ is the intersection of the three $P_{ij}, 1\le i\le 3$.

As a result, $F_Y\cap \cU$ is a union of open segments $I_1, I_2,\cdots$, each of which is parallel to one of the $D_{1i}, 1\le i\le 3$ and $D_{2j},1\le j\le 3$, and every endpoint is either a point on the boundary $\pa\cU$, or a point in $F_S\bs F_Y$ (singular points of type other than $\Y$). Moreover, 
\be \begin{split}\mbox{for each }&x\in F_Y\mbox{ such that }T_xF_Y=D_{1i}\mbox{(resp. }D_{2j}\mbox{), there exists }r>0\\
&\mbox{such that, in }B(x,r), F\mbox{ is a }\Y-\mbox{set whose spine is }x+D_{1i}\mbox{(resp. }D_{2j}).\end{split}\ee 

\begin{lem}If $x\in F_S\bs F_Y$, then $T_xF=Z$, and there exists $r>0$ so that $F\cap (x+r\cU)=(x+Z)\cap (x+r\cU)$.
\end{lem}

\nd Let $x\in F_S\bs F_Y$, and let $X$ be a blow up limit of $F$ at $x$. Then for each $n$, there exists $r_n>0$ so that $B(x,r_n)\subset \cU$, and $d_{x,r_n}(F,x+X)<\frac 1n$. By the bi-H\"older regularity Theorem \tb{2.17}, there exists a neighborhood $U_n\subset B(x,r_n)$ of $x$, and a bi-h\"older map : $f_n: B=B(0,1)\to U_n$, so that $F\cap U_n=f_n(X\cap B)$ and $f_n(0)=x$. By the structure theorem \tb{2.22} for  2-dimensional minimal cones, the set of $\Y$ points $X_Y$ of $X\cap B$ is a union of disjoint open segments $J_k,1\le k\le m$, so that $m$ is even, and for each $1\le k\le m$, $J_k=(0, x_k)$ where $x_k\in \pa B$.

As a result, the set of $\Y$ points $F_Y\cap U_n$ is the disjoint union of $f_n(J_k)$, with each $f_n(J_k)$ connected and admits $x$ as an endpoint. By \tb{(7.11)}, we know that each $f_n(J_k)$ is a segment contained in $x+D_{1i}$ for some $1\le i\le 3$, or in $x+D_{2j}$ for some $1\le j\le 3$. In otherwords, $F_Y\cap U_n\subset (x+Z)_Y\cap U_n$. In particular, the even number $m\le 6$. 

Since $F_Y\cap U_n\subset  (x+Z)_Y\cap U_n$, and $\frac {1}{r_n} d_{U_n}(F,x+X)<\frac 1n$, we know that $\sup\{d(p, (x+Z)_Y): p\in (x+X)Y\cap U_n\} <\frac 1nr_n$. Since both $Z$ and $X$ are cones, we know that $X_Y\subset Z_Y$.

With the same argument for regular points and  \tb{(7.10}), we get 
\be F\cap U_n\subset (x+Z)\cap U_n,\ee
 and hence
\be X\subset Z.\ee

Recall that $m$ is the number of $\Y$ points of $X\cap \pa B$: 

If $m=2$, then by the structure Theorem \tb{2.22}, $X$ is a $\Y$ set, this is impossible because we have supposed that $x$ is not of type $\Y$;

If $m=4$, again by Theorem \tb{2.22}, we know that $\Y$ is a set $T$ which is the cone over the 1-skeleton of a regular tetrahedron in $\R^3$. But the set $Z$ does not contain any such set;

Hence the only possibility is that $m=6$, and since $X\subset Z$, we must have $X=Z$.

As a result, since $F\cap U_1$ is a bi-H\"older image of $(x+Z)\cap B$, by \tb{(7.11)} we know that there exists $r>0$, with $x+r\cU\subset U_1$, $F\cap (x+r\cU)=(x+Z)\cap (x+r\cU)$.
\qed

For each $y\in \R^4$ and each $r>0$, set $\cU(y,r)=y+r\cU$. Then

After the above lemma, we are going to discuss two cases: when there exists at least a point in $F_S\bs F_Y$, or there is no such points.

\textbf{Case 1:} There exists a point $x\in F_S\bs F_Y$. 

\begin{lem} If there exists a point $x\in F_S\bs F_Y$, then $Z\cap \cU =F\cap \cU$.
\end{lem}

\nd Fix such a point $x$. By Lemma \tb{7.5}, there exists $r>0$ such that 
\be F\cap \cU(x,r)=(x+Z)\cap \cU(x,r).\ee

Recall that $F_Y\cap \cU$ is a union of open segments $I_1,\cdots, I_n,\cdots$, each of which is parallel to one of the $D_{1i}, 1\le i\le 3$ and $D_{2j},1\le j\le 3$, and every endpoint is either a point on the boundary $\pa\cU$, or a point in $F_S\bs F_Y$.

 We claim that : 
 \be\mbox{ for each segment }I_i\mbox{, at least one of its endpoints is in }\pa\cU.\ee
 
 In fact, suppose there is some $I_i$ so that neither endpoint lies in $\bar\cU$. As a result, both of them belong to $F_S\bs F_Y$. Let $p$ and $q$ denote its endpoints. They by Lemma \tb{7.5}, there exists $r_p>0$ and $r_q>0$ so that $F\cap \cU(p,r_p)=(p+Z)\cap \cU(p,r_p)$, $F\cap \cU(q,r_q)=(q+Z)\cap \cU(q,r_q)$. Hence $[p,q]\cap \cU(q,r_q)=(q+Z_Y)\cap \cU(q,r_q)$, and $[p,q]\cap \cU(p,r_p)=(p+Z_Y)\cap \cU(p,r_p)$. As a result, the half lines $R_{p,q}-p=R_{0,q-p}$ and $R_{q,p}-q=R_{0,p-q}$ both lie in the cone $Z_Y$. Hence the line generated by $p-q$ is part of $Z$. This is impossible because, according to the structure of $Z$, $Z_Y$ does not contain any line. 
 
 Thus we get Claim \tb{(7.15)}.

%
 
 Denote by $L_i,1\le i\le 6$, the six spines (which are half lines issued from $x$) of $Z+x$. Then $L_i\cap \cU\subset F_Y$. By \tb{(7.14)}, $L_i\cap \cU(x,r)$ is part of some $I_j\subset F_Y$. Hence $I_j$ already has an endpoint $x$ that does not belong to $\pa\cU$, therefore the other endpoint must lie in $\pa\cU$, which yields $I_j=L_i\cap \cU$.
 
 Now we take a one parameter family of regions $\cU_s=\cU(y_s, s)$, $r\le s\le 1$, with $\cU_r=\cU(x,r)$, $\cU_1=\cU$, such that
 
 $1^\circ$ $\cU_s\subsetneqq \cU_{s'}$ for all $s<s'$;
 
 $2^\circ$ $\cap_{1>t>s} \cU_t=\bar \cU_s$ and $\cup_{t<s}\cU_t=\cU_s$ for all $r\le s\le 1$.
 
 Set $R=\inf\{s>r, (Z+x)\cap \cU_s\ne F\}$. We claim that $R=1$.
 
 Suppose this is not true. By definition of $\cU_s$, we know that the six spines and the nine faces of $Z+x$ are never tangent to $\pa \cU_s$ for any $r<s<1$. Then we claim that 
 \be \pa \cU_R\cap F=\pa \cU_R\cap (Z+x)\subset F_P\cup F_Y.\ee 
 
 In fact, if $y$ belong to one of the $L_i$, then, $y\in L_i\cap \cU\bs\{x\}\subset E_Y$; otherwise, suppose $y$ does not lie in the six $L_i,1\le i\le 6$. Then $y$ belong to  $x+P_{ij}$ for some $1\le i,j\le 3$. As a result, for any $t>0$ small, we know that $E\cap B(y,t)\cap \cU_R=(x+P_{ij})\cap B(y,t)\cap \cU_R$. Note that the set $(x+P_{ij})\cap B(y,t)\cap \cU_R$ is almost a half disk when $t$ is sufficiently small, hence in particular, $F\cap B(y,t)$ cannot coincide with a $\Y$ set or a $Z$ set. After \tb{(7.11)} and Lemma \tb{7.5}, we must have $y\in F_P$.
 
 Thus we have Claim \tb{(7.16)}.
 
 Now take any $y\in \pa \cU_R\cap F$.
 
 If $y\in F_P$, then $y\in x+P_{ij}$ for some $i\ne j$. Thus $T_yF=P_{ij}$. By \tb{(7.9)}, and the fact that $R<1$, there exists $r_y>0$ such that $B(y,r_y)\subset \cU$ and $F\cap B(y,r_y)=(P_{ij}+y)\cap B(y,r_y)$. In other words, 
 \be\mbox{ there exists }r_y>0\mbox{ such that }F\mbox{ coincides with }Z+x\mbox{ in }B(y,r_y).\ee
 
 If $y$ is a $\Y$ point, then it lies in one of the $L_i$. By the same argument as above, using \tb{(7.11)}, we also have \tb{(7.17)}.
 
 Thus \tb{(7.17)} holds for all $y\in \pa \cU_R\cap F$. Since $\pa \cU_R\cap F$ is compact, we get an $r>0$, such that $F\cap B(\cU_R,r)=(Z+x)\cap B(\cU_R,r)$. By the continuous condition $2^\circ$ for the family $\cU_s$, there exists $R'\in (R,1)$ such that $\cU_{R'}\subset B(\cU_R,r)$. As consequence, $F\cap \cU_{R'}=(Z+x)\cap \cU_{R'}$, this contradicts the definition of $R$.
 
 Hence $R=1$, and by definition of $R$, we have $(Z+x)\cap \cU=F\cap \cU$. Since $F\cap \pa \cU=Z\cap \pa \cU$, and $F$ is closed and reduced, $x$ must be the origin. Thus we get the conclusion of Lemma \tb{7.6}.\qed
 
 \textbf{Case 2:} $F_S\bs F_Y=\emptyset$. In this case, the same kind of argument as in Lemma \tb{7.6} gives the following:
 
 \begin{lem}Let $x$ be a $\Y$ point in $F$. Then $F$ coincides with the intersection of $\cU$ with a $\Y$ set centered at the origin. 
 \end{lem}
 
 But this is impossible, because $E\cap \pa\cU=Z\cap \pa\cU$, which is not the intersection of a $\Y$ set with $\pa\cU$.
 
 Hence we have $E\cap \cU=Z\cap \bar \cU$, and thus $Z$ is topological unique in $\cU$. We thus get the $\Z_2$-topological uniqueness of $Y\times Y$ in $\R^4$.
 
By Proposition \tb{7.3}, $Y\times Y$ is also $\Z_2$ topological unique and Almgren unique of dimension 2 in $\R^n$ for all $n\ge 4$.\qed

\renewcommand\refname{References}
\bibliographystyle{plain}
\bibliography{reference}

\end{document}